\documentclass[11pt,twoside]{article}
\usepackage{theorem}
\usepackage{epsfig,amsmath,amssymb}
\usepackage{color}
\allowdisplaybreaks[1]

\newtheorem{thm}{\bf Theorem}[section]

\newtheorem{lem}[thm]{\bf Lemma}
\newtheorem{defn}{\bf Definition}[section]
\newtheorem{rem}{\bf Remark}[section]
\newtheorem{prop}[rem]{\bf Proposition}
\numberwithin{equation}{section}

\title{\LARGE\bf The conjecture of Ulam on the invariance of
                 measure on Hilbert cube}

\author{Soon-Mo Jung}

\date{\small\it Mathematics Section,
                College of Science and Technology, \\
                Hongik University, 30016 Sejong, Korea \\
                E-mail: {\tt smjung@hongik.ac.kr} \\
                \hspace*{12.5mm} {\tt smjung747@gmail.com}}

\begin{document}

\markboth{Ulam's Conjecture on Invariance of Measure}
         {S.-M. Jung}
\maketitle

\begin{abstract}
A conjecture of Ulam states that the standard probability
measure $\pi$ on the Hilbert cube $I^\omega$ is invariant
under the induced metric $d_a$ when the sequence
$a = \{ a_i \}$ of positive numbers satisfies the condition
$\sum\limits_{i=1}^\infty a_i^2 < \infty$.
This conjecture was proved in \cite{jung1} when $E_1$ is a
non-degenerate subset of $M_a$.
In this paper, we prove the conjecture of Ulam completely by
classifying cylinders as non-degenerate and degenerate
cylinders and by treating the degenerate case that was
overlooked in the previous paper.
\end{abstract}
\vspace{5mm}

\vspace{5mm}

\normalsize


\section{Introduction}

Throughout this paper, using the symbol $\mathbb{R}^\omega$, 
we represents an infinite dimensional real vector space defined 
as
\begin{eqnarray*}
\mathbb{R}^\omega = \big\{ (x_1, x_2, \ldots ) :
x_i \in \mathbb{R} ~\,\mbox{for all}\,~ i \in \mathbb{N} \,\big\}.
\end{eqnarray*}
From now on, we denote by $(\mathbb{R}^\omega, \mathcal{T})$
the product space $\prod\limits_{i=1}^\infty \mathbb{R}$,
where $(\mathbb{R}, \mathcal{T}_{\mathbb{R}})$ is the usual
topological space.
Since $(\mathbb{R}, \mathcal{T}_{\mathbb{R}})$ is a Hausdorff
space, so is $(\mathbb{R}^\omega, \mathcal{T})$
(see \cite[Theorem 111.7]{kasriel}).

Let $I = [0,1]$ be the unit closed interval,
$I^\omega = \prod\limits_{i=1}^\infty I$ the \emph{Hilbert cube},
and let $\pi$ be the \emph{standard product probability measure}
on $I^\omega$.
We denote by $(I^\omega, \mathcal{T}_\omega)$ the (topological)
subspace of $(\mathbb{R}^\omega, \mathcal{T})$.
Then, $\mathcal{T}_\omega$ is the relative topology for
$I^\omega$ induced by $\mathcal{T}$.

In this paper, let $a = \{ a_i \}_{i \in \mathbb{N}}$ be a
sequence of positive real numbers satisfying
\begin{eqnarray} \label{eq:20110922-1}
\sum_{i=1}^\infty a_i^2 < \infty.
\end{eqnarray}
Using this sequence $a = \{ a_i \}_{i \in \mathbb{N}}$, we
define the metric on $I^\omega$ by
\begin{eqnarray} \label{eq:20150928}
d_a (x,y)
= \Bigg( \sum_{i=1}^\infty a_i^2 (x_i - y_i)^2 \Bigg)^{1/2}
\end{eqnarray}
for all $x = (x_1, x_2, \ldots) \in I^\omega$ and
$y = (y_1, y_2, \ldots) \in I^\omega$.


\begin{rem} \label{rem:1.1}
It is to be noted that
\begin{itemize}
\item[$(i)$]   $d_a$ is consistent with the topology
               $\mathcal{T}_\omega$ and it is invariant under
               translation $($see {\rm \cite{5}}$)$;
\item[$(ii)$]  $(I^\omega, \mathcal{T}_\omega)$ is a Hausdorff
               space 
               as a subspace of the Hausdorff space
               $(\mathbb{R}^\omega, \mathcal{T})$;
\item[$(iii)$] $(I^\omega, \mathcal{T}_\omega)$ is a compact
               subspace of $(\mathbb{R}^\omega, \mathcal{T})$
               by Tychonoff's theorem.
\end{itemize}
\end{rem}
\vspace*{5mm}

S. M. Ulam raised the conjecture on the invariance of measures
defined in the compact metric space (see \cite {ulam}):
\begin{quote}
\emph{Let $E$ be a compact metric space.
Does there exist a finitely additive measure $m$ defined for
at least all the Borel subsets of $E$, such that $m(E) = 1$,
$m(p) = 0$ for all points $p$ of $E$, and such that congruent
sets have equal measure\/?}
\end{quote}
Thereafter, J. Mycielski \cite{5} confined the question of Ulam
to the Hilbert cube $I^\omega$ and reformulated it using modern
mathematical terms:
\begin{center}
\emph{The standard probability measure $\pi$ on $I^\omega$ is
$d_a$-invariant.}
\end{center}
The above statement is widely known today as the conjecture of
Ulam.

In 1974, by using the axiom of choice, J. Mycielski \cite{4,5}
answered the question of Ulam affirmatively under the additional
assumption that the sets are open.
In addition, he asked in \cite{4} whether one can prove the
conjecture of Ulam under the assumption that the sets are
closed.
J. W. Fickett \cite{1}, one step further, showed in 1982 that
Ulam's conjecture is true when the sequence
$a = \{ a_i \}_{i \in \mathbb{N}}$ decreases very rapidly to
$0$ such that
\begin{eqnarray*}
a_{i+1} = o \Big( a_i^{2^{i+1}} \Big).
\end{eqnarray*}

In 2018, the author and E. Kim proved in their paper
\cite{jungkim} that the Ulam's conjecture is true when the
sequence $a = \{ a_i \}_{i \in \mathbb{N}}$ of positive real
numbers is monotone decreasing and satisfies the condition
\begin{eqnarray*}
a_{i+1} = o \bigg( \frac{a_i}{\sqrt{i}} \bigg)
\end{eqnarray*}
(see also \cite{jung,jungkim2}).
It is evident that the last condition is much weaker than that
of Fickett.

Most recently, the author proved in \cite{jung1} that the
above conjecture is true when the involved sets are
non-degenerate.
For example, in the proof of \cite[Theorem 3.1]{jung1}, the
sentence ``Using the $\ldots\,$, Theorem 2.2 states that
$F : M_a \to M_a$ is a $d_a$-isometry $\ldots$ in a natural
way'' is true only if $E_1$ is a non-degenerate subset of
$M_a$.

In this paper, we completely prove that Ulam's conjecture is
true, taking into account the degenerate case that has been
overlooked in \cite{jung1}.
More precisely, under the assumption that the sequence
$a = \{ a_i \}_{i \in \mathbb{N}}$ of positive real numbers
satisfies the condition (\ref{eq:20110922-1}), we prove that
$\pi(E_1) = \pi(E_2)$ for any Borel subsets $E_1$ and $E_2$ of
$I^\omega$ which are $d_a$-isometric to each other, where
$\pi$ is the standard probability measure on $I^\omega$
(see also \cite{jung2}).

A notable improvement of this paper over the previous paper
\cite{jung1} is the introduction of the
\emph{generalized linear span} concept, by which we can easily
explain the intermediate results of the previous paper with
minimal changes.


\section{Preliminaries}
\label{sec:2}

We define
\begin{eqnarray*}
M_a = \Bigg\{ (x_1, x_2, \ldots) \in \mathbb{R}^\omega :\,
              \sum_{i=1}^\infty a_i^2 x_i^2 < \infty
      \Bigg\},
\end{eqnarray*}
where $a = \{ a_i \}_{i \in \mathbb{N}}$ is a sequence of
positive real numbers that satisfies the condition
(\ref{eq:20110922-1}).
Then $M_a$ is a vector space over $\mathbb{R}$, and we can
define an inner product $\langle \cdot, \cdot \rangle_a$ on
$M_a$ by
\begin{eqnarray*}
\langle x, y \rangle_a = \sum_{i=1}^\infty a_i^2 x_i y_i
\end{eqnarray*}
for all $x = (x_1, x_2, \ldots)$ and $y = (y_1, y_2, \ldots)$
of $M_a$.
This inner product induces the norm
\begin{eqnarray*}
\| x \|_a = \sqrt{\langle x, x \rangle_a}
\end{eqnarray*}
for all $x \in M_a$.


\begin{rem} \label{rem:new2.1}
$M_a$ is the set of all elements $x \in \mathbb{R}^\omega$
satisfying $\| x \|_a^2 < \infty$, \textit{i.e.},
\begin{eqnarray*}
M_a = \big\{ (x_1, x_2, \ldots) \in \mathbb{R}^\omega :\,
             \| x \|_a^2 < \infty
      \big\}.
\end{eqnarray*}
\end{rem}

In view of definition (\ref{eq:20150928}), the metric $d_a$ on
$I^\omega$ can be extended to the metric on $M_a$,
\textit{i.e.},
\begin{eqnarray*}
d_a(x,y) = \sqrt{\langle x-y, x-y \rangle_a}
\end{eqnarray*}
for all $x, y \in M_a$.

What follows is a basic definition we are familiar with, but
for the sake of completeness of the paper, we now define
precisely the $d_a$-isometry between subsets of $M_a$.


\begin{defn} \label{def:1.1}
Let $E_1$ and $E_2$ be nonempty subsets of $M_a$.
\begin{itemize}
\item[$(i)$]  A function $f : E_1 \to E_2$ is called a
              $d_a$-\emph{isometry} provided
              $d_a(f(x), f(y)) = d_a(x,y)$ for all
              $x, y \in E_1$;
\item[$(ii)$] $E_1$ is said to be $d_a$-\emph{isometric} to
              $E_2$ provided there exists a surjective
              $d_a$-isometry $f : E_1 \to E_2$.
\end{itemize}
\end{defn}

Assuming that the sequence $a = \{ a_i \}_{i \in \mathbb{N}}$
satisfies the condition (\ref{eq:20110922-1}), we can prove
that $(M_a, \langle \cdot, \cdot \rangle_a)$ is a Hilbert space
in the same way as \cite[Theorem 2.1]{jungkim}.
Let $(M_a, \mathcal{T}_a)$ be the topological space generated by
the metric $d_a$.
In view of Remark \ref{rem:1.1} $(ii), (iii)$ and using
\cite[Theorem 91.2]{kasriel}, it is easy to prove the following
remarks.
Since $[-M,M]^\omega \subset M_a$ for any fixed $M > 0$, we may
consider the families of open sets which are included in $M_a$
only to prove Remark \ref{rem:2.1} $(iii)$.
This idea, together with Remark \ref{rem:1.1} $(iii)$, implies
the validity of Remark \ref{rem:2.1} $(iii)$.


\begin{rem} \label{rem:2.1}
We note that
\begin{itemize}
\item[$(i)$]   $(M_a, \langle \cdot, \cdot \rangle_a)$ is a
               Hilbert space over $\mathbb{R}$;
\item[$(ii)$]  $(M_a, \mathcal{T}_a)$ is a Hausdorff space as
               a subspace of the Hausdorff space
               $(\mathbb{R}^\omega, \mathcal{T})$;
\item[$(iii)$] $(I^\omega, \mathcal{T}_\omega)$ is a compact
               subspace of $(M_a, \mathcal{T}_a)$;
\item[$(iv)$]  $(I^\omega, \mathcal{T}_\omega)$ is a closed subset
               of $(M_a, \mathcal{T}_a)$.
\end{itemize}
\end{rem}


\begin{defn} \label{def:2.3}
Given $c \in M_a$ the \emph{translation} by $c$ is the mapping
$T_c : M_a \to M_a$ defined by $T_c (x) = x + c$ for all
$x \in M_a$.
\end{defn}


\section{First-order generalized linear span}
\label{sec:3}

In \cite[Theorem 2.5]{jungkim}, we were able to extend the
domain of a $d_a$-isometry $f$ to the whole space when the
domain of $f$ is a \emph{non-degenerate basic cylinder}
(see Definition \ref{def:basic} for the exact definition of
non-degenerate basic cylinders).
However, we shall see in Definition \ref{def:3.5} and Theorem
\ref{thm:3.17} that the domain of a $d_a$-isometry $f$ can be
extended to the whole space whenever $f$ is defined on a
bounded set which contains more than one element.

From now on, it is assumed that $E$, $E_1$, and $E_2$ are
subsets of $M_a$, each of them contains more than one element,
unless specifically stated for their cardinalities, and that
they are bounded because the Hilbert cube $I^\omega$ is a
bounded subset of $M_a$ and the involved sets are indeed
(Borel) subsets of $I^\omega$ in the main theorems.

If the set has only one element or no element, this case will
not be covered here because the results derived from this case
are trivial and uninteresting.


\begin{defn} \label{def:gls}
Assume that $E$ is a nonempty bounded subset of $M_a$ and $p$
is a fixed element of $E$.
We define the \emph{first-order generalized linear span} of
$E$ with respect to $p$ as
\begin{eqnarray*}
\begin{split}
{\rm GS}(E,p)
& = \Bigg\{ p + \sum_{i=1}^m \sum_{j=1}^\infty
            \alpha_{ij} (x_{ij} - p) \in M_a :
            m \in \mathbb{N}\,;
            \vspace*{1mm}\\
& \hspace*{9mm}
            x_{ij} \in E ~\mbox{and}~
            \alpha_{ij} \in \mathbb{R}~
            \mbox{for~all~$i$~and~$j$}
    \Bigg\}.
\end{split}
\end{eqnarray*}
\end{defn}

We remark that if a bounded subset $E$ of $M_a$ contains more
than one element, then $E$ is a proper subset of its first-order
generalized linear span ${\rm GS}(E,p)$, because
$x = p + (x - p) \in {\rm GS}(E,p)$ for any $x \in E$ and
$p + \alpha (x - p) \in {\rm GS}(E,p)$ for any
$\alpha \in \mathbb{R}$, which implies that ${\rm GS}(E,p)$ is
unbounded.
Moreover, we note that $\alpha x + \beta y \in M_a$ for all
$x, y \in M_a$ and $\alpha, \beta \in \mathbb{R}$, because
$\| \alpha x + \beta y \|_a
 \leq | \alpha | \| x \|_a + | \beta | \| y \|_a < \infty$.
Therefore, ${\rm GS}(E,p) - p$ is a real vector space, because
the double sum in the definition of ${\rm GS}(E,p)$ guarantees
$\alpha x + \beta y \in {\rm GS}(E,p) - p$ for all
$x, y \in {\rm GS}(E,p) - p$ and $\alpha, \beta \in \mathbb{R}$
and because ${\rm GS}(E,p) - p$ is a subset of a real vector
space $M_a$ (\textit{cf.} Lemma \ref{lem:3.8} $(i)$ below).

We remark that the smallest flat containing $E$ was introduced
in \cite{1} for any subset $E$ of an $n$-dimensional Euclidean
space somewhat similarly to the first-order generalized linear
span as follows:
\begin{eqnarray*}
\begin{split}
H(E) & = \Bigg\{ \sum_{i=0}^n \alpha_i x_i :
                 x_i \in E ~\mbox{and}~ \alpha_i \in \mathbb{R}~
                 \mbox{for~all~$i$~with}~
                 \sum_{i=0}^n \alpha_i = 1
         \Bigg\} \\
     & = \Bigg\{ p + \sum_{i=1}^n \alpha_i (x_i - p) :
                 p, x_i \in E ~\mbox{and}~
                 \alpha_i \in \mathbb{R}~ \mbox{for~all~$i$}
         \Bigg\}.
\end{split}
\end{eqnarray*}

Given a $p = (p_1, p_2, \ldots) \in M_a$ and an
$n \in \mathbb{N}$, if we set
$\mathbb{R}_p^n = \{ (x_1, x_2, \ldots) \in M_a :
 x_i \in \mathbb{R} ~\mbox{for}~ 1 \leq i \leq n ~\mbox{and}~
 x_i = p_i ~\mbox{for}~ i > n \}$, then it is easy to see that
$H(E) \subset {\rm GS}(E,p)$ for any set
$E \subset \mathbb{R}_p^n$.
However, it is obvious that ${\rm GS}(E,p) \not\subset H(E)$
for some $E \subset \mathbb{R}_p^n$.

For each $i \in \mathbb{N}$, we set
$e_i = (0, \ldots, 0, 1, 0, \ldots)$, where $1$ is in the
$i$th position.
Then $\{ \frac{1}{a_i} e_i \}_{i \in \mathbb{N}}$ is a complete
orthonormal sequence in $M_a$.
The following definition introduces the concept of index based
on the `standard' coordinate system
$\{ \frac{1}{a_i} e_i \}_{i \in \mathbb{N}}$.


\begin{defn} \label{def:Lambda}
Let $E$ be a nonempty subset of $M_a$.
\begin{itemize}
\item[$(i)$]  We define the \emph{index set} of $E$ by
\begin{eqnarray*}
\Lambda(E) = \Big\{ i \in \mathbb{N} : \mbox{there are an}~
                    x \in E ~\mbox{and an}~
                    \alpha \in \mathbb{R} \!\setminus\! \{ 0 \}
                    ~\mbox{satisfying}~ x + \alpha e_i \in E
             \Big\}.
\end{eqnarray*}
Each $i \in \Lambda(E)$ is called an \emph{index} of $E$.
If $\Lambda(E) \neq \mathbb{N}$, then the set $E$ is called
\emph{degenerate}.
Otherwise, $E$ is called \emph{non-degenerate}.

\item[$(ii)$] Let $\beta = \{ \beta_i \}_{i \in \mathbb{N}}$
              be another complete orthonormal sequence in $M_a$.
              We define the \emph{$\beta$-index set} of $E$ by
\begin{eqnarray*}
\Lambda_\beta(E)
= \Big\{ i \in \mathbb{N} : \mbox{there are an}~ x \in E
         ~\mbox{and an}~
         \alpha \in \mathbb{R} \!\setminus\! \{ 0 \}
         ~\mbox{satisfying}~ x + \alpha \beta_i \in E
  \Big\}.
\end{eqnarray*}
Each $i \in \Lambda_\beta(E)$ is called a
\emph{$\beta$-index} of $E$.
\end{itemize}
\end{defn}

We will find that the concept of index set in Hilbert space
sometimes takes over the role that the concept of dimension
plays in vector space.
According to the definition above, if $i$ is an index of $E$,
\textit{i.e.}, $i \in \Lambda(E)$, then there are $x \in E$
and $x + \alpha e_i \in E$ for some $\alpha \neq 0$.
Since $x \neq x + \alpha e_i$, we remark that if
$\Lambda(E) \neq \emptyset$, then the set $E$ contains at least
two elements.

In the following lemma, we prove that if $i$ is an index of
$E$ and $p \in E$, then the first-order generalized linear span
${\rm GS}(E,p)$ contains the line through $p$ in the direction
$e_i$.


\begin{lem} \label{lem:2018-9-18}
Assume that $E$ is a bounded subset of $M_a$ and
${\rm GS}(E,p)$ is the first-order generalized linear span of
$E$ with respect to a fixed element $p \in E$.
If $i \in \Lambda(E)$, then $p + \alpha e_i \in {\rm GS}(E,p)$
for all $\alpha \in \mathbb{R}$.
\end{lem}

\noindent
\emph{Proof}.
By Definition \ref{def:Lambda} $(i)$, if $i \in \Lambda(E)$
then there exist an $x \in E$ and an $\alpha_0 \neq 0$, which
satisfy $x + \alpha_0 e_i \in E$.
Since $x \in E$ and $x + \alpha_0 e_i \in E$, by Definition
\ref{def:gls}, we get
\begin{eqnarray*}
p + \alpha_0 \beta e_i
=   p + \beta (x + \alpha_0 e_i - p) - \beta (x - p)
\in {\rm GS}(E,p)
\end{eqnarray*}
for all $\beta \in \mathbb{R}$.
Setting $\alpha = \alpha_0 \beta$ in the above relation, we
obtain $p + \alpha e_i \in {\rm GS}(E,p)$ for any
$\alpha \in \mathbb{R}$.
\hfill$\Box$
\vspace*{5mm}

We now introduce a lemma, which is a generalized version of
\cite[Lemma 2.3]{jungkim} and whose proof runs in the same way.
We prove that the function
$T_{-q} \circ f \circ T_p : E_1 - p \to E_2 - q$ preserves the
inner product.
This property is important in proving the following theorems
as a necessary condition for $f$ to be a $d_a$-isometry.


\begin{lem} \label{lem:new2}
Assume that $E_1$ and $E_2$ are bounded subsets of $M_a$ that
are $d_a$-isometric to each other via a surjective
$d_a$-isometry $f : E_1 \to E_2$.
Assume that $p$ is an element of $E_1$ and $q$ is an element
of $E_2$ with $q = f(p)$.
Then the function $T_{-q} \circ f \circ T_p : E_1 - p \to E_2 - q$
preserves the inner product, \textit{i.e.},
\begin{eqnarray*}
\big\langle ( T_{-q} \circ f \circ T_p )(x-p),\,
            ( T_{-q} \circ f \circ T_p )(y-p)
\big\rangle_a = \langle x-p, y-p \rangle_a
\end{eqnarray*}
for all $x, y \in E_1$.
\end{lem}

\noindent
\emph{Proof}.
Since $T_{-q} \circ f \circ T_p : E_1 - p \to E_2 - q$ is a
$d_a$-isometry, we have
\begin{eqnarray*}
\| ( T_{-q} \circ f \circ T_p )(x-p) -
   ( T_{-q} \circ f \circ T_p )(y-p)
\|_a = \| (x-p) - (y-p) \|_a
\end{eqnarray*}
for any $x, y \in E_1$.
If we put $y = p$ in the last equality, then we get
\begin{eqnarray*}
\| ( T_{-q} \circ f \circ T_p )(x-p) \|_a = \| x-p \|_a
\end{eqnarray*}
for each $x \in E_1$.
Moreover, it follows from the previous equality that
\begin{eqnarray*}
\begin{split}
& \| ( T_{-q} \circ f \circ T_p )(x-p) -
     ( T_{-q} \circ f \circ T_p )(y-p)
  \|_a^2 \\
& \hspace*{3mm}
  = \big\langle ( T_{-q} \circ f \circ T_p )(x-p) -
                ( T_{-q} \circ f \circ T_p )(y-p), \\
& \hspace*{10mm}
                ( T_{-q} \circ f \circ T_p )(x-p) -
                ( T_{-q} \circ f \circ T_p )(y-p)
    \big\rangle_a \\
& \hspace*{3mm}
  = \| x-p \|_a^2 -
    2 \big\langle ( T_{-q} \circ f \circ T_p )(x-p),\,
                  ( T_{-q} \circ f \circ T_p )(y-p)
      \big\rangle_a + \| y-p \|_a^2
\end{split}
\end{eqnarray*}
and
\begin{eqnarray*}
\begin{split}
\| (x-p) - (y-p) \|_a^2
& = \big\langle (x-p) - (y-p), (x-p) - (y-p) \big\rangle_a \\
& = \| x-p \|_a^2 - 2 \langle x-p, y-p \rangle_a +
    \| y-p \|_a^2.
\end{split}
\end{eqnarray*}
Finally, comparing the last two equalities yields the validity
of our assertion.
\hfill$\Box$
\vspace*{5mm}


\section{First-order extension of isometries}

In the previous section, we made all the necessary preparations
to extend the domain $E_1$ of the surjective $d_a$-isometry
$f : E_1 \to E_2$ to its first-order generalized linear span
${\rm GS}(E_1,p)$.

Although $E_1$ is a bounded set, ${\rm GS}(E_1,p) - p$ is a
real vector space.
Now we will extend the $d_a$-isometry $T_{-q} \circ f \circ T_p$
defined on the bounded set $E_1 - p$ to the $d_a$-isometry
$T_{-q} \circ F \circ T_p$ defined on the vector space
${\rm GS}(E_1,p) - p$.
Comparing their `sizes' of $E_1 - p$ and ${\rm GS}(E_1,p) - p$,
or considering that ${\rm GS}(E_1,p) - p$ is an algebraically
closed space, it is a great achievement to extend the
$d_a$-isometry $T_{-q} \circ f \circ T_p$ defined on the
bounded set $E_1 - p$ to the $d_a$-isometry defined on the
vector space ${\rm GS}(E_1,p) - p$.


\begin{defn} \label{def:ext}
Assume that $E_1$ and $E_2$ are nonempty bounded subsets of
$M_a$ that are $d_a$-isometric to each other via a surjective
$d_a$-isometry $f : E_1 \to E_2$.
Let $p$ be a fixed element of $E_1$ and let $q$ be an element
of $E_2$ that satisfies $q = f(p)$.
We define a function $F : {\rm GS}(E_1,p) \to M_a$ as
\begin{eqnarray*} 
(T_{-q} \circ F \circ T_p)
\Bigg( \sum_{i=1}^m \sum_{j=1}^\infty
       \alpha_{ij} (x_{ij} - p) \Bigg)
= \sum_{i=1}^m \sum_{j=1}^\infty
  \alpha_{ij} (T_{-q} \circ f \circ T_p)(x_{ij} - p)
\end{eqnarray*}
for any $m \in \mathbb{N}$, $x_{ij} \in E_1$, and for all
$\alpha_{ij} \in \mathbb{R}$ satisfying
$\sum\limits_{i=1}^m \sum\limits_{j=1}^\infty
 \alpha_{ij} (x_{ij} - p) \in M_a$.
\end{defn}

We note that in the definition above, it is important for the
argument of $T_{-q} \circ F \circ T_p$ to belong to $M_a$.
Now we show that the function $F : {\rm GS}(E_1,p) \to M_a$ is
well defined.


\begin{lem} \label{lem:3.2}
Assume that $E_1$ and $E_2$ are bounded subsets of $M_a$ that
are $d_a$-isometric to each other via a surjective
$d_a$-isometry $f : E_1 \to E_2$.
Let $p$ be an element of $E_1$ and let $q$ be an element of
$E_2$ that satisfy $q = f(p)$.
The function $F : {\rm GS}(E_1,p) \to M_a$ given in Definition
\ref{def:ext} is well defined.
\end{lem}

\noindent
\emph{Proof}.
First, we will check that the range of $F$ is a subset of $M_a$.
For any $m, n_1, n_2 \in \mathbb{N}$ with $n_2 > n_1$,
$x_{ij} \in E_1$, and for all $\alpha_{ij} \in \mathbb{R}$, it
follows from Lemma \ref{lem:new2} that
\begin{eqnarray} \label{eq:2020-1-10}
\begin{split}
& \bigg\| \sum_{i=1}^m \sum_{j=1}^{n_2}
          \alpha_{ij} (T_{-q} \circ f \circ T_p)(x_{ij} - p) -
          \sum_{i=1}^m \sum_{j=1}^{n_1}
          \alpha_{ij} (T_{-q} \circ f \circ T_p)(x_{ij} - p)
  \bigg\|_a^2 \\
& \hspace*{3mm} =
  \bigg\langle
  \sum_{i=1}^m \sum_{j = n_1 + 1}^{n_2}
  \alpha_{ij} (T_{-q} \circ f \circ T_p)(x_{ij} - p),\,
  \sum_{k=1}^m \sum_{\ell = n_1 + 1}^{n_2}
  \alpha_{k\ell} (T_{-q} \circ f \circ T_p)(x_{k\ell} - p)
  \bigg\rangle_a \\
& \hspace*{3mm} =
  \sum_{i=1}^m \sum_{k=1}^m
  \sum_{j = n_1 + 1}^{n_2} \alpha_{ij}
  \sum_{\ell = n_1 + 1}^{n_2} \alpha_{k\ell}
  \big\langle (T_{-q} \circ f \circ T_p)(x_{ij} - p),\,
              (T_{-q} \circ f \circ T_p)(x_{k\ell} - p)
  \big\rangle_a \\
& \hspace*{3mm} =
  \sum_{i=1}^m \sum_{k=1}^m
  \sum_{j = n_1 + 1}^{n_2} \alpha_{ij}
  \sum_{\ell = n_1 + 1}^{n_2} \alpha_{k\ell}
  \big\langle x_{ij} - p,\, x_{k\ell} - p \big\rangle_a \\
& \hspace*{3mm} =
  \bigg\langle
  \sum_{i=1}^m \sum_{j = n_1 + 1}^{n_2} \alpha_{ij}
  (x_{ij} - p),\,
  \sum_{k=1}^m \sum_{\ell = n_1 + 1}^{n_2} \alpha_{k\ell}
  (x_{k\ell} - p)
  \bigg\rangle_a \\
& \hspace*{3mm} =
  \bigg\|
  \sum_{i=1}^m \sum_{j = n_1 + 1}^{n_2} \alpha_{ij} (x_{ij} - p)
  \bigg\|_a^2 \\
& \hspace*{3mm} =
  \bigg\| \sum_{i=1}^m \sum_{j=1}^{n_2}
          \alpha_{ij} (x_{ij} - p) -
          \sum_{i=1}^m \sum_{j=1}^{n_1}
          \alpha_{ij} (x_{ij} - p)
  \bigg\|_a^2.
\end{split}
\end{eqnarray}
Indeed, the equality (\ref{eq:2020-1-10}) holds for all
$m, n_1, n_2 \in \mathbb{N}$.

We now assume that
$\sum\limits_{i=1}^m \sum\limits_{j=1}^\infty
 \alpha_{ij} (x_{ij} - p) \in M_a$ for some $x_{ij} \in E_1$
and $\alpha_{ij} \in \mathbb{R}$, where $m$ is a fixed positive
integer.
Then since $(M_a, \mathcal{T}_a)$ is a Hausdorff space on
account of Remark \ref{rem:2.1} $(ii)$ and the topology
$\mathcal{T}_a$ is consistent with the metric $d_a$ and with
the norm $\| \cdot \|_a$ (\textit{cf.} Remark \ref{rem:1.1}
$(i)$), the sequence
$\big\{ \sum\limits_{i=1}^m \sum\limits_{j=1}^n
        \alpha_{ij} (x_{ij} - p) \big\}_n$ converges to
$\sum\limits_{i=1}^m \sum\limits_{j=1}^\infty
 \alpha_{ij} (x_{ij} - p)$ (in $M_a$) and hence, the sequence
$\big\{ \sum\limits_{i=1}^m \sum\limits_{j=1}^n
        \alpha_{ij} (x_{ij} - p) \big\}_n$ is a Cauchy sequence
in $M_a$.

We know by (\ref{eq:2020-1-10}) and the definition of Cauchy
sequences that for each $\varepsilon > 0$ there exists an
integer $N_\varepsilon > 0$ such that
\begin{eqnarray*}
\begin{split}
& \bigg\| \sum_{i=1}^m \sum_{j=1}^{n_2}
          \alpha_{ij} (T_{-q} \circ f \circ T_p)(x_{ij} - p) -
          \sum_{i=1}^m \sum_{j=1}^{n_1}
          \alpha_{ij} (T_{-q} \circ f \circ T_p)(x_{ij} - p)
  \bigg\|_a \\
& \hspace*{3mm} =
  \bigg\| \sum_{i=1}^m \sum_{j=1}^{n_2}
          \alpha_{ij} (x_{ij} - p) -
          \sum_{i=1}^m \sum_{j=1}^{n_1}
          \alpha_{ij} (x_{ij} - p)
  \bigg\|_a < \varepsilon
\end{split}
\end{eqnarray*}
for all integers $n_1, n_2 > N_\varepsilon$, which implies that
$\big\{ \sum\limits_{i=1}^m \sum\limits_{j=1}^n
        \alpha_{ij} (T_{-q} \circ f \circ T_p)(x_{ij} - p)
 \big\}_n$ is also a Cauchy sequence in $M_a$.
As we proved in \cite[Theorem 2.1]{jungkim} or by Remark
\ref{rem:2.1} $(i)$, we observe that
$(M_a, \langle \cdot, \cdot \rangle_a)$ is a real Hilbert space
when the sequence $a = \{ a_i \}$ satisfies the condition
(\ref{eq:20110922-1}).
Thus, $M_a$ is not only complete, but also a Hausdorff space,
so the Cauchy sequence
$\big\{ \sum\limits_{i=1}^m \sum\limits_{j=1}^n
        \alpha_{ij} (T_{-q} \circ f \circ T_p)(x_{ij} - p)
 \big\}_n$ converges in $M_a$, \textit{i.e.}, by Definition
\ref{def:ext}, we have
\begin{eqnarray*}
\begin{split}
(T_{-q} \circ F \circ T_p)
\Bigg( \sum_{i=1}^m \sum_{j=1}^\infty \alpha_{ij} (x_{ij} - p)
\Bigg)
& =   \sum_{i=1}^m \sum_{j=1}^\infty \alpha_{ij}
      (T_{-q} \circ f \circ T_p)(x_{ij} - p) \\
& =   \lim_{n \to \infty}
      \sum_{i=1}^m \sum_{j=1}^n
      \alpha_{ij} (T_{-q} \circ f \circ T_p)(x_{ij} - p) \\
& \in M_a,
\end{split}
\end{eqnarray*}
which implies
\begin{eqnarray*}
F \Bigg( p + \sum_{i=1}^m \sum_{j=1}^\infty
         \alpha_{ij} (x_{ij} - p)
  \Bigg) \in M_a + q = M_a
\end{eqnarray*}
for all $x_{ij} \in E_1$ and $\alpha_{ij} \in \mathbb{R}$ with
$\sum\limits_{i=1}^m \sum\limits_{j=1}^\infty
 \alpha_{ij} (x_{ij} - p) \in M_a$, \textit{i.e.}, the image
of each element of ${\rm GS}(E_1,p)$ under $F$ belongs to
$M_a$.

We now assume that
$\sum\limits_{i=1}^{m_1} \sum\limits_{j=1}^\infty
 \alpha_{ij} (x_{ij} - p)
 = \sum\limits_{i=1}^{m_2} \sum\limits_{j=1}^\infty
   \beta_{ij} (y_{ij} - p) \in M_a$ for some
$m_1, m_2 \in \mathbb{N}$, $x_{ij}, y_{ij} \in E_1$, and
for some $\alpha_{ij}, \beta_{ij} \in \mathbb{R}$.
It then follows from Definition \ref{def:ext} and Lemma
\ref{lem:new2} that
\begin{eqnarray*}
\begin{split}
& \bigg\|
  \big( T_{-q} \circ F \circ T_p \big)
  \bigg( \sum_{i=1}^{m_1} \sum_{j=1}^\infty
         \alpha_{ij} (x_{ij} - p)
  \bigg) -
  \big( T_{-q} \circ F \circ T_p \big)
  \bigg( \sum_{i=1}^{m_2} \sum_{j=1}^\infty
         \beta_{ij} (y_{ij} - p)
  \bigg)
  \bigg\|_a^2 \\
& \hspace*{3mm}
  = \bigg\|
    \sum_{i=1}^{m_1} \sum_{j=1}^\infty \alpha_{ij}
    \big( T_{-q} \circ f \circ T_p \big) (x_{ij} - p) -
    \sum_{i=1}^{m_2} \sum_{j=1}^\infty \beta_{ij}
    \big( T_{-q} \circ f \circ T_p \big) (y_{ij} - p)
    \bigg\|_a^2 \\
& \hspace*{3mm}
  = \bigg\langle
    \sum_{i=1}^{m_1} \sum_{j=1}^\infty \alpha_{ij}
    \big( T_{-q} \circ f \circ T_p \big) (x_{ij} - p) -
    \sum_{i=1}^{m_2} \sum_{j=1}^\infty \beta_{ij}
    \big( T_{-q} \circ f \circ T_p \big) (y_{ij} - p), \\
& \hspace*{11mm}
    \sum_{k=1}^{m_1} \sum_{\ell=1}^\infty \alpha_{k\ell}
    \big( T_{-q} \circ f \circ T_p \big) (x_{k\ell} - p) -
    \sum_{k=1}^{m_2} \sum_{\ell=1}^\infty \beta_{k\ell}
    \big( T_{-q} \circ f \circ T_p \big) (y_{k\ell} - p)
    \bigg\rangle_a \\
& \hspace*{3mm}
  = \cdots \\
& \hspace*{3mm}
  = \bigg\langle
    \sum_{i=1}^{m_1} \sum_{j=1}^\infty \alpha_{ij}
    (x_{ij} - p) -
    \sum_{i=1}^{m_2} \sum_{j=1}^\infty \beta_{ij}
    (y_{ij} - p),\, \\
& \hspace*{11mm}
    \sum_{k=1}^{m_1} \sum_{\ell=1}^\infty \alpha_{k\ell}
    (x_{k\ell} - p) -
    \sum_{k=1}^{m_2} \sum_{\ell=1}^\infty \beta_{k\ell}
    (y_{k\ell} - p)
    \bigg\rangle_a \\
& \hspace*{3mm}
  = \bigg\|
    \sum_{i=1}^{m_1} \sum_{j=1}^\infty \alpha_{ij} (x_{ij} - p) -
    \sum_{i=1}^{m_2} \sum_{j=1}^\infty \beta_{ij} (y_{ij} - p)
    \bigg\|_a^2 \\
& \hspace*{3mm}
  = 0,
\end{split}
\end{eqnarray*}
which implies that
\begin{eqnarray*}
\big( T_{-q} \circ F \circ T_p \big)
\bigg( \sum_{i=1}^{m_1} \sum_{j=1}^\infty
       \alpha_{ij} (x_{ij} - p)
\bigg)
= \big( T_{-q} \circ F \circ T_p \big)
  \bigg( \sum_{i=1}^{m_2} \sum_{j=1}^\infty
         \beta_{ij} (y_{ij} - p)
  \bigg)
\end{eqnarray*}
for all $m_1, m_2 \in \mathbb{N}$, $x_{ij}, y_{ij} \in E_1$,
and for all $\alpha_{ij}, \beta_{ij} \in \mathbb{R}$ satisfying
$\sum\limits_{i=1}^{m_1} \sum\limits_{j=1}^\infty
 \alpha_{ij} (x_{ij} - p)
 = \sum\limits_{i=1}^{m_2} \sum\limits_{j=1}^\infty
   \beta_{ij} (y_{ij} - p) \in M_a$.
\hfill$\Box$
\vspace*{5mm}

In \cite[Theorem 2.2]{jung1}, we were able to extend the domain
of a $d_a$-isometry $f : J \to K$ to the whole space $M_a$ when
$J$ is a non-degenerate basic cylinder, while we prove in the
following theorem that the domain of a $d_a$-isometry
$f : E_1 \to E_2$ can be extended to the first-order generalized
linear span ${\rm GS}(E_1,p)$ whenever $E_1$ is a nonempty
bounded subset of $M_a$, whether degenerate or non-degenerate.
Therefore, Theorem \ref{thm:new3} is a generalization of
\cite[Theorem 2.2]{jung1}.

In the proof, we use the fact that ${\rm GS}(E_1,p) - p$ is a
real vector space.
This fact is self-evident, as briefly mentioned earlier.


\begin{thm} \label{thm:new3}
Assume that $E_1$ and $E_2$ are bounded subsets of $M_a$ that
are $d_a$-isometric to each other via a surjective
$d_a$-isometry $f : E_1 \to E_2$.
Assume that $p$ is an element of $E_1$ and $q$ is an element
of $E_2$ with $q = f(p)$.
The function $F : {\rm GS}(E_1,p) \to M_a$ defined in Definition
\ref{def:ext} is a $d_a$-isometry and the function
$T_{-q} \circ F \circ T_p : {\rm GS}(E_1,p) - p \to M_a$ is a
linear $d_a$-isometry.
In particular, $F$ is an extension of $f$.
\end{thm}

\noindent
\emph{Proof}.
$(a)$
Let $u$ and $v$ be arbitrary elements of the first-order
generalized linear span ${\rm GS}(E_1,p)$ of $E_1$ with
respect to $p$.
Then
\begin{eqnarray} \label{eq:2018-8-23-1}
u - p = \sum_{i=1}^m \sum_{j=1}^\infty
        \alpha_{ij} (x_{ij} - p) \in M_a
~~~\mbox{and}~~~
v - p = \sum_{i=1}^n \sum_{j=1}^\infty
        \beta_{ij} (y_{ij} - p) \in M_a
\end{eqnarray}
for some $m, n \in \mathbb{N}$, some
$x_{ij}, y_{ij} \in E_1$, and for some
$\alpha_{ij}, \beta_{ij} \in \mathbb{R}$.
Then, according to Definition \ref{def:ext}, we have
\begin{eqnarray} \label{eq:2018-8-23-2}
\begin{split}
(T_{-q} \circ F \circ T_p)(u-p)
& = \sum_{i=1}^m \sum_{j=1}^\infty
    \alpha_{ij} (T_{-q} \circ f \circ T_p)(x_{ij} - p), \\
(T_{-q} \circ F \circ T_p)(v-p)
& = \sum_{i=1}^n \sum_{j=1}^\infty
    \beta_{ij} (T_{-q} \circ f \circ T_p)(y_{ij} - p).
\end{split}
\end{eqnarray}

$(b)$
By Lemma \ref{lem:new2}, (\ref{eq:2018-8-23-1}), and
(\ref{eq:2018-8-23-2}), we get
\begin{eqnarray} \label{eq:2018-9-28}
\begin{split}
& \big\langle ( T_{-q} \circ F \circ T_p )(u-p),\,
              ( T_{-q} \circ F \circ T_p )(v-p)
  \big\rangle_a \\
& \hspace*{3mm}
  = \bigg\langle \,\sum_{i=1}^m \sum_{j=1}^\infty \alpha_{ij}
                 (T_{-q} \circ f \circ T_p)(x_{ij} - p),\,
                 \sum_{k=1}^n \sum_{\ell=1}^\infty \beta_{k\ell}
                 (T_{-q} \circ f \circ T_p)(y_{k\ell} - p)
    \bigg\rangle_a \\
& \hspace*{3mm}
  = \sum_{i=1}^m \sum_{k=1}^n
    \sum_{j=1}^\infty \alpha_{ij}
    \sum_{\ell=1}^\infty \beta_{k\ell}
    \big\langle (T_{-q} \circ f \circ T_p)(x_{ij} - p),\,
                (T_{-q} \circ f \circ T_p)(y_{k\ell} - p)
    \big\rangle_a \\
& \hspace*{3mm}
  = \sum_{i=1}^m \sum_{k=1}^n
    \sum_{j=1}^\infty \alpha_{ij}
    \sum_{\ell=1}^\infty \beta_{k\ell}
    \langle x_{ij} - p, y_{k\ell} - p \rangle_a \\
& \hspace*{3mm}
  = \bigg\langle \,\sum_{i=1}^m \sum_{j=1}^\infty \alpha_{ij}
                 (x_{ij} - p),\,
                 \sum_{k=1}^n \sum_{\ell=1}^\infty \beta_{k\ell}
                 (y_{k\ell} - p)
    \bigg\rangle_a \\
& \hspace*{3mm}
  = \langle u - p, v - p \rangle_a
\end{split}
\end{eqnarray}
for all $u, v \in {\rm GS}(E_1,p)$.
That is, $T_{-q} \circ F \circ T_p$ preserves the inner product.
Indeed, equality (\ref{eq:2018-9-28}) is an extended version
of Lemma \ref{lem:new2}.

$(c)$
By using equality (\ref{eq:2018-9-28}), we further obtain
\begin{eqnarray*}
\begin{split}
d_a \big( F(u), F(v) \big)^2
& = \| F(u) - F(v) \|_a^2 \\
& = \big\| ( T_{-q} \circ F \circ T_p )(u-p) -
           ( T_{-q} \circ F \circ T_p )(v-p)
    \big\|_a^2 \\
& = \big\langle ( T_{-q} \circ F \circ T_p )(u-p) -
                ( T_{-q} \circ F \circ T_p )(v-p), \\
&   ~~~~~~      ( T_{-q} \circ F \circ T_p )(u-p) -
                ( T_{-q} \circ F \circ T_p )(v-p)
    \big\rangle_a \\
& = \langle u-p, u-p \rangle_a - \langle u-p, v-p \rangle_a
    - \langle v-p, u-p \rangle_a + \langle v-p, v-p \rangle_a \\
& = \big\langle (u-p) - (v-p), (u-p) - (v-p) \big\rangle_a \\
& = \| (u-p) - (v-p) \|_a^2 \\
& = \| u - v \|_a^2 \\
& = d_a(u, v)^2
\end{split}
\end{eqnarray*}
for all $u, v \in {\rm GS}(E_1,p)$, \textit{i.e.}, $F$ is a
$d_a$-isometry.

$(d)$
Now, let $u$ and $v$ be arbitrary elements of ${\rm GS}(E_1,p)$.
Then, it holds that $u - p \in {\rm GS}(E_1,p) - p$,
$v - p \in {\rm GS}(E_1,p) - p$, and
$\alpha (u - p) + \beta (v - p) \in {\rm GS}(E_1,p) - p$ for
any $\alpha, \beta \in \mathbb{R}$, because
${\rm GS}(E_1,p) - p$ is a real vector space.

We get
\begin{eqnarray*}
\begin{split}
& \big\| (T_{-q} \circ F \circ T_p)
         \left( \alpha (u - p) + \beta (v - p) \right) \\
& \hspace*{3mm}
         -\alpha (T_{-q} \circ F \circ T_p) (u - p)
         -\beta (T_{-q} \circ F \circ T_p) (v - p)
  \big\|_a^2 \\
& \hspace*{3mm}
  = \Big\langle
    (T_{-q} \circ F \circ T_p)
    \left( \alpha (u - p) + \beta (v - p) \right) \\
& \hspace*{11mm}
    -\alpha (T_{-q} \circ F \circ T_p) (u - p)
    -\beta (T_{-q} \circ F \circ T_p) (v - p), \\
& \hspace*{11mm}
    (T_{-q} \circ F \circ T_p)
    \left( \alpha (u - p) + \beta (v - p) \right) \\
& \hspace*{11mm}
    -\alpha (T_{-q} \circ F \circ T_p) (u - p)
    -\beta (T_{-q} \circ F \circ T_p) (v - p)
    \Big\rangle_a. 
\end{split}
\end{eqnarray*}
Since $\alpha (u - p) + \beta (v - p) = w - p$ for some
$w \in {\rm GS}(E_1,p)$, we further use (\ref{eq:2018-9-28})
to obtain
\begin{eqnarray*}
\begin{split}
& \big\| (T_{-q} \circ F \circ T_p)
         \left( \alpha (u - p) + \beta (v - p) \right) \\
& \hspace*{3mm}
         -\alpha (T_{-q} \circ F \circ T_p) (u - p)
         -\beta (T_{-q} \circ F \circ T_p) (v - p)
  \big\|_a^2 \\
& \hspace*{3mm}
  = \langle w - p,\, w - p \rangle_a -
    \alpha \langle w - p,\, u - p \rangle_a -
    \beta \langle w - p,\, v - p \rangle_a \\
& \hspace*{8mm}
    -\alpha \langle u - p,\, w - p \rangle_a
    +\alpha^2 \langle u - p,\, u - p \rangle_a
    +\alpha \beta \langle u - p,\, v - p \rangle_a \\
& \hspace*{8mm}
    -\beta \langle v - p,\, w - p \rangle_a
    +\alpha \beta \langle v - p,\, u - p \rangle_a
    +\beta^2 \langle v - p,\, v - p \rangle_a \\
& \hspace*{3mm}
  = 0,
\end{split}
\end{eqnarray*}
which implies that the function
$T_{-q} \circ F \circ T_p : {\rm GS}(E_1,p) - p \to M_a$ is
linear.

$(e)$
Finally, we set $\alpha_{11} = 1$, $\alpha_{ij} = 0$ for any
$(i,j) \neq (1,1)$, and $x_{11} = x$ in (\ref{eq:2018-8-23-1})
and (\ref{eq:2018-8-23-2}) to see
\begin{eqnarray*}
(T_{-q} \circ F \circ T_p)(x-p)
= (T_{-q} \circ f \circ T_p)(x-p)
\end{eqnarray*}
for every $x \in E_1$, which implies that $F(x) = f(x)$ for
every $x \in E_1$, \textit{i.e.}, $F$ is an extension of $f$.
\hfill$\Box$
\vspace*{5mm}


\section{Second-order generalized linear span}

For any element $x$ of $M_a$ and $r > 0$, we denote by $B_r(x)$
the open ball defined by
$B_r(x) = \{ y \in M_a : \| y - x \|_a < r \}$.

Definitions \ref{def:gls} and \ref{def:ext} will be generalized
to the cases of $n \geq 2$ in the following definition.
We introduce the concept of $n$th-order generalized linear
span ${\rm GS}^n(E_1,p)$, which generalizes the concept of
first-order generalized linear span ${\rm GS}(E,p)$.
Moreover, we define the $d_a$-isometry $F_n$ which extends the
domain of a $d_a$-isometry $f$ to ${\rm GS}^n(E_1,p)$.

It is surprising, however, that this process of generalization
does not go far.
Indeed, we will find in Proposition \ref{rem:new3.1} and
Theorem \ref{thm:setyol} that ${\rm GS}^2(E_1,p)$ and $F_2$
are their limits.


\begin{defn} \label{def:6.1}
Let $E_1$ be a nonempty bounded subset of $M_a$ that is
$d_a$-isometric to a subset $E_2$ of $M_a$ via a surjective
$d_a$-isometry $f : E_1 \to E_2$.
Let $p$ be an element of $E_1$ and $q$ an element of $E_2$
with $q = f(p)$.
Assume that $r$ is a positive real number satisfying
$E_1 \subset B_r(p)$.
\begin{itemize}
\item[$(i)$]  We define ${\rm GS}^0(E_1,p) = E_1$ and
              ${\rm GS}^1(E_1,p) = {\rm GS}(E_1,p)$.
              In general, we define the
              \emph{$n$th-order generalized linear span} of
              $E_1$ with respect to $p$ as
              ${\rm GS}^n(E_1,p) = {\rm GS}
               ({\rm GS}^{n-1}(E_1,p) \cap B_r(p), p)$ for all
              $n \in \mathbb{N}$.
\item[$(ii)$] We define $F_0 = f$ and $F_1 = F$, where $F$ is
              defined in Definition \ref{def:ext}.
              Moreover, for any $n \in \mathbb{N}$, we define
              the function $F_n : {\rm GS}^n(E_1,p) \to M_a$
              by
              \begin{eqnarray*}
              (T_{-q} \circ F_n \circ T_p)
              \bigg( \sum_{i=1}^m \sum_{j=1}^\infty
                     \alpha_{ij} (x_{ij} - p) \bigg)
              = \sum_{i=1}^m \sum_{j=1}^\infty \alpha_{ij}
                (T_{-q} \circ F_{n-1} \circ T_p) (x_{ij} - p)
              \end{eqnarray*}
              for all $m \in \mathbb{N}$,
              $x_{ij} \in {\rm GS}^{n-1}(E_1,p) \cap B_r(p)$,
              and $\alpha_{ij} \in \mathbb{R}$ with
              $\sum\limits_{i=1}^m \sum\limits_{j=1}^\infty
               \alpha_{ij} (x_{ij} - p) \in M_a$.
\end{itemize}
\end{defn}


\begin{prop} \label{rem:31}
Let $E$ be a nonempty bounded subset of $M_a$.
If $s$ and $t$ are positive real numbers that satisfy
$E \subset B_s(p) \cap B_t(p)$, then
\begin{eqnarray*}
{\rm GS} \big( {\rm GS}(E,p) \cap B_s(p), p \big)
= {\rm GS} \big( {\rm GS}(E,p) \cap B_t(p), p \big).
\end{eqnarray*}
\end{prop}

\noindent
\emph{Proof}.
Assume that $0 < s < t$.
Then, there exists a real number $c > 1$ with $s > \frac{t}{c}$
and it is obvious that $B_{t/c}(p) \subset B_s(p)$.
Assume that $x$ is an arbitrary element of
${\rm GS}({\rm GS}(E,p) \cap B_t(p), p)$.
Then there exist some $m \in \mathbb{N}$, some
$u_{ij} \in {\rm GS}(E,p) \cap B_t(p)$ and some
$\alpha_{ij} \in \mathbb{R}$ such that
$x = p + \sum\limits_{i=1}^m \sum\limits_{j=1}^\infty
 \alpha_{ij} (u_{ij} - p) \in M_a$.
We notice that
\begin{eqnarray*}
\big( {\rm GS}(E,p) - p \big) \cap \big( B_t(p) - p \big)
= \big\{ u - p \in M_a : u \in {\rm GS}(E,p) \cap B_t(p) \big\}.
\end{eqnarray*}

Since ${\rm GS}(E,p) - p$ is a real vector space,
$\frac{t}{c} < s$, and since
$u_{ij} - p \in ({\rm GS}(E,p) - p) \cap (B_t(p) - p)$ for
any $i$ and $j$, we have
\begin{eqnarray*}
\frac{1}{c} (u_{ij} - p)
\in ({\rm GS}(E,p) - p) \cap (B_s(p) - p).
\end{eqnarray*}
Hence, we can choose a $v_{ij} \in {\rm GS}(E,p) \cap B_s(p)$
such that $\frac{1}{c} (u_{ij} - p) = v_{ij} - p$.
Thus, we get
\begin{eqnarray*}
x =   p + \sum_{i=1}^m \sum_{j=1}^\infty
      \alpha_{ij} (u_{ij} - p)
  =   p + \sum_{i=1}^m \sum_{j=1}^\infty
      c \alpha_{ij} (v_{ij} - p)
  \in {\rm GS}({\rm GS}(E,p) \cap B_s(p), p),
\end{eqnarray*}
which implies that
${\rm GS}({\rm GS}(E,p) \cap B_t(p), p)
 \subset {\rm GS}({\rm GS}(E,p) \cap B_s(p), p)$.

The reverse inclusion is obvious, since
$B_s(p) \subset B_t(p)$.
\hfill$\Box$
\vspace*{5mm}

We generalize Lemma \ref{lem:new2} and formula
(\ref{eq:2018-9-28}) in the following lemma.
Indeed, we prove that the function
$T_{-q} \circ F_n \circ T_p : {\rm GS}^n(E_1,p) - p \to M_a$
preserves the inner product.
This property is important in proving the following theorems
as a necessary condition for $F_n$ to be a $d_a$-isometry.


\begin{lem} \label{lem:6.1}
Let $E_1$ be a bounded subset of $M_a$ that is $d_a$-isometric
to a subset $E_2$ of $M_a$ via a surjective $d_a$-isometry
$f : E_1 \to E_2$.
Assume that $p$ and $q$ are elements of $E_1$ and $E_2$, which
satisfy $q = f(p)$.
If $n \in \mathbb{N}$, then
\begin{eqnarray*}
\big\langle (T_{-q} \circ F_n \circ T_p) (u - p),\,
            (T_{-q} \circ F_n \circ T_p) (v - p)
\big\rangle_a =
\langle u - p, v - p \rangle_a
\end{eqnarray*}
for all $u, v \in {\rm GS}^n(E_1,p)$.
\end{lem}

\noindent
\emph{Proof}.
Our assertion for $n = 1$ was already proved in
(\ref{eq:2018-9-28}).
Considering Proposition \ref{rem:31}, assume that $r$ is a
positive real number satisfying $E_1 \subset B_r(p)$.
Now we assume that the assertion is true for some
$n \in \mathbb{N}$.
Let $u, v$ be arbitrary elements of ${\rm GS}^{n+1}(E_1,p)$.
Then there exist some $m_1, m_2 \in \mathbb{N}$, some
$x_{ij}, y_{k\ell} \in {\rm GS}^n(E_1,p) \cap B_r(p)$ and some
$\alpha_{ij}, \beta_{k\ell} \in \mathbb{R}$ such that
\begin{eqnarray*}
u - p = \sum_{i=1}^{m_1} \sum_{j=1}^\infty
        \alpha_{ij} (x_{ij} - p) \in M_a ~~~\mbox{and}~~~
v - p = \sum_{k=1}^{m_2} \sum_{\ell=1}^\infty
        \beta_{k\ell} (y_{k\ell} - p) \in M_a.
\end{eqnarray*}

Using Definition \ref{def:6.1} $(ii)$ and our assumption, we
get
\begin{eqnarray*}
\begin{split}
& \big\langle ( T_{-q} \circ F_{n+1} \circ T_p )(u - p),\,
              ( T_{-q} \circ F_{n+1} \circ T_p )(v - p)
  \big\rangle_a \\
& \hspace*{3mm}
  = \bigg\langle \,\sum_{i=1}^{m_1} \sum_{j=1}^\infty
                 \alpha_{ij}
                 (T_{-q} \circ F_n \circ T_p)(x_{ij} - p),\,
                 \sum_{k=1}^{m_2} \sum_{\ell=1}^\infty
                 \beta_{k\ell}
                 (T_{-q} \circ F_n \circ T_p)(y_{k\ell} - p)
    \bigg\rangle_a \\
& \hspace*{3mm}
  = \sum_{i=1}^{m_1} \sum_{k=1}^{m_2}
    \sum_{j=1}^\infty \alpha_{ij}
    \sum_{\ell=1}^\infty \beta_{k\ell}
    \big\langle (T_{-q} \circ F_n \circ T_p)(x_{ij} - p),\,
                (T_{-q} \circ F_n \circ T_p)(y_{k\ell} - p)
    \big\rangle_a \\
& \hspace*{3mm}
  = \sum_{i=1}^{m_1} \sum_{k=1}^{m_2}
    \sum_{j=1}^\infty \alpha_{ij}
    \sum_{\ell=1}^\infty \beta_{k\ell}
    \langle x_{ij} - p, y_{k\ell} - p \rangle_a \\
& \hspace*{3mm}
  = \bigg\langle \,\sum_{i=1}^{m_1} \sum_{j=1}^\infty
                 \alpha_{ij} (x_{ij} - p),\,
                 \sum_{k=1}^{m_2} \sum_{\ell=1}^\infty
                 \beta_{k\ell} (y_{k\ell} - p)
    \bigg\rangle_a \\
& \hspace*{3mm}
  = \langle u - p, v - p \rangle_a
\end{split}
\end{eqnarray*}
for all $u, v \in {\rm GS}^{n+1}(E_1,p)$.
By mathematical induction, we may then conclude that our
assertion is true for all $n \in \mathbb{N}$.
\hfill$\Box$
\vspace*{5mm}

When $n = 1$ and $p = p'$, the first assertion in $(i)$ of the
following lemma is self-evident, so we have used that fact
several times before, omitting the proof.
The assertion $(iv)$ in the following lemma seems to be related
in some way to Proposition \ref{rem:31}.


\begin{lem} \label{lem:3.8}
Let $E$ be a bounded subset of $M_a$ and $p, p' \in E$.
Assume that $r$ is a positive real number satisfying
$E \subset B_r(p)$.
\begin{itemize}
\item[$(i)$]   ${\rm GS}^n(E,p) - p'$ is a vector space over
               $\mathbb{R}$ for each $n \in \mathbb{N}$.
\item[$(ii)$]  ${\rm GS}^n(E,p) \subset {\rm GS}^{n+1}(E,p)$
               for each $n \in \mathbb{N}$.
\item[$(iii)$] ${\rm GS}^2(E,p) = \overline{{\rm GS}(E,p)}$,
               where $\overline{{\rm GS}(E,p)}$ is the closure
               of\/ ${\rm GS}(E,p)$ in $M_a$.
\item[$(iv)$]  $\Lambda({\rm GS}^n(E,p)) =
                \Lambda({\rm GS}^n(E,p) \cap B_r(p))$ for all
               $n \in \mathbb{N}$.
\end{itemize}
\end{lem}

\noindent
\emph{Proof}.
$(i)$ By using Definitions \ref{def:gls} and \ref{def:6.1}, we
prove that ${\rm GS}(E,p) - p'$ is a real vector space.
(We can prove similarly for the case of $n > 1$.)
Given $x, y \in {\rm GS}(E,p) - p'$, we may choose some
$m_1, m_2 \in \mathbb{N}$, some $u_{ij}, v_{ij} \in E$, and
some $\alpha_{ij}, \beta_{ij} \in \mathbb{R}$ such that
$x = (p - p') + \sum\limits_{i=1}^{m_1} \sum\limits_{j=1}^\infty
     \alpha_{ij} (u_{ij} - p) \in M_a$ and
$y = (p - p') + \sum\limits_{i=1}^{m_2} \sum\limits_{j=1}^\infty
     \beta_{ij} (v_{ij} - p) \in M_a$.
Since $M_a$ is a real vector space,
$\alpha \sum\limits_{i=1}^{m_1} \sum\limits_{j=1}^\infty
 \alpha_{ij} (u_{ij} - p) +
 \beta \sum\limits_{i=1}^{m_2} \sum\limits_{j=1}^\infty
 \beta_{ij} (v_{ij} - p) \in M_a$ for all
$\alpha, \beta \in \mathbb{R}$.

Moreover, we see that
\begin{eqnarray*}
\begin{split}
\alpha x + \beta y
& =   \left(
      p + (1 - \alpha - \beta)(p' - p) +
      \sum_{i=1}^{m_1} \sum_{j=1}^\infty
      \alpha \alpha_{ij} (u_{ij} - p) +
      \sum_{i=1}^{m_2} \sum_{j=1}^\infty
      \beta \beta_{ij} (v_{ij} - p)
      \right) - p' \\
& \in {\rm GS}(E,p) - p'
\end{split}
\end{eqnarray*}
for all $\alpha, \beta \in \mathbb{R}$.
Hence, ${\rm GS}(E,p) - p'$ is a real vector space as a
subspace of real vector space $M_a$.

$(ii)$
Let $r$ be a positive real number with $E \subset B_r(p)$.
If $x \in {\rm GS}^n(E,p)$ for some $n \in \mathbb{N}$, then
$x - p \in {\rm GS}^n(E,p) - p$.
Since ${\rm GS}^n(E,p) - p$ is a real vector space by $(i)$
and $B_r(p) - p = B_r(0)$, we can choose a (sufficiently small)
real number $\mu \neq 0$ such that
$\mu (x - p) \in ({\rm GS}^n(E,p) - p) \cap (B_r(p) - p)$.
We notice that
\begin{eqnarray} \label{eq:2021-11-2}
\big( {\rm GS}^n(E,p) - p \big) \cap \big( B_r(p) - p \big)
= \big\{ v - p \in M_a : v \in {\rm GS}^n(E,p) \cap B_r(p)
  \big\}.
\end{eqnarray}
Thus, we see that $\mu (x - p) = v - p$ for some
$v \in {\rm GS}^n(E,p) \cap B_r(p)$.
Since $x = p + \frac{1}{\mu} (v - p)$, it holds that
$x \in {\rm GS}^{n+1}(E,p)$.
Therefore, we conclude that
${\rm GS}^n(E,p) \subset {\rm GS}^{n+1}(E,p)$ for every
$n \in \mathbb{N}$.

$(iii)$ Let $x$ be an arbitrary element of
$\overline{{\rm GS}(E,p)}$.
Then there exists some sequence $\{ x_n \}$ that converges to
$x$, where $x_n \in {\rm GS}(E,p) \!\setminus\! \{ x \}$ for
all $n \in \mathbb{N}$.
We now set $y_1 = x_1$ and $y_i = x_i - x_{i-1}$ for each
integer $i \geq 2$.
Then we have
\begin{eqnarray*}
x_n = \sum_{i=1}^n y_i,
\end{eqnarray*}
where $y_i = (x_i - p) - (x_{i-1} - p) \in {\rm GS}(E,p) - p$
for $i \geq 2$.
Since ${\rm GS}(E,p) - p$ is a real vector space and
$B_r(p) - p = B_r(0)$, we can select a real number
$\mu_i \neq 0$ such that
\begin{eqnarray*}
\mu_i y_i \in {\rm GS}(E,p) - p ~~~\mbox{and}~~~
\mu_i y_i \in B_r(p) - p
\end{eqnarray*}
for every integer $i \geq 2$.
Thus, it follows from (\ref{eq:2021-11-2}) that
\begin{eqnarray*}
x_n = \sum_{i=1}^n y_i
    = y_1 + \sum_{i=2}^n \frac{1}{\mu_i} (\mu_i y_i)
    = x_1 + \sum_{i=2}^n \frac{1}{\mu_i} (v_i - p),
\end{eqnarray*}
where $v_i \in {\rm GS}(E,p) \cap B_r(p)$ for $i \geq 2$.
Since the sequence $\{ x_n \}$ is assumed to converge to $x$,
the sequence
$\big\{ x_1 + \sum\limits_{i=2}^n \frac{1}{\mu_i} (v_i - p)
 \big\}_n$ converges to $x$.
Hence, we have
\begin{eqnarray} \label{eq:2021-10-16-1}
x_1 + \sum_{i=2}^\infty \frac{1}{\mu_i} (v_i - p)
=   \lim_{n \to \infty} x_n
=   x
\in M_a.
\end{eqnarray}
(Since $M_a$ is a Huasdorff space, $x$ is the unique limit
point of the sequence $\{ x_n \}$.)

Furthermore, there exists a real number $\mu_1 \neq 0$ that
satisfies $\mu_1 (x_1 - p) \in {\rm GS}(E,p) - p$ and
$\mu_1 (x_1 - p) \in B_r(p) - p$, \textit{i.e.},
$\mu_1 (x_1 - p) \in ({\rm GS}(E,p) - p) \cap (B_r(p) - p)$.
Thus, there exists a $v_1 \in {\rm GS}(E,p) \cap B_r(p)$ such
that $\mu_1 (x_1 - p) = v_1 - p$ or
$x_1 - p = \frac{1}{\mu_1} (v_1 - p)$.
Therefore,
\begin{eqnarray} \label{eq:2021-10-16-2}
x =   p + (x_1 - p) + \sum_{i=2}^\infty
      \frac{1}{\mu_i} (v_i - p)
  =   p + \sum_{i=1}^\infty \frac{1}{\mu_i} (v_i - p),
\end{eqnarray}
where $v_i \in {\rm GS}(E,p) \cap B_r(p)$ for each
$i \in \mathbb{N}$.
It follows from (\ref{eq:2021-10-16-1}) that
$\sum\limits_{i=1}^\infty \frac{1}{\mu_i} (v_i - p) \in M_a$.
Thus, by (\ref{eq:2021-10-16-2}), we see that
$x \in {\rm GS}^2(E,p)$, which implies that
$\overline{{\rm GS}(E,p)} \subset {\rm GS}^2(E,p)$.

On the other hand, let $y \in {\rm GS}^2(E,p)$.
Then there are some $m \in \mathbb{N}$, some
$v_{ij} \in {\rm GS}(E,p) \cap B_r(p)$, and some
$\alpha_{ij} \in \mathbb{R}$ such that
$y = p + \sum\limits_{i=1}^m \sum\limits_{j=1}^\infty
 \alpha_{ij} (v_{ij} - p) \in M_a$.
Let us define
$y_n = p + \sum\limits_{i=1}^m \sum\limits_{j=1}^n
 \alpha_{ij} (v_{ij} - p)$ for every $n \in \mathbb{N}$.
Since $v_{ij} - p \in {\rm GS}(E,p) - p$ for all $i$ and $j$
and ${\rm GS}(E,p) - p$ is a real vector space, we know that
$y_n - p = \sum\limits_{i=1}^m \sum\limits_{j=1}^n
 \alpha_{ij} (v_{ij} - p) \in {\rm GS}(E,p) - p$ and hence,
$y_n \in {\rm GS}(E,p)$ for all $n \in \mathbb{N}$.
Since ${\rm GS}(E,p)$ is a Hausdorff space, $y$ is the unique
limit point of the sequence $\{ y_n \}$.
Thus, we see that
\begin{eqnarray*}
y =   p + \sum_{i=1}^m \sum_{j=1}^\infty
      \alpha_{ij} (v_{ij} - p)
  =   \lim_{n \to \infty} y_n
  \in \overline{{\rm GS}(E,p)},
\end{eqnarray*}
which implies that
${\rm GS}^2(E,p) \subset \overline{{\rm GS}(E,p)}$.

$(iv)$ Let $i \in \Lambda({\rm GS}^n(E,p))$.
In view of Definition \ref{def:Lambda}, there exist
$x \in {\rm GS}^n(E,p)$ and $\alpha \neq 0$ with
$x + \alpha e_i \in {\rm GS}^n(E,p)$.
Further,
$x = p + \sum\limits_{i=1}^m \sum\limits_{j=1}^\infty
 \alpha_{ij} (u_{ij} - p)$ for some $m \in \mathbb{N}$, some
$u_{ij} \in {\rm GS}^{n-1}(E,p) \cap B_r(p)$, and for some
$\alpha_{ij} \in \mathbb{R}$.
Since
$\sum\limits_{i=1}^m \sum\limits_{j=1}^\infty
 \alpha_{ij} (u_{ij} - p) + \alpha e_i
 =   x - p + \alpha e_i
 \in {\rm GS}^n(E,p) - p$ and $B_r(p) - p = B_r(0)$, it holds
that
\begin{eqnarray*}
\mu \left( \sum\limits_{i=1}^m \sum\limits_{j=1}^\infty
           \alpha_{ij} (u_{ij} - p) + \alpha e_i
    \right)
\in ({\rm GS}^n(E,p) - p) \cap (B_r(p) - p)
\end{eqnarray*}
for any sufficiently small $\mu \neq 0$, or equivalently, it
follows from (\ref{eq:2021-11-2}) that
\begin{eqnarray} \label{eq:2021-4-28}
\left( p + \sum\limits_{i=1}^m \sum\limits_{j=1}^\infty
       \mu \alpha_{ij} (u_{ij} - p)
\right) + \mu \alpha e_i
\in {\rm GS}^n(E,p) \cap B_r(p).
\end{eqnarray}

On the other hand, since
$\sum\limits_{i=1}^m \sum\limits_{j=1}^\infty
 \alpha_{ij} (u_{ij} - p) = x - p
 \in {\rm GS}^n(E,p) - p$, it holds that
$\sum\limits_{i=1}^m \sum\limits_{j=1}^\infty
 \mu \alpha_{ij} (u_{ij} - p) \in ({\rm GS}^n(E,p) - p)
 \cap (B_r(p) - p)$ for any sufficiently small $\mu \neq 0$.
Hence, it follows from (\ref{eq:2021-11-2}) that
$p + \sum\limits_{i=1}^m \sum\limits_{j=1}^\infty
 \mu \alpha_{ij} (u_{ij} - p)  \in {\rm GS}^n(E,p)
 \cap B_r(p)$ for any sufficiently small $\mu \neq 0$.
Thus, by Definition \ref{def:Lambda} and (\ref{eq:2021-4-28}),
it holds that $i \in \Lambda({\rm GS}^n(E,p) \cap B_r(p))$,
which implies that
$\Lambda({\rm GS}^n(E,p))
 \subset \Lambda({\rm GS}^n(E,p) \cap B_r(p))$.
Obviously, the inverse inclusion is true.
\hfill$\Box$
\vspace*{5mm}

As we mentioned earlier, we will see that the second-order
generalized linear span is the last step in this kind of domain
extension.


\begin{prop} \label{rem:new3.1}
If $E$ is a bounded subset of $M_a$ and $p \in E$, then
\begin{eqnarray*}
E \subset {\rm GS}(E,p) \subset \overline{{\rm GS}(E,p)}
= {\rm GS}^2(E,p) = {\rm GS}^n(E,p)
\end{eqnarray*}
for any integer $n \geq 3$.
Indeed, ${\rm GS}^n(E,p) - p$ is a real Hilbert space for
$n \geq 2$.
\end{prop}

\noindent
\emph{Proof}.
$(a)$
Considering Proposition \ref{rem:31}, we can choose a real
number $r > 0$ that satisfies $E \subset B_r(p)$.
Assume that $x \in {\rm GS}^3(E,p)$.
Then there exist some $m_0 \in \mathbb{N}$, some
$u_{ij} \in {\rm GS}^2(E,p) \cap B_r(p)$, and some
$\alpha_{ij} \in \mathbb{R}$ such that
$x = p + \sum\limits_{i=1}^{m_0} \sum\limits_{j=1}^\infty
 \alpha_{ij} (u_{ij} - p) \in M_a$.

We define
$x_m = p + \sum\limits_{i=1}^{m_0} \sum\limits_{j=1}^m
 \alpha_{ij} (u_{ij} - p)$ for each $m \in \mathbb{N}$.
Since $u_{ij} \in {\rm GS}^2(E,p)$, there exist some
$m_{ij} \in \mathbb{N}$, some
$v_{ijk\ell} \in {\rm GS}(E,p) \cap B_r(p)$, and some
$\beta_{ijk\ell} \in \mathbb{R}$ such that
$u_{ij} = p + \sum\limits_{k=1}^{m_{ij}}
 \sum\limits_{\ell=1}^\infty
 \beta_{ijk\ell} (v_{ijk\ell} - p) \in M_a$.
Hence, it holds that
\begin{eqnarray*}
x_m =   p + \sum_{i=1}^{m_0} \sum_{j=1}^m \sum_{k=1}^{m_{ij}}
            \sum_{\ell=1}^\infty
            \alpha_{ij} \beta_{ijk\ell}
            \big( v_{ijk\ell} - p \big)
    \in M_a,
\end{eqnarray*}
which implies that $x_m \in {\rm GS}^2(E,p)$ for all
$m \in \mathbb{N}$.
Thus, $\{ x_m \}$ is a sequence in ${\rm GS}^2(E,p)$ that
converges to $x$.
Therefore, $x \in {\rm GS}^2(E,p)$ because ${\rm GS}^2(E,p)$
is closed.
Thus, ${\rm GS}^3(E,p) \subset {\rm GS}^2(E,p)$.
The inverse inclusion is of course true due to Lemma
\ref{lem:3.8} $(ii)$.
We have proved that ${\rm GS}^2(E,p) = {\rm GS}^3(E,p)$.

$(b)$ Assume that
${\rm GS}^2(E,p) = \cdots = {\rm GS}^n(E,p) = {\rm GS}^{n+1}(E,p)$
for some integer $n \geq 2$.

$(c)$ If we replace ${\rm GS}(E,p)$, ${\rm GS}^2(E,p)$, and
${\rm GS}^3(E,p)$ in the previous part $(a)$ with
${\rm GS}^n(E,p)$, ${\rm GS}^{n+1}(E,p)$, and
${\rm GS}^{n+2}(E,p)$, respectively, and if we consider the
fact that ${\rm GS}^{n+1}(E,p) = {\rm GS}^2(E,p)$ is closed in
$M_a$ by Lemma \ref{lem:3.8} $(iii)$ and our assumption $(b)$,
then we arrive at the conclusion that
${\rm GS}^{n+1}(E,p) = {\rm GS}^{n+2}(E,p)$.

$(d)$ With the conclusion of mathematical induction we prove
that ${\rm GS}^n(E,p) = {\rm GS}^2(E,p)$ for every integer
$n \geq 3$.
Moreover, when $n \geq 2$, ${\rm GS}^n(E,p)$ is complete as a
closed subset of a real Hilbert space $M_a$
(ref. Remark \ref{rem:2.1} and \cite[Theorem 63.3]{kasriel}).
Therefore, ${\rm GS}^n(E,p) - p$ is a real Hilbert space for
$n \geq 2$.
\hfill$\Box$
\vspace*{5mm}

The following lemma is an extension of Lemma \ref{lem:2018-9-18}
for the second-order generalized linear span ${\rm GS}^2(E,p)$.
Indeed, we prove that if $i \in \Lambda({\rm GS}^2(E,p))$, then
the second-order generalized linear span of $E$ contains all
the lines through ${\rm GS}(E,p)$ in the direction $e_i$.


\begin{lem} \label{lem:setahob}
Assume that a bounded subset $E$ of $M_a$ contains at least
two elements and $p \in E$.
If $i \in \Lambda({\rm GS}^2(E,p))$ and $p' \in {\rm GS}(E,p)$,
then $p' + \alpha_i e_i \in {\rm GS}^2(E,p)$ for any
$\alpha_i \in \mathbb{R}$.
\end{lem}

\noindent
\emph{Proof}.
Let $r$ be a positive real number with $E \subset B_r(p)$.
Assume that $i \in \Lambda({\rm GS}^2(E,p))$.
Considering Lemma \ref{lem:3.8} $(iv)$ and Proposition
\ref{rem:new3.1}, if we substitute
${\rm GS}^2(E,p) \cap B_{r}(p)$ for $E$ in Lemma
\ref{lem:2018-9-18}, then
$p + \alpha_i e_i \in {\rm GS}^3(E,p) = {\rm GS}^2(E,p)$ for
all $\alpha_i \in \mathbb{R}$.
Thus, there are some $m \in \mathbb{N}$, some
$w_{ij} \in {\rm GS}(E,p) \cap B_{r}(p)$, and some
$\beta_{ij} \in \mathbb{R}$ with
$\sum\limits_{i=1}^m \sum\limits_{j=1}^\infty
 \beta_{ij} (w_{ij} - p) \in M_a$ such that
$p + \alpha_i e_i
 = p + \sum\limits_{i=1}^m \sum\limits_{j=1}^\infty
   \beta_{ij} (w_{ij} - p)$, and hence, we have
\begin{eqnarray} \label{eq:2020-1-21-1}
p' + \alpha_i e_i
= p + \alpha_i e_i + (p' - p)
= p + \sum_{i=1}^m \sum_{j=1}^\infty
  \beta_{ij} (w_{ij} - p) + (p' - p).
\end{eqnarray}
Because $p' - p$ belongs to ${\rm GS}(E,p) - p$, which is a
real vector space by Lemma \ref{lem:3.8} $(i)$, and
$B_r(p) - p = B_r(0)$, we can choose some sufficiently small
real number $\mu \neq 0$ such that
\begin{eqnarray} \label{eq:2020-1-21-2}
\mu (p' - p) \in {\rm GS}(E,p) - p ~~~\mbox{and}~~~
\mu (p' - p) \in B_{r}(p) - p.
\end{eqnarray}

Considering (\ref{eq:2021-11-2}), (\ref{eq:2020-1-21-1}) and
(\ref{eq:2020-1-21-2}), if we put $\mu (p' - p) = w - p$ with
a $w \in {\rm GS}(E,p) \cap B_{r}(p)$, then we have
\begin{eqnarray*}
p' + \alpha_i e_i
=   p + \sum_{i=1}^m \sum_{j=1}^\infty
    \beta_{ij} (w_{ij} - p) + \frac{1}{\mu} (w - p)
\in {\rm GS}^2(E,p)
\end{eqnarray*}
for all $\alpha_i \in \mathbb{R}$.
\hfill$\Box$


\section{Basic cylinders and basic intervals}

In \cite[Theorem 2.5]{jungkim}, we could extend the domain of
a $d_a$-isometry $f$ to the whole space when the domain of $f$
is a non-degenerate basic cylinder (see the definition below
for the exact definition of non-degenerate basic cylinders).

Now we will define the infinite dimensional intervals more
precisely divided into non-degenerate basic cylinders,
degenerate basic cylinders, and basic intervals.


\begin{defn} \label{def:basic}
For any positive integer $n$, we define the infinite
dimensional interval by
\begin{eqnarray*}
J = \prod_{i=1}^\infty J_i, ~~~\mbox{where}~~~
J_i = \left\{\begin{array}{ll}
                \mbox{$[ 0, p_{2i} ]$}
                & (\mbox{for}~ i \in \Lambda_1),
                \vspace*{2mm}\\
                \mbox{$[ p_{1i}, p_{2i} ]$}
                & (\mbox{for}~ i \in \Lambda_2),
                \vspace*{2mm}\\
                \mbox{$[ p_{1i}, 1 ]$}
                & (\mbox{for}~ i \in \Lambda_3),
                \vspace*{2mm}\\
                \mbox{$\{ p_{1i} \}$}
                & (\mbox{for}~ i \in \Lambda_4),
                \vspace*{2mm}\\
                \mbox{$[ 0, 1 ]$}
                & (\mbox{otherwise})
             \end{array}
      \right.
\end{eqnarray*}
for some disjoint finite subsets $\Lambda_1$, $\Lambda_2$,
$\Lambda_3$ of $\{ 1, 2, \ldots, n \}$ and
$0 < p_{1i} < p_{2i} < 1$ for
$i \in \Lambda_1 \cup \Lambda_2 \cup \Lambda_3$ and
$0 \leq p_{1i} \leq 1$ for $i \in \Lambda_4$.
If $\Lambda_4 = \emptyset$, then $J$ is called a
\emph{non-degenerate basic cylinder}.
When $\Lambda_4$ is a nonempty finite set, $J$ is called a
\emph{degenerate basic cylinder}.
If $\Lambda_4$ is an infinite set, then $J$ will be called a
\emph{basic interval}.
\end{defn}


\begin{rem} \label{rem:6.1}
\begin{itemize}
\item[$(i)$]   In order for an infinite dimensional interval
               $J$ to become a basic cylinder, $\Lambda_4$
               must be a finite set.
\item[$(ii)$]  We remark that
               $\Lambda_4 = \mathbb{N} \!\setminus\! \Lambda(J)$
               and
               $\Lambda(J) = \mathbb{N} \!\setminus\! \Lambda_4$.
               That is, $\mathbb{N}$ is the disjoint union of
               $\Lambda(J)$ and $\Lambda_4$.
\item[$(iii)$] If $p = (p_1, p_2, \ldots, p_i, \ldots)$ is an
               element of an infinite dimensional interval $J$,
               then $J_i = \{ p_i \}$ for each
               $i \not\in \Lambda(J)$.
\end{itemize}
\end{rem}

We note that the basic cylinder or the basic interval $J$
defined in Definition \ref{def:basic} can be expressed as
\begin{eqnarray*}
J = \left\{ \sum_{i=1}^\infty \alpha_i\!
            \left( \frac{1}{a_i} e_i \right)\! :
            \alpha_i \in a_i J_i
            ~\mbox{for all}~ i \in \mathbb{N}
    \right\},
\end{eqnarray*}
where $J_i$ is the interval defined in Definition
\ref{def:basic}.


\begin{defn} \label{def:betabasic}
Let $\beta = \{ \beta_i \}_{i \in \mathbb{N}}$ be a complete
orthonormal sequence in $M_a$, $J_i$ the interval given in
Definition \ref{def:basic}, and let $n$ be a positive integer.
We define
\begin{eqnarray*}
J_\beta = \left\{ \sum_{i=1}^\infty \alpha_i \beta_i :
                  \alpha_i \in a_i J_i
                  ~\mbox{for all}~ i \in \mathbb{N}
          \right\},
\end{eqnarray*}
for some disjoint finite subsets $\Lambda_1$, $\Lambda_2$,
$\Lambda_3$ of $\{ 1, 2, \ldots, n \}$ and
$0 < p_{1i} < p_{2i} < 1$ for
$i \in \Lambda_1 \cup \Lambda_2 \cup \Lambda_3$ and
$0 \leq p_{1i} \leq 1$ for $i \in \Lambda_4$.
If $\Lambda_4 = \emptyset$, then $J_\beta$ is called a
\emph{non-degenerate $\beta$-basic cylinder}.
When $\Lambda_4$ is a nonempty finite set, $J_\beta$ is called
a \emph{degenerate $\beta$-basic cylinder}.
If $\Lambda_4$ is an infinite set, then $J_\beta$ will be
called a \emph{$\beta$-basic interval}.
\end{defn}

Using Definitions \ref{def:basic} and \ref{def:betabasic},
Remark \ref{rem:6.1} $(ii)$ is generalized to:


\begin{rem} \label{rem:6.2}
Let $\beta = \{ \beta_i \}_{i \in \mathbb{N}}$ be a complete
orthonormal sequence in $M_a$ and let $J_\beta$ be a
$\beta$-basic cylinder or a $\beta$-basic interval.
It holds that
$\Lambda_\beta(J_\beta) = \mathbb{N} \!\setminus\! \Lambda_4$,
where $\Lambda_4$ is given in Definitions \ref{def:basic} and
\ref{def:betabasic}.
\end{rem}

\noindent
\emph{Proof}.
In general, if $i \in \Lambda_4$, then it follows from
Definition \ref{def:betabasic} that
\begin{eqnarray*}
\langle x, \beta_i \rangle_a
=   \left\langle \sum_{j=1}^\infty \alpha_j \beta_j, \beta_i
    \right\rangle_a
=   \alpha_i
\in a_i J_i
=   \{ a_i p_{1i} \}
\end{eqnarray*}
for all $x \in J_\beta$.
That is, $\langle x, \beta_i \rangle_a = \alpha_i = a_i p_{1i}$
for all $x \in J_\beta$ and $i \in \Lambda_4$.
If $i \in \Lambda_4$, then
$\langle x + \alpha \beta_i, \beta_i \rangle_a
 = \langle x , \beta_i \rangle_a + \alpha
 = a_i p_{1i} + \alpha \neq a_i p_{1i}$ for all $x \in J_\beta$
and $\alpha \neq 0$, which implies that
$x + \alpha \beta_i \not\in J_\beta$.
That is, in view of Definition \ref{def:Lambda} $(ii)$, we
conclude that $i \not\in \Lambda_\beta(J_\beta)$.

We now assume that $i \not\in \Lambda_\beta(J_\beta)$.
Then by Definition \ref{def:Lambda} $(ii)$, it holds that
\begin{eqnarray} \label{eq:2022-4-28-1}
x + \alpha \beta_i \not\in J_\beta
\end{eqnarray}
for any $x \in J_\beta$ and $\alpha \neq 0$.
Using Definition \ref{def:betabasic} again, we have
\begin{eqnarray} \label{eq:2022-4-28-2}
x + \alpha \beta_i
= \sum_{j \not\in \Lambda_4} \alpha_j \beta_j +
  \sum_{j \in \Lambda_4} a_j p_{1j} \beta_j +
  \alpha \beta_i
\end{eqnarray}
for all $x \in J_\beta$ and $\alpha \neq 0$.
We assume on the contrary that $i \not\in \Lambda_4$.
In view of (\ref{eq:2022-4-28-2}) and by the structure of
$J_i$ ($a_i J_i$ is indeed a non-degenerate interval for
$i \not\in \Lambda_4$), it holds that
\begin{eqnarray*}
x + \alpha \beta_i
=   \sum_{j \not\in \Lambda_4 \cup \{ i \}} \alpha_j \beta_j +
    (\alpha_i + \alpha) \beta_i +
    \sum_{j \in \Lambda_4} a_j p_{1j} \beta_j
\in J_\beta
\end{eqnarray*}
for some $x \in J_\beta$ and $\alpha \neq 0$, which is contrary
to (\ref{eq:2022-4-28-1}).
(We note that, for each $i \not\in \Lambda_4$,
$\alpha_i \in a_i J_i$ and there exists a real number
$\alpha \neq 0$ satisfying $\alpha_i + \alpha \in a_i J_i$.)
Therefore, we conclude that if
$i \not\in \Lambda_\beta(J_\beta)$, then $i \in \Lambda_4$.
\hfill$\Box$


\begin{thm} \label{thm:set9}
Let $\beta = \{ \beta_i \}_{i \in \mathbb{N}}$ be a complete
orthonormal sequence in $M_a$ and let $J_\beta$ be either a
translation of a $\beta$-basic cylinder or a translation of a
$\beta$-basic interval and $p \in J_\beta$.
Then
\begin{eqnarray*}
{\rm GS}(J_\beta,p)
= \Bigg\{ p + \sum_{i \in \Lambda_\beta(J_\beta)}
          \alpha_i \beta_i \in M_a :
          \alpha_i \in \mathbb{R} ~\mbox{for all}~
          i \in \Lambda_\beta(J_\beta)
  \Bigg\}.
\end{eqnarray*}
\end{thm}

\noindent
\emph{Proof}.
Assume that $x$ is an arbitrary element of
${\rm GS}(J_\beta,p)$.
By Definition \ref{def:gls}, we have
\begin{eqnarray*}
x - p =   \sum_{i=1}^m \sum_{j=1}^\infty
          \varepsilon_{ij} (x_{ij} - p)
      \in M_a
\end{eqnarray*}
for some $m \in \mathbb{N}$, $\varepsilon_{ij} \in \mathbb{R}$,
and $x_{ij} \in J_\beta$.
Furthermore, since $x_{ij}, p \in J_\beta$, by Definition
\ref{def:betabasic}, we get
\begin{eqnarray*}
x_{ij} = \sum_{k=1}^\infty \gamma_k \beta_k
       = \sum_{k \in \mathbb{N} \setminus \Lambda_4}
         \gamma_k \beta_k +
         \sum_{k \in \Lambda_4} a_k p_{1k} \beta_k
\end{eqnarray*}
and
\begin{eqnarray*}
p = \sum_{k=1}^\infty \delta_k \beta_k
  = \sum_{k \in \mathbb{N} \setminus \Lambda_4}
    \delta_k \beta_k +
    \sum_{k \in \Lambda_4} a_k p_{1k} \beta_k
\end{eqnarray*}
for some $\gamma_k, \delta_k \in a_k J_k$.

Since $\{ \beta_i \}_{i \in \mathbb{N}}$ is a complete
orthonormal sequence in $M_a$, it follows from Definition
\ref{def:betabasic} and Remark \ref{rem:6.2} that
\begin{eqnarray*}
\begin{split}
x - p & = \sum_{i=1}^m \sum_{j=1}^\infty
          \varepsilon_{ij} (x_{ij} - p) \\
      & = \sum_{i=1}^m \sum_{j=1}^\infty \varepsilon_{ij}
          \sum_{k \in \mathbb{N} \setminus \Lambda_4}
          (\gamma_k - \delta_k) \beta_k \\
      & = \sum_{k \in \mathbb{N} \setminus \Lambda_4}
          \omega_k \beta_k
        = \sum_{i \in \Lambda_\beta(J_\beta)} \omega_i \beta_i
\end{split}
\end{eqnarray*}
for some real numbers $\omega_i$.
Hence, we see that $x \in {\rm GS}(J_\beta,p) \subset M_a$ and
\begin{eqnarray*}
\begin{split}
x & =   p + \sum_{i \in \Lambda_\beta(J_\beta)} \omega_i \beta_i
        ~(\in M_a) \\
  & \in \Bigg\{ p + \sum_{i \in \Lambda_\beta(J_\beta)}
                \alpha_i \beta_i \in M_a :
                \alpha_i \in \mathbb{R}
                ~\mbox{for all}~ i \in \Lambda_\beta(J_\beta)
        \Bigg\},
\end{split}
\end{eqnarray*}
which implies that
\begin{eqnarray*}
\hspace*{3mm}
{\rm GS}(J_\beta,p)
\subset \Bigg\{ p + \sum_{i \in \Lambda_\beta(J_\beta)}
                \alpha_i \beta_i \in M_a :
                \alpha_i \in \mathbb{R}
                ~\mbox{for all}~ i \in \Lambda_\beta(J_\beta)
        \Bigg\}.
\end{eqnarray*}

It remains to prove the reverse inclusion.
According to the structure of $J_\beta$ given in Definition
\ref{def:betabasic}, for each $i \in \Lambda_\beta(J_\beta)$,
there exists a real number $\gamma_i \neq 0$ such that
$p + \gamma_i \beta_i \in J_\beta$.
In other words, for each $i \in \Lambda_\beta(J_\beta)$, there
exists a $u_i \in J_\beta$ such that
$\gamma_i \beta_i = u_i - p$.
Thus, if we assume that
\begin{eqnarray*}
p + \sum_{i \in \Lambda_\beta(J_\beta)} \alpha_i \beta_i \in M_a
\end{eqnarray*}
for some $\alpha_i \in \mathbb{R}$, then
\begin{eqnarray*}
p + \sum_{i \in \Lambda_\beta(J_\beta)} \alpha_i \beta_i
=   p + \sum_{i \in \Lambda_\beta(J_\beta)}
    \frac{\alpha_i}{\gamma_i}
    (\gamma_i \beta_i)
=   p + \sum_{i \in \Lambda_\beta(J_\beta)}
    \frac{\alpha_i}{\gamma_i} (u_i - p)
\in {\rm GS}(J_\beta,p),
\end{eqnarray*}
since $u_i \in J_\beta$ for all $i \in \Lambda_\beta(J_\beta)$,
which implies that
\begin{eqnarray*}
\hspace*{3mm}
{\rm GS}(J_\beta,p)
\supset \Bigg\{ p + \sum_{i \in \Lambda_\beta(J_\beta)}
                \alpha_i \beta_i \in M_a :
                \alpha_i \in \mathbb{R}
                ~\mbox{for all}~ i \in \Lambda_\beta(J_\beta)
        \Bigg\}.
\end{eqnarray*}
We end the proof in this way.
\hfill$\Box$
\vspace*{5mm}

In the following theorem, we introduce an interesting inclusion
property of the second-order generalized linear span.


\begin{thm} \label{thm:setten}
Assume that a bounded subset $E$ of $M_a$ contains at least
two elements and $E \subset H$, where $H$ is a closed subspace
of $M_a$.
Let $\beta = \{ \beta_i \}_{i \in \mathbb{N}}$ and
$\{ \beta_i \}_{i \in \Lambda}$ be complete orthonormal
sequences in the Hilbert spaces $M_a$ and $H$, respectively,
and let $J_\beta$ be either a translation of a $\beta$-basic
cylinder or a translation of a $\beta$-basic interval that
possesses the following properties:
\begin{itemize}
\item[$(i)$]  There exists a $p \in J_\beta \cap E$;
\item[$(ii)$] $\Lambda_\beta(J_\beta) \subset \Lambda$.
\end{itemize}
Then $J_\beta \subset {\rm GS}(J_\beta,p) \subset H$.
\end{thm}

\noindent
\emph{Proof}.
Assume that $x \in {\rm GS}(J_\beta,p)$.
Then, according to Theorem \ref{thm:set9}, there exist some
real numbers $\alpha_i$ that satisfy
\begin{eqnarray} \label{eq:2021-10-3-1}
x =   p + \sum_{i \in \Lambda_\beta(J_\beta)} \alpha_i \beta_i
  \in M_a.
\end{eqnarray}
Assume that $i \in \Lambda_\beta(J_\beta)$.
Then, $i \in \Lambda$ by $(ii)$.
Since $H$ is a real vector space and $\beta_i \in H$, it holds
that
\begin{eqnarray} \label{eq:2021-10-3-2}
\alpha_i \beta_i \in H
\end{eqnarray}
for all $i \in \Lambda_\beta(J_\beta)$.

Now we define
\begin{eqnarray*}
x_n := p + \sum_{i \in \Lambda_n} \alpha_i \beta_i
\end{eqnarray*}
for any $n \in \mathbb{N}$, where we set
$\Lambda_n = \{ i \in \Lambda_\beta(J_\beta) : i < n \}$.
Since $H$ is a real vector space, it follows from
(\ref{eq:2021-10-3-2}) that
\begin{eqnarray*}
\sum_{i \in \Lambda_n} \alpha_i \beta_i \in H
\end{eqnarray*}
for each $n \in \mathbb{N}$.
Furthermore, our assumption that $p \in E \subset H$ yields
\begin{eqnarray*}
x_n = p + \sum_{i \in \Lambda_n} \alpha_i \beta_i \in H
\end{eqnarray*}
for each $n \in \mathbb{N}$.
Since $H$ is closed, it follows from (\ref{eq:2021-10-3-1})
that
\begin{eqnarray*}
x = \lim_{n \to \infty} x_n \in H,
\end{eqnarray*}
which implies that ${\rm GS}(J_\beta,p) \subset H$.
\hfill$\Box$
\vspace*{5mm}

Since in some ways index sets have some properties of dimensions
in vector space, the following theorem may seem obvious.


\begin{thm} \label{thm:setsib}
Assume that a bounded subset $E$ of $M_a$ contains at least
two elements and $p \in E$.
Then, $\Lambda({\rm GS}^2(E,p)) = \mathbb{N}$ if and only if\/
${\rm GS}^2(E,p) = M_a$.
\end{thm}

\noindent
\emph{Proof}.
Let $x$ be an arbitrary element of $M_a$.
According to Remark \ref{rem:new2.1}, there exist some real
numbers $\alpha_i$ such that
\begin{eqnarray} \label{eq:2021-10-9}
x = \sum\limits_{i=1}^\infty \alpha_i e_i \in M_a.
\end{eqnarray}

If $\Lambda({\rm GS}^2(E,p)) = \mathbb{N}$, then it follows
from Lemma \ref{lem:setahob} that
\begin{eqnarray*}
p + \alpha_i e_i \in {\rm GS}^2(E,p)
\end{eqnarray*}
for all $i \in \mathbb{N}$.
In other words,
\begin{eqnarray*}
\alpha_i e_i \in {\rm GS}^2(E,p) - p
\end{eqnarray*}
for all $i \in \mathbb{N}$.

By Lemma \ref{lem:3.8} $(i)$, we get
\begin{eqnarray*}
x_n := \sum_{i=1}^n \alpha_i e_i \in {\rm GS}^2(E,p) - p
\end{eqnarray*}
for any $n \in \mathbb{N}$.
Due to Lemma \ref{lem:3.8} $(iii)$ and (\ref{eq:2021-10-9}),
we further obtain
\begin{eqnarray*}
x =   \sum_{i=1}^\infty \alpha_i e_i
  =   \lim_{n \to \infty} x_n
  \in {\rm GS}^2(E,p) - p,
\end{eqnarray*}
which implies that $M_a \subset {\rm GS}^2(E,p) - p$, or
equivalently, $M_a \subset {\rm GS}^2(E,p)$.

The reverse inclusion is trivial.
\hfill$\Box$


\section{Second-order extension of isometries}

It was proved in Theorem \ref{thm:new3} that the domain of a
$d_a$-isometry $f : E_1 \to E_2$ can be extended to the
first-order generalized linear span ${\rm GS}(E_1,p)$ whenever
$E_1$ is a nonempty bounded subset of $M_a$, whether degenerate
or non-degenerate.

Now we generalize Theorem \ref{thm:new3} in the following
theorem.
More precisely, we prove that the domain of $f$ can be extended
to its second-order generalized linear span ${\rm GS}^2(E_1,p)$.
We note that ${\rm GS}^2(E_1,p) = \overline{{\rm GS}(E_1,p)}$
by Lemma \ref{lem:3.8} $(iii)$.
Therefore, Theorem \ref{thm:setahob} is a further generalization
of \cite[Theorem 2.2]{jung1}.

In the proof, we use the fact that ${\rm GS}^n(E_1,p) - p$ is
a real vector space.


\begin{thm} \label{thm:setahob}
Let $E_1$ be a bounded subset of $M_a$ that is $d_a$-isometric
to a subset $E_2$ of $M_a$ via a surjective $d_a$-isometry
$f : E_1 \to E_2$.
Assume that $p$ and $q$ are elements of $E_1$ and $E_2$, which
satisfy $q = f(p)$.
The function $F_2 : {\rm GS}^2(E_1,p) \to M_a$ is a
$d_a$-isometry and the function
$T_{-q} \circ F_2 \circ T_p : {\rm GS}^2(E_1,p) - p \to M_a$
is linear.
In particular, $F_2$ is an extension of $F$.
\end{thm}

\noindent
\emph{Proof}.
$(a)$
Suppose $r$ is a positive real number satisfying
$E_1 \subset B_r(p)$.
Referring to the changes presented in the table below
and following the first part of proof of Theorem \ref{thm:new3},
we can easily prove that $F_2$ is a $d_a$-isometry.

\begin{center}
\begin{tabular}{ccccccc}
\hline\hline
Theorem \ref{thm:new3}: & $E_1$                         &
${\rm GS}(E_1,p)$   & $f$ & $F$   &
Definition \ref{def:ext}        & Lemma \ref{lem:new2} \\
\hline
Here:                   & ${\rm GS}(E_1,p) \cap B_r(p)$ &
${\rm GS}^2(E_1,p)$ & $F$ & $F_2$ &
Definition \ref{def:6.1} & Lemma \ref{lem:6.1}  \\
\hline\hline
\end{tabular}
\end{center}

$(b)$
Referring to the changes presented in the table below and
following $(d)$ of the proof of Theorem \ref{thm:new3}, we can
prove the linearity of
$T_{-q} \circ F_n \circ T_p : {\rm GS}^n(E_1,p) - p \to M_a$
in the more general setting for $n \geq 2$.

\begin{center}
\begin{tabular}{cccc}
\hline\hline
Theorem \ref{thm:new3}: & ${\rm GS}(E_1,p)$   & $F$
& (\ref{eq:2018-9-28}) \\
\hline
Here:                   & ${\rm GS}^n(E_1,p)$ & $F_n$
& Lemma \ref{lem:6.1}  \\
\hline\hline
\end{tabular}
\end{center}
\vspace*{2mm}

$(c)$ According to Definition \ref{def:6.1} $(i)$, for any
$m \in \mathbb{N}$, $x_{ij} \in {\rm GS}(E_1,p) \cap B_r(p)$,
and any $\alpha_{ij} \in \mathbb{R}$ with
$\sum\limits_{i=1}^m \sum\limits_{j=1}^\infty \alpha_{ij}
 (x_{ij} - p) \in M_a$, there exists a
$u \in {\rm GS}^2(E_1,p)$ satisfying
\begin{eqnarray} \label{eq:3.6a}
u - p = \sum_{i=1}^m \sum_{j=1}^\infty
        \alpha_{ij} (x_{ij} - p) \in M_a.
\end{eqnarray}
Due to Definition \ref{def:6.1} $(ii)$, we further have
\begin{eqnarray} \label{eq:3.7a}
(T_{-q} \circ F_2 \circ T_p)(u - p)
= \sum_{i=1}^m \sum_{j=1}^\infty
  \alpha_{ij} (T_{-q} \circ F \circ T_p)(x_{ij} - p).
\end{eqnarray}

If we set $\alpha_{11} = 1$, $\alpha_{ij} = 0$ for each
$(i,j) \neq (1,1)$, and $x_{11} = x$ in (\ref{eq:3.6a}) and
(\ref{eq:3.7a}) to see
\begin{eqnarray} \label{eq:2021-9-26}
(T_{-q} \circ F_2 \circ T_p)(x-p)
= (T_{-q} \circ F \circ T_p)(x-p)
\end{eqnarray}
for all $x \in {\rm GS}(E_1,p) \cap B_r(p)$.

Let $w$ be an arbitrary element of ${\rm GS}(E_1,p)$.
Then, $w - p \in {\rm GS}(E_1,p) - p$.
Since ${\rm GS}(E_1,p) - p$ is a real vector space and
$B_r(p) - p = B_r(0)$, there exists a real number $\mu \neq 0$
such that
\begin{eqnarray*}
\mu (w - p) \in ({\rm GS}(E_1,p) - p) \cap (B_r(p) - p).
\end{eqnarray*}
Hence, by (\ref{eq:2021-11-2}), we can choose a
$v \in {\rm GS}(E_1,p) \cap B_r(p)$ such that
$\mu (w - p) = v - p$.
Since both $T_{-q} \circ F_2 \circ T_p$ and
$T_{-q} \circ F \circ T_p$ are linear and
${\rm GS}(E_1,p) \subset {\rm GS}^2(E_1,p)$, it follows from
(\ref{eq:2021-9-26}) that
\begin{eqnarray*}
\begin{split}
\mu (T_{-q} \circ F_2 \circ T_p)(w - p)
& = (T_{-q} \circ F_2 \circ T_p)(\mu (w - p)) \\
& = (T_{-q} \circ F_2 \circ T_p)(v - p) \\
& = (T_{-q} \circ F \circ T_p)(v - p) \\
& = (T_{-q} \circ F \circ T_p)(\mu (w - p)) \\
& = \mu (T_{-q} \circ F \circ T_p)(w - p).
\end{split}
\end{eqnarray*}
Therefore, it follows that
$(T_{-q} \circ F_2 \circ T_p)(w-p)
 = (T_{-q} \circ F \circ T_p)(w-p)$ for all
$w \in {\rm GS}(E_1,p)$, \textit{i.e.}, $F_2(w) = F(w)$ for
all $w \in {\rm GS}(E_1,p)$.
In other words, $F_2$ is an extension of $F$.
Also, because of Theorem \ref{thm:new3}, we see that $F_2$ is
obviously an extension of $f$.
\hfill$\Box$
\vspace*{5mm}

On account of Proposition \ref{rem:new3.1}, it holds that
\begin{eqnarray*}
{\rm GS}^2(E_1,p) = \cdots = {\rm GS}^{n-1}(E_1,p)
= {\rm GS}^n(E_1,p)
\end{eqnarray*}
for every integer $n \geq 3$.
According to this formula, the assertion of the following
theorem seems obvious, but since the proof is not long, we
introduce the proof here.


\begin{thm} \label{thm:setyol}
Let $E_1$ be a bounded subset of $M_a$ that is $d_a$-isometric
to a subset $E_2$ of $M_a$ via a surjective $d_a$-isometry
$f : E_1 \to E_2$.
Assume that $p$ and $q$ are elements of $E_1$ and $E_2$, which
satisfy $q = f(p)$.
Then $F_n$ is identically the same as $F_2$ for any integer
$n \geq 3$, where $F_2$ and $F_n$ are defined in Definition
\ref{def:6.1}.
\end{thm}

\noindent
\emph{Proof}.
Let $r$ be a fixed positive real number satisfying
$E_1 \subset B_r(p)$.
We assume that $F_2 \equiv F_3 \equiv \cdots \equiv F_{n-1}$
on ${\rm GS}^2(E_1,p)$.
Let $x$ be an arbitrary element of ${\rm GS}^n(E_1,p)$.
Then, in view of (\ref{eq:2021-11-2}), there exist a real
number $\mu \neq 0$ and an element $u$ of
${\rm GS}^n(E_1,p) \cap B_r(p)$ such that
\begin{eqnarray*}
u - p = \mu (x - p)
\in ({\rm GS}^n(E_1,p) - p) \cap (B_r(p) - p).
\end{eqnarray*}
If we put $\alpha_{11} = 1$, $\alpha_{ij} = 0$ for all
$(i,j) \neq (1,1)$, and $x_{11} = v$ in Definition \ref{def:6.1}
$(ii)$, then we get
\begin{eqnarray} \label{eq:2021-10-29}
(T_{-q} \circ F_n \circ T_p)(v - p)
= (T_{-q} \circ F_{n-1} \circ T_p)(v - p)
\end{eqnarray}
for all
$v \in {\rm GS}^{n-1}(E_1,p) \cap B_r(p)
   =   {\rm GS}^n(E_1,p) \cap B_r(p)$.
We note by Proposition \ref{rem:new3.1} that
${\rm GS}^n(E_1,p) = {\rm GS}^{n-1}(E_1,p) = \cdots =
 {\rm GS}^2(E_1,p)$.

Since $T_{-q} \circ F_n \circ T_p$ is linear by $(b)$ in the
proof of Theorem \ref{thm:setahob}, it follows from
(\ref{eq:2021-10-29}) that
\begin{eqnarray*}
\begin{split}
\mu (T_{-q} \circ F_n \circ T_p)(x - p)
& = (T_{-q} \circ F_n \circ T_p)(u - p) \\
& = (T_{-q} \circ F_{n-1} \circ T_p)(u - p) \\
& = (T_{-q} \circ F_2 \circ T_p)(u - p) \\
& = \mu (T_{-q} \circ F_2 \circ T_p)(x - p),
\end{split}
\end{eqnarray*}
\textit{i.e.}, $F_n(x) = F_2(x)$ for every
$x \in {\rm GS}^n(E_1,p) = {\rm GS}^2(E_1,p)$.
By mathematical induction, we conclude that $F_n$ is
identically the same as $F_2$ for every integer $n \geq 3$.
\hfill$\Box$
\vspace*{5mm}

Assume that $J$ is either a translation of a basic cylinder or
a translation of a basic interval, and $p$ is an element of
$J$.
Due to Definition \ref{def:basic}, Remark \ref{rem:6.1}, and
Theorem \ref{thm:set9}, ${\rm GS}(J,p)$ is a closed subset of
$M_a$.


\begin{rem} \label{rem:7.1}
${\rm GS}(J,p)$ is a closed subset of $M_a$.
\end{rem}

\noindent
\emph{Proof}.
Assume that $p = (p_1, p_2, \ldots, p_i, \ldots)$ is a fixed
element of $J$, where $J$ is a translation of a basic cylinder
or a translation of a basic interval.
In view of Definition \ref{def:gls} and Remark \ref{rem:6.1}
$(iii)$, we note that $x_i = p_i$ for each
$x = (x_1, x_2, \ldots, x_i, \ldots) \in {\rm GS}(J,p)$ and
each $i \not\in \Lambda(J)$.

Assume that $\{ z_n \}_{n \in \mathbb{N}}$ is a sequence of
elements in ${\rm GS}(J,p)$, which converges to an element
$z = (z_1, z_2, \ldots, z_i, \ldots)$ of $M_a$.
Let us denote by $z_{ni}$ the $i$th component of $z_n$ for any
$i, n \in \mathbb{N}$.
Since $z_n \in {\rm GS}(J,p)$ for every $n \in \mathbb{N}$,
the previous argument implies that $z_{ni} = p_i$ for each
$i \not\in \Lambda(J)$.
Thus, we conclude that $z_i = p_i$ for each
$i \not\in \Lambda(J)$.
This fact, together with Theorem \ref{thm:set9}, implies that
$z \in {\rm GS}(J,p)$.
Therefore, we conclude that ${\rm GS}(J,p)$ is a closed subset
of $M_a$.
\hfill$\Box$
\vspace*{5mm}

We note that $\{ \frac{1}{a_i} e_i \}_{i \in \mathbb{N}}$ is a
complete orthonormal sequence in $M_a$.
On account of Theorem \ref{thm:set9}, we notice that
$\Lambda(J) = \Lambda({\rm GS}(J,p))$.


\begin{rem} \label{rem:7.2}
${\rm GS}^2(J,p) = {\rm GS}(J,p)$.
\end{rem}

\noindent
\emph{Proof}.
Referring to the changes presented in the table below

\begin{center}
\begin{tabular}{cccccc}
\hline\hline
Proposition \ref{rem:new3.1}: & ${\rm GS}(E,p) \cap B_r(p)$ &
${\rm GS}^2(E,p)$ & ${\rm GS}^3(E,p)$ & $x$ & $x_m$ \\
\hline
Here:                         & $J$                         &
${\rm GS}(J,p)$   & ${\rm GS}^2(J,p)$ & $u$ & $u_m$ \\
\hline\hline
\end{tabular}
\end{center}
\vspace*{2mm}

\noindent
and following the part $(a)$ in the proof of Proposition
\ref{rem:new3.1}, we can easily show that
${\rm GS}^2(J,p) = {\rm GS}(J,p)$.
\hfill$\Box$
\vspace*{5mm}

Hence, by Theorem \ref{thm:set9} and Remark \ref{rem:7.2}, we
have
\begin{eqnarray} \label{eq:2021-10-21-2}
\begin{split}
u - p
& = \sum_{i=1}^\infty
    \left\langle u - p,\, \frac{1}{a_i} e_i \right\rangle_a
    \frac{1}{a_i} e_i \\
& = \sum_{i=1}^\infty a_i (u_i - p_i) \frac{1}{a_i} e_i \\
& = \sum_{i \in \Lambda(J)}
    a_i (u_i - p_i) \frac{1}{a_i} e_i \\
& = \sum_{i \in \Lambda(J)}
    \left\langle u - p,\, \frac{1}{a_i} e_i \right\rangle_a
    \frac{1}{a_i} e_i
\end{split}
\end{eqnarray}
for all $u \in {\rm GS}^2(J,p) = {\rm GS}^n(J,p)$, where
$n \in \mathbb{N}$.

Using a similar approach to the proof of
\cite[Theorem 2.4]{jungkim}, we can apply Lemma \ref{lem:6.1}
to prove the following theorem.


\begin{thm} \label{thm:new4}
Assume that $J$ is either a translation of a basic cylinder or
a translation of a basic interval, $K$ is a subset of $M_a$,
and that there exists a surjective $d_a$-isometry $f : J \to K$.
Suppose $p$ is an element of $J$ and $q$ is an element of $K$
with $q = f(p)$.
For any $n \in \mathbb{N}$, the $d_a$-isometry
$F_n : {\rm GS}^n(J,p) \to M_a$ given in Definition
\ref{def:6.1} satisfies
\begin{eqnarray*}
(T_{-q} \circ F_n \circ T_p)(u-p)
= \sum_{i \in \Lambda(J)}
  \bigg\langle u-p, \frac{1}{a_i} e_i \bigg\rangle_a
  \frac{1}{a_i} (T_{-q} \circ F_n \circ T_p)(e_i)
\end{eqnarray*}
for all $u \in {\rm GS}^n(J,p)$.
\end{thm}

\noindent
\emph{Proof}.
Since $p + e_i \in {\rm GS}^n(J,p)$ for each $i \in \Lambda(J)$,
it follows from Lemma \ref{lem:6.1} that
\begin{eqnarray}
\begin{split}
& \Bigg\langle (T_{-q} \circ F_n \circ T_p)(u-p) -
               \sum_{i \in \Lambda(J)}
               \left\langle u-p, \frac{1}{a_i} e_i \right\rangle_a
               \frac{1}{a_i}
               (T_{-q} \circ F_n \circ T_p)(e_i), \\
& \hspace*{3mm}
    \left.     (T_{-q} \circ F_n \circ T_p)(u-p) -
               \sum_{j \in \Lambda(J)}
               \left\langle u-p, \frac{1}{a_j} e_j \right\rangle_a
               \frac{1}{a_j}
               (T_{-q} \circ F_n \circ T_p)(e_j)
  \right\rangle_a \\
& \hspace*{3mm}
  = \big\langle (T_{-q} \circ F_n \circ T_p)(u-p),\,
                (T_{-q} \circ F_n \circ T_p)(u-p)
    \big\rangle_a \\
& \hspace*{8mm}
    -\,\sum_{j \in \Lambda(J)}
    \left\langle u-p, \frac{1}{a_j} e_j \right\rangle_a
    \frac{1}{a_j}
    \big\langle (T_{-q} \circ F_n \circ T_p)(u-p),\,
                (T_{-q} \circ F_n \circ T_p)(e_j)
    \big\rangle_a \\
& \hspace*{8mm}
    -\,\sum_{i \in \Lambda(J)}
    \left\langle u-p, \frac{1}{a_i} e_i \right\rangle_a
    \frac{1}{a_i}
    \big\langle (T_{-q} \circ F_n \circ T_p)(e_i),\,
                (T_{-q} \circ F_n \circ T_p)(u-p)
    \big\rangle_a \\
& \hspace*{8mm}
    +\,\sum_{i \in \Lambda(J)} \sum_{j \in \Lambda(J)}
    \left\langle u-p, \frac{1}{a_i} e_i \right\rangle_a
    \left\langle u-p, \frac{1}{a_j} e_j \right\rangle_a
    \times
    \label{eq:3.3} \\
& \hspace*{34mm} \times
    \frac{1}{a_i a_j}
    \big\langle (T_{-q} \circ F_n \circ T_p)(e_i),\,
                (T_{-q} \circ F_n \circ T_p)(e_j)
    \big\rangle_a \\
& \hspace*{3mm}
  = \big\langle u-p, u-p \big\rangle_a -
    \sum_{j \in \Lambda(J)}
    \left\langle u-p, \frac{1}{a_j} e_j \right\rangle_a
    \left\langle u-p, \frac{1}{a_j} e_j \right\rangle_a \\
& \hspace*{8mm}
    -\,\sum_{i \in \Lambda(J)}
    \left\langle u-p, \frac{1}{a_i} e_i \right\rangle_a
    \left\langle \frac{1}{a_i} e_i, u-p \right\rangle_a \\
& \hspace*{8mm}
    +\,\sum_{i \in \Lambda(J)} \sum_{j \in \Lambda(J)}
    \left\langle u-p, \frac{1}{a_i} e_i \right\rangle_a
    \left\langle u-p, \frac{1}{a_j} e_j \right\rangle_a
    \left\langle \frac{1}{a_i} e_i,\, \frac{1}{a_j} e_j
    \right\rangle_a \\
& \hspace*{3mm}
  = \big\langle u-p, u-p \big\rangle_a -
    \sum_{j \in \Lambda(J)}
    \left\langle u-p, \frac{1}{a_j} e_j \right\rangle_a
    \left\langle u-p, \frac{1}{a_j} e_j \right\rangle_a
\end{split}
\end{eqnarray}
for all $u \in {\rm GS}^n(J,p)$, since
$\{ \frac{1}{a_i} e_i \}_{i \in \mathbb{N}}$ is an orthonormal
sequence in $M_a$.

Furthermore, we note that each $u \in {\rm GS}^n(J,p)$ has the
expression given in (\ref{eq:2021-10-21-2}).
Hence, if we replace $u-p$ in (\ref{eq:3.3}) with the expression
(\ref{eq:2021-10-21-2}), then we have
\begin{eqnarray*}
\Bigg\| (T_{-q} \circ F_n \circ T_p)(u-p) -
        \sum_{i \in \Lambda(J)}
        \left\langle u-p, \frac{1}{a_i} e_i \right\rangle_a
        \frac{1}{a_i} (T_{-q} \circ F_n \circ T_p)(e_i)
\Bigg\|_a^2 = 0
\end{eqnarray*}
for all $u \in {\rm GS}^n(J,p)$, which implies the validity of our
assertion.
\hfill$\Box$
\vspace*{5mm}

According to the following theorem, the image of the
first-order generalized linear span of $E_1$ with respect to
$p$ under the $d_a$-isometry $F$ is just the first-order
generalized linear span of $F(E_1)$ with respect to $F(p)$.
This assertion holds also for the second-order generalized
linear span and $F_2$.
According to Proposition \ref{rem:new3.1} and Theorem
\ref{thm:setyol}, the argument of the following theorem only
makes sense when $n = 1$ or $2$.


\begin{thm} \label{thm:2018-11-7}
Assume that $E_1$ and $E_2$ are bounded subsets of $M_a$ that
are $d_a$-isometric to each other via a surjective
$d_a$-isometry $f : E_1 \to E_2$.
Suppose $p$ is an element of $E_1$ and $q$ is an element of
$E_2$ with $q = f(p)$.
If $F_n : {\rm GS}^n(E_1,p) \to M_a$ is the extension of $f$
defined in Definition \ref{def:6.1}, then
${\rm GS}^n(E_2,q) = F_n({\rm GS}^n(E_1,p))$ for every
$n \in \mathbb{N}$.
\end{thm}

\noindent
\emph{Proof}.
$(a)$
First, we prove that our assertion is true for $n = 1$,
\textit{i.e.}, we prove that
${\rm GS}(E_2,q) = F({\rm GS}(E_1,p))$.
Let $r$ be a fixed positive real number satisfying
$E_1 \subset B_r(p)$.

$(b)$
Due to Definition \ref{def:gls}, for any
$y \in F({\rm GS}(E_1,p))$, there exists an element
$x \in {\rm GS}(E_1,p)$ with
\begin{eqnarray*}
y = F(x)
  = F \bigg( p + \sum_{i=1}^m \sum_{j=1}^\infty
             \alpha_{ij} \big( u_{ij} - p \big) \bigg)
\end{eqnarray*}
for some $m \in \mathbb{N}$, $u_{ij} \in E_1 \cap B_r(p)$, and
some $\alpha_{ij} \in \mathbb{R}$ with
$x = p + \sum\limits_{i=1}^m \sum\limits_{j=1}^\infty
 \alpha_{ij} (u_{ij} - p) \in M_a$.

On the other hand, by Definition \ref{def:ext}, we have
\begin{eqnarray*}
(T_{-q} \circ F \circ T_p)
\bigg( \sum_{i=1}^m \sum_{j=1}^\infty
       \alpha_{ij} \big( u_{ij} - p \big)
\bigg)
= \sum_{i=1}^m \sum_{j=1}^\infty
  \alpha_{ij} (T_{-q} \circ f \circ T_p) (u_{ij} - p)
\end{eqnarray*}
which is equivalent to
\begin{eqnarray*}
F(x) - q
= F \bigg( p + \sum_{i=1}^m \sum_{j=1}^\infty
           \alpha_{ij} \big( u_{ij} - p \big)
    \bigg) - q
= \sum_{i=1}^m \sum_{j=1}^\infty
  \alpha_{ij} \big( f(u_{ij}) - q \big).
\end{eqnarray*}
Since $u_{ij} \in E_1$ for all $i$ and $j$, it holds that
$f(u_{ij}) \in f(E_1) = E_2$ for each $i$ and $j$.
Moreover, since $u_{ij} \in E_1 \cap B_r(p)$ for all $i$ and
$j$, it follows from Lemma \ref{lem:new2} that
\begin{eqnarray*}
\begin{split}
\| f(u_{ij}) - q \|_a^2
& = \| (T_{-q} \circ f \circ T_p)(u_{ij} - p) \|_a^2 \\
& = \big\langle (T_{-q} \circ f \circ T_p)(u_{ij} - p),\,
                (T_{-q} \circ f \circ T_p)(u_{ij} - p)
    \big\rangle_a \\
& = \langle u_{ij} - p, u_{ij} - p \rangle_a \\
& = \| u_{ij} - p \|_a^2 \\
& < r^2
\end{split}
\end{eqnarray*}
for all $i$ and $j$.
Hence, $f(u_{ij}) \in E_2 \cap B_r(q)$ for all $i$ and $j$.

Furthermore, it follows from Lemma \ref{lem:new2} that
\begin{eqnarray*}
\begin{split}
& \left\| \sum_{i=1}^m \sum_{j=1}^\infty
          \alpha_{ij} \big( f(u_{ij}) - q \big)
  \right\|_a^2 \\
& \hspace*{3mm}
  = \left\| \sum_{i=1}^m \sum_{j=1}^\infty
            \alpha_{ij} (T_{-q} \circ f \circ T_p)(u_{ij} - p)
    \right\|_a^2 \\
& \hspace*{3mm}
  = \left\langle
    \sum_{i=1}^m \sum_{j=1}^\infty
    \alpha_{ij} (T_{-q} \circ f \circ T_p)(u_{ij} - p),\,
    \sum_{k=1}^m \sum_{\ell=1}^\infty
    \alpha_{k\ell} (T_{-q} \circ f \circ T_p)(u_{k\ell} - p)
    \right\rangle_a \\
& \hspace*{3mm}
  = \sum_{i=1}^m \sum_{k=1}^m
    \sum_{j=1}^\infty \alpha_{ij}
    \sum_{\ell=1}^\infty \alpha_{k\ell}
    \big\langle (T_{-q} \circ f \circ T_p)(u_{ij} - p),\,
                (T_{-q} \circ f \circ T_p)(u_{k\ell} - p)
    \big\rangle_a \\
& \hspace*{3mm}
  = \sum_{i=1}^m \sum_{k=1}^m
    \sum_{j=1}^\infty \alpha_{ij}
    \sum_{\ell=1}^\infty \alpha_{k\ell}
    \langle u_{ij} - p, u_{k\ell} - p \rangle_a \\
& \hspace*{3mm}
  = \left\langle
    \sum_{i=1}^m \sum_{j=1}^\infty \alpha_{ij} (u_{ij} - p),\,
    \sum_{k=1}^m \sum_{\ell=1}^\infty
    \alpha_{k\ell} (u_{k\ell} - p)
    \right\rangle_a \\
& \hspace*{3mm}
  = \left\|
    \sum_{i=1}^m \sum_{j=1}^\infty \alpha_{ij} (u_{ij} - p)
    \right\|_a^2 \\
& \hspace*{3mm}
  < \infty,
\end{split}
\end{eqnarray*}
since
$\sum\limits_{i=1}^m \sum\limits_{j=1}^\infty
 \alpha_{ij} (u_{ij} - p) = x - p \in M_a$.
Thus, on account of Remark \ref{rem:new2.1}, we see that
$\sum\limits_{i=1}^m \sum\limits_{j=1}^\infty
 \alpha_{ij} (f(u_{ij}) - q) \in M_a$.
Therefore, in view of Definition \ref{def:gls}, we get
\begin{eqnarray*}
y =   F(x)
  =   q + \sum_{i=1}^m \sum_{j=1}^\infty
      \alpha_{ij} \big( f(u_{ij}) - q \big)
  \in {\rm GS}(E_2,q)
\end{eqnarray*}
and we conclude that
$F({\rm GS}(E_1,p)) \subset {\rm GS}(E_2,q)$.

$(c)$ Now we assume that $y \in {\rm GS}(E_2,q)$.
By Definition \ref{def:gls}, there exist some
$m \in \mathbb{N}$, $v_{ij} \in E_2 \cap B_r(q)$, and some
$\alpha_{ij} \in \mathbb{R}$ such that
$y - q = \sum\limits_{i=1}^m \sum\limits_{j=1}^\infty
 \alpha_{ij} (v_{ij} - q) \in M_a$.
Since $f : E_1 \to E_2$ is surjective, there exists a
$u_{ij} \in E_1$ satisfying $v_{ij} = f(u_{ij})$ for any $i$
and $j$.
Moreover, by Lemma \ref{lem:new2}, we have
\begin{eqnarray*}
\begin{split}
\| u_{ij} - p \|_a^2
& = \langle u_{ij} - p, u_{ij} - p \rangle_a \\
& = \big\langle (T_{-q} \circ f \circ T_p)(u_{ij} - p),\,
                (T_{-q} \circ f \circ T_p)(u_{ij} - p)
    \big\rangle_a \\
& = \big\langle f(u_{ij}) - q, f(u_{ij}) - q \big\rangle_a \\
& = \langle v_{ij} - q, v_{ij} - q \rangle_a \\
& = \| v_{ij} - q \|_a^2 \\
& < r^2
\end{split}
\end{eqnarray*}
for any $i$ and $j$.
So we conclude that $u_{ij} \in E_1 \cap B_r(p)$ and
$v_{ij} = f(u_{ij})$ for all $i$ and $j$.

On the other hand, using Lemma \ref{lem:new2}, we have
\begin{eqnarray*}
\begin{split}
& \bigg\|
  \sum_{i=1}^m \sum_{j=1}^\infty \alpha_{ij} (u_{ij} - p)
  \bigg\|_a^2 \\
& \hspace*{3mm}
  = \left\langle
    \sum_{i=1}^m \sum_{j=1}^\infty \alpha_{ij} (u_{ij} - p),\,
    \sum_{k=1}^m \sum_{\ell=1}^\infty
    \alpha_{k\ell} (u_{k\ell} - p)
    \right\rangle_a \\
& \hspace*{3mm}
  = \sum_{i=1}^m \sum_{k=1}^m
    \sum_{j=1}^\infty \alpha_{ij}
    \sum_{\ell=1}^\infty \alpha_{k\ell}
    \big\langle u_{ij} - p,\, u_{k\ell} - p \big\rangle_a \\
& \hspace*{3mm}
  = \sum_{i=1}^m \sum_{k=1}^m
    \sum_{j=1}^\infty \alpha_{ij}
    \sum_{\ell=1}^\infty \alpha_{k\ell}
    \big\langle (T_{-q} \circ f \circ T_p)(u_{ij} - p),\,
                (T_{-q} \circ f \circ T_p)(u_{k\ell} - p)
    \big\rangle_a \\
& \hspace*{3mm}
  = \left\langle
    \sum_{i=1}^m \sum_{j=1}^\infty \alpha_{ij}
    (T_{-q} \circ f \circ T_p)(u_{ij} - p),\,
    \sum_{k=1}^m \sum_{\ell=1}^\infty \alpha_{k\ell}
    (T_{-q} \circ f \circ T_p)(u_{k\ell} - p)
    \right\rangle_a \\
& \hspace*{3mm}
  = \bigg\|
    \sum_{i=1}^m \sum_{j=1}^\infty
    \alpha_{ij} (T_{-q} \circ f \circ T_p)(u_{ij} - p)
    \bigg\|_a^2 \\
& \hspace*{3mm}
  = \bigg\|
    \sum_{i=1}^m \sum_{j=1}^\infty
    \alpha_{ij} \big( f(u_{ij}) - q \big)
    \bigg\|_a^2 \\
& \hspace*{3mm}
  = \bigg\|
    \sum_{i=1}^m \sum_{j=1}^\infty \alpha_{ij} (v_{ij} - q)
    \bigg\|_a^2 \\
& \hspace*{3mm}
  < \infty,
\end{split}
\end{eqnarray*}
since
$\sum\limits_{i=1}^m \sum\limits_{j=1}^\infty
 \alpha_{ij} (v_{ij} - q) = y - q \in M_a$.
Thus, Remark \ref{rem:new2.1} implies that
$\sum\limits_{i=1}^m \sum\limits_{j=1}^\infty
 \alpha_{ij} (u_{ij} - p) \in M_a$.

Hence, it follows from Definition \ref{def:ext} that
\begin{eqnarray*}
\begin{split}
y & =   q + \sum_{i=1}^m \sum_{j=1}^\infty
        \alpha_{ij} \big( f(u_{ij}) - q \big) \\
  & =   q + \sum_{i=1}^m \sum_{j=1}^\infty
        \alpha_{ij} (T_{-q} \circ f \circ T_p)(u_{ij} - p) \\
  & =   q + (T_{-q} \circ F \circ T_p)
        \bigg(
        \sum_{i=1}^m \sum_{j=1}^\infty \alpha_{ij} (u_{ij} - p)
        \bigg) \\
  & =   F \bigg( p + \sum_{i=1}^m \sum_{j=1}^\infty
                 \alpha_{ij} (u_{ij} - p)
          \bigg) \\
  & \in F({\rm GS}(E_1,p)).
\end{split}
\end{eqnarray*}
Thus, we conclude that
${\rm GS}(E_2,q) \subset F({\rm GS}(E_1,p))$.

$(d)$
Similarly, referring to the changes presented in the tables
below and following the previous parts $(b)$ and $(c)$ in this
proof, we can prove that
${\rm GS}^2(E_2,q) = F_2({\rm GS}^2(E_1,p))$.

\begin{center}
\begin{tabular}{ccccccc}
\hline\hline
The case $n = 1$:   & $E_1$               & $E_2$             &
${\rm GS}(E_1,p)$   & ${\rm GS}(E_2,q)$   & $f$ & $F$   \\
\hline
The case $n = 2$:   & ${\rm GS}(E_1,p)$   & ${\rm GS}(E_2,q)$ &
${\rm GS}^2(E_1,p)$ & ${\rm GS}^2(E_2,q)$ & $F$ & $F_2$ \\
\hline\hline
\end{tabular}
\end{center}

\begin{center}
\begin{tabular}{cccc}
\hline\hline
The case $n = 1$:   & Definition \ref{def:gls}
& Definition \ref{def:ext}        & Lemma \ref{lem:new2}   \\
\hline
The case $n = 2$:   & Definition \ref{def:6.1} $(i)$
& Definition \ref{def:6.1} $(ii)$ & $(\ref{eq:2018-9-28})$ \\
\hline\hline
\end{tabular}
\end{center}

$(e)$
Finally, according to Proposition \ref{rem:new3.1}, Theorem
\ref{thm:setyol}, and $(d)$, we further have
\begin{eqnarray*}
{\rm GS}^n(E_2,q) = {\rm GS}^2(E_2,q) = F_2({\rm GS}^2(E_1,p))
= F_n({\rm GS}^n(E_1,p))
\end{eqnarray*}
for any integer $n \geq 3$.
\hfill$\Box$
\vspace*{5mm}


\section{Extension of isometries to the entire space}
\label{sec:8}

Let $I^\omega = \prod\limits_{i=1}^\infty I$ be the
\emph{Hilbert cube}, where $I = [0,1]$ is the unit closed
interval.
From now on, we assume that $E_1$ and $E_2$ are nonempty
subsets of $I^\omega$.
They are bounded, of course.

In Theorem \ref{thm:3.17}, we will prove that the domain of a
local $d_a$-isometry $f : E_1 \to E_2$ can be extended to any
real Hilbert space including the domain $E_1$.


\begin{defn} \label{def:3.5}
Let $E_1$ be a nonempty subset of $I^\omega$ that is
$d_a$-isometric to a subset $E_2$ of $I^\omega$ via a
surjective $d_a$-isometry $f : E_1 \to E_2$.
Let $p$ be an element of $E_1$ and $q$ an element of $E_2$
with $q = f(p)$.
Assume that $\{ \frac{1}{a_i} e_i \}_{i \in \Lambda_\alpha}$
is a complete orthonormal sequence in the Hilbert space
${\rm GS}^2(E_1,p) - p$, where $\Lambda_\alpha$ is a nonempty
proper subset of $\mathbb{N}$.
Moreover, assume that $\{ \beta_i \}_{i \in \mathbb{N}}$ is a
complete orthonormal sequence in the Hilbert space $M_a$ such
that $\beta_i = \frac{1}{a_i} (T_{-q} \circ F_2 \circ T_p)(e_i)$
for each $i \in \Lambda_\alpha$, where
$F_2 : {\rm GS}^2(E_1,p) \to M_a$ is defined in Definition
\ref{def:6.1}.
Let $p_i$ be the $i$th component of $p$, \textit{i.e.},
$p = \sum\limits_{i=1}^\infty p_i e_i$.
For any set $\Lambda$ satisfying
$\Lambda_\alpha \subset \Lambda \subset \mathbb{N}$, we define
a basic cylinder or a basic interval $\tilde{J}$ by
\begin{eqnarray*}
\tilde{J} = \prod_{i=1}^\infty\, \tilde{J}_i,
~~~\mbox{where}~~~
\tilde{J}_i = \left\{\begin{array}{ll}
                        [0, 1]
                        & (\mbox{for}~i \in \Lambda),
                        \vspace*{1mm}\\
                        \{ p_i \}
                        & (\mbox{for}~i \not\in \Lambda).
                     \end{array}
              \right.
\end{eqnarray*}
Moreover, referring to Theorem \ref{thm:new4}, we define the
function $G_2 : {\rm GS}^2(\tilde{J},p) \to M_a$ by
\begin{eqnarray} \label{eq:2021-10-21-1}
(T_{-q} \circ G_2 \circ T_p)(u - p)
= \sum_{i \in \Lambda(\tilde{J})}
  \left\langle u - p,\, \frac{1}{a_i} e_i \right\rangle_a
  \beta_i
\end{eqnarray}
for all $u \in {\rm GS}^2(\tilde{J},p)$.
\end{defn}

The following theorem states that the domain of a local
$d_a$-isometry can be extended to any real Hilbert space
including the domain of the local $d_a$-isometry.


\begin{thm} \label{thm:3.17}
Let $E_1$ be a bounded subset of $I^\omega$ that contains at
least two elements.
Suppose $E_1$ is $d_a$-isometric to a subset $E_2$ of
$I^\omega$ via a surjective $d_a$-isometry $f : E_1 \to E_2$.
Let $p$ and $q$ be elements of $E_1$ and $E_2$ satisfying
$q = f(p)$.
Assume that $\{ \frac{1}{a_i} e_i \}_{i \in \Lambda_\alpha}$
is a complete orthonormal sequence in the Hilbert space
${\rm GS}^2(E_1,p) - p$, where $\Lambda_\alpha$ is a nonempty
proper subset of $\mathbb{N}$.
Moreover, assume that $\{ \beta_i \}_{i \in \mathbb{N}}$ is a
complete orthonormal sequence in the Hilbert space $M_a$ such
that $\beta_i = \frac{1}{a_i} (T_{-q} \circ F_2 \circ T_p)(e_i)$
for each $i \in \Lambda_\alpha$.
Let $\Lambda$ be a set satisfying
$\Lambda_\alpha \subset \Lambda \subset \mathbb{N}$ and let
$\tilde{J}$ be defined as in Definition \ref{def:3.5}.
Then the function $G_2 : {\rm GS}^2(\tilde{J},p) \to M_a$ is a
$d_a$-isometry and the function
$T_{-q} \circ G_2 \circ T_p : {\rm GS}^2(\tilde{J},p) - p \to M_a$
is linear.
In particular, $G_2$ is an extension of $F_2$.
\end{thm}

\noindent
\emph{Proof}.
$(a)$
First, we assert that the function
$T_{-q} \circ G_2 \circ T_p : {\rm GS}^2(\tilde{J},p) - p \to M_a$
preserves the inner product.
Assume that $u$ and $v$ are arbitrary elements of
${\rm GS}^2(\tilde{J},p)$.
Since $\Lambda = \Lambda(\tilde{J})$, it follows from
(\ref{eq:2021-10-21-2}), (\ref{eq:2021-10-21-1}), and the
orthonormality of $\{ \frac{1}{a_i} e_i \}_{i \in \mathbb{N}}$
and $\{ \beta_i \}_{i \in \mathbb{N}}$ that
\begin{eqnarray*}
\begin{split}
& \big\langle (T_{-q} \circ G_2 \circ T_p)(u - p),\,
              (T_{-q} \circ G_2 \circ T_p)(v - p)
  \big\rangle_a \\
& \hspace*{3mm}
  = \left\langle \sum_{i \in \Lambda}
                 \left\langle u - p,\, \frac{1}{a_i} e_i
                 \right\rangle_a \beta_i,\;
                 \sum_{j \in \Lambda}
                 \left\langle v - p,\, \frac{1}{a_j} e_j
                 \right\rangle_a \beta_j
    \right\rangle_a \\
& \hspace*{3mm}
  = \sum_{i \in \Lambda}
    \left\langle u - p,\, \frac{1}{a_i} e_i \right\rangle_a
    \sum_{j \in \Lambda}
    \left\langle v - p,\, \frac{1}{a_j} e_j \right\rangle_a
    \langle \beta_i,\, \beta_j \rangle_a \\
& \hspace*{3mm}
  = \sum_{i \in \Lambda}
    \left\langle u - p,\, \frac{1}{a_i} e_i \right\rangle_a
    \sum_{j \in \Lambda}
    \left\langle v - p,\, \frac{1}{a_j} e_j \right\rangle_a
    \left\langle \frac{1}{a_i} e_i,\, \frac{1}{a_j} e_j
    \right\rangle_a \\
& \hspace*{3mm}
  = \left\langle \sum_{i \in \Lambda}
                 \left\langle u - p,\, \frac{1}{a_i} e_i
                 \right\rangle_a
                 \frac{1}{a_i} e_i,\;
                 \sum_{j \in \Lambda}
                 \left\langle v - p,\, \frac{1}{a_j} e_j
                 \right\rangle_a
                 \frac{1}{a_j} e_j
    \right\rangle_a \\
& \hspace*{3mm}
  = \langle u - p, v - p \rangle_a
\end{split}
\end{eqnarray*}
for all $u, v \in {\rm GS}^2(\tilde{J},p)$, \textit{i.e.},
$T_{-q} \circ G_2 \circ T_p$ preserves the inner product.

$(b)$
We assert that $G_2$ is a $d_a$-isometry.
Let $u$ and $v$ be arbitrary elements of
${\rm GS}^2(\tilde{J},p)$.
Since $T_{-q} \circ G_2 \circ T_p$ preserves the inner product
by $(a)$, we have
\begin{eqnarray*}
\begin{split}
& d_a \big( G_2(u), G_2(v) \big)^2 \\
& \hspace*{3mm}
  = \big\| ( T_{-q} \circ G_2 \circ T_p )(u-p) -
           ( T_{-q} \circ G_2 \circ T_p )(v-p)
    \big\|_a^2 \\
& \hspace*{3mm}
  = \big\langle ( T_{-q} \circ G_2 \circ T_p )(u-p) -
                ( T_{-q} \circ G_2 \circ T_p )(v-p), \\
& \hspace*{3mm}
    ~~~~~~      ( T_{-q} \circ G_2 \circ T_p )(u-p) -
                ( T_{-q} \circ G_2 \circ T_p )(v-p)
    \big\rangle_a \\
& \hspace*{3mm}
  = \langle u-p, u-p \rangle_a - \langle u-p, v-p \rangle_a
    - \langle v-p, u-p \rangle_a + \langle v-p, v-p \rangle_a \\
& \hspace*{3mm}
  = \big\langle (u-p) - (v-p), (u-p) - (v-p) \big\rangle_a \\
& \hspace*{3mm}
  = \| (u-p) - (v-p) \|_a^2 \\
& \hspace*{3mm}
  = \| u - v \|_a^2 \\
& \hspace*{3mm}
  = d_a(u, v)^2
\end{split}
\end{eqnarray*}
for all $u, v \in {\rm GS}^2(\tilde{J},p)$, \textit{i.e.},
$G_2 : {\rm GS}^2(\tilde{J},p) \to M_a$ is a $d_a$-isometry.

$(c)$
Now we assert that the function
$T_{-q} \circ G_2 \circ T_p : {\rm GS}^2(\tilde{J},p) - p \to M_a$
is linear.
Assume that $u$ and $v$ are arbitrary elements of
${\rm GS}^2(\tilde{J},p)$ and $\alpha$ and $\beta$ are real
numbers.
Since ${\rm GS}^2(\tilde{J},p) - p$ is a real vector space, it
holds that 
$\alpha (u - p) + \beta (v - p) \in {\rm GS}^2(\tilde{J},p) - p$.
Thus, $\alpha (u - p) + \beta (v - p) = w - p$ for some
$w \in {\rm GS}^2(\tilde{J},p)$.
Hence, referring to the changes presented in the table below
and following $(d)$ of the proof of Theorem \ref{thm:new3}, we
can easily prove that $T_{-q} \circ G_2 \circ T_p$ is linear.

\begin{center}
\begin{tabular}{cccc}
\hline\hline
Theorem \ref{thm:new3}: & ${\rm GS}(E_1,p)$         & $F$
& (\ref{eq:2018-9-28}) \\
\hline
Here:                   & ${\rm GS}^2(\tilde{J},p)$ & $G_2$
& $(a)$  \\
\hline\hline
\end{tabular}
\end{center}

$(d)$
Finally, we assert that $G_2$ is an extension of $F_2$.
Let $\hat{J}$ be either a basic cylinder or a basic interval
defined by
\begin{eqnarray*}
\hat{J} = \prod_{i=1}^\infty \hat{J}_i, ~~~\mbox{where}~~~
\hat{J}_i
= \left\{\begin{array}{ll}
            [0, 1]    & (\mbox{for}~i \in \Lambda_\alpha),
            \vspace*{1mm}\\
            \{ p_i \} & (\mbox{for}~i \not\in \Lambda_\alpha).
         \end{array}
  \right.
\end{eqnarray*}
We see that $p = (p_1, p_2, \ldots) \in \hat{J} \cap E_1$ and
$\Lambda(\hat{J}) = \Lambda_\alpha = \Lambda({\rm GS}^2(E_1,p))$.

According to Lemma \ref{lem:setahob}, if
$i \in \Lambda({\rm GS}^2(E_1,p))$, then
$\alpha_i e_i \in {\rm GS}^2(E_1,p) - p$ for all
$\alpha_i \in \mathbb{R}$.
Since ${\rm GS}^2(E_1,p) - p$ is a real vector space, if we set
$\Lambda_n = \{ i \in \Lambda({\rm GS}^2(E_1,p)) : i < n \}$,
then we have
\begin{eqnarray*}
\sum_{i \in \Lambda_n} \alpha_i e_i \in {\rm GS}^2(E_1,p) - p
\end{eqnarray*}
for all $n \in \mathbb{N}$ and all $\alpha_i \in \mathbb{R}$.
For now, with all $\alpha_i$ fixed, we define
$x_n = p + \sum\limits_{i \in \Lambda_n} \alpha_i e_i$ for any
$n \in \mathbb{N}$.
Then $\{ x_n \}$ is a sequence in ${\rm GS}^2(E_1,p)$.
When $\{ x_n \}$ converges in $M_a$, it holds that
\begin{eqnarray*}
p + \sum\limits_{i \in \Lambda({\rm GS}^2(E_1,p))} \alpha_i e_i
=   \lim_{n \to \infty} x_n
\in {\rm GS}^2(E_1,p),
\end{eqnarray*}
because ${\rm GS}^2(E_1,p)$ is closed by Lemma \ref{lem:3.8}
$(iii)$.
That is,
\begin{eqnarray*}
\Bigg\{ p + \sum_{i \in \Lambda({\rm GS}^2(E_1,p))}
        \alpha_i e_i \in M_a :
        \alpha_i \in \mathbb{R} ~\mbox{for all}~
        i \in \Lambda \big( {\rm GS}^2(E_1,p) \big)
\Bigg\} \subset {\rm GS}^2(E_1,p).
\end{eqnarray*}

Hence, by Theorem \ref{thm:set9} with
$\beta = \{ \frac{1}{a_i} e_i \}_{i \in \mathbb{N}}$ and
$J_\beta = \hat{J}$, we get
\begin{eqnarray*}
\begin{split}
{\rm GS}(\hat{J},p) - p
& =       \Bigg\{ \sum_{i \in \Lambda(\hat{J})}
                  \alpha_i e_i \in M_a :
                  \alpha_i \in \mathbb{R} ~\mbox{for all}~
                  i \in \Lambda(\hat{J})
          \Bigg\} \\
& =       \Bigg\{ \sum_{i \in \Lambda({\rm GS}^2(E_1,p))}
                  \alpha_i e_i \in M_a :
                  \alpha_i \in \mathbb{R} ~\mbox{for all}~
                  i \in \Lambda \big( {\rm GS}^2(E_1,p) \big)
          \Bigg\} \\
& \subset {\rm GS}^2(E_1,p) - p.
\end{split}
\end{eqnarray*}

So we have
\begin{eqnarray*}
\hat{J} \cap B_r(p) \subset {\rm GS}(\hat{J},p) \cap B_r(p)
\subset {\rm GS}^2(E_1,p) \cap B_r(p)
\end{eqnarray*}
for some real number $r > 0$ and hence, we further have
\begin{eqnarray*}
{\rm GS}(\hat{J},p) \subset {\rm GS}^2(\hat{J},p)
\subset {\rm GS}^3(E_1,p) = {\rm GS}^2(E_1,p).
\end{eqnarray*}
Moreover, by Remark \ref{rem:7.2}, we know that
${\rm GS}^2(\hat{J},p) = {\rm GS}(\hat{J},p)$.
Hence, we have
\begin{eqnarray*}
{\rm GS}(\hat{J},p) = {\rm GS}^2(\hat{J},p)
\subset {\rm GS}^2(E_1,p).
\end{eqnarray*}

On the other hand, since
$\{ \frac{1}{a_i} e_i \}_{i \in \Lambda_\alpha}$ is a complete
orthonormal sequence in ${\rm GS}^2(E_1,p) - p$, it follows
from Theorem \ref{thm:set9} with
$\beta = \{ \frac{1}{a_i} e_i \}_{i \in \mathbb{N}}$ that
\begin{eqnarray*}
x =   \sum_{i \in \Lambda_\alpha}
      \left\langle x, \frac{1}{a_i} e_i \right\rangle_a
      \frac{1}{a_i} e_i
  \in {\rm GS}(\hat{J},p) - p
  =   {\rm GS}^2(\hat{J},p) - p
\end{eqnarray*}
for all $x \in {\rm GS}^2(E_1,p) - p$, which implies that
${\rm GS}^2(E_1,p) = {\rm GS}^2(\hat{J},p) = {\rm GS}(\hat{J},p)$.

Let $u$ be an arbitrary element of ${\rm GS}^2(E_1,p)$.
Then by (\ref{eq:2021-10-21-2}) with $\hat{J}$ instead of $J$,
we have
\begin{eqnarray} \label{eq:2021-11-8}
u - p = \sum_{i \in \Lambda(\hat{J})}
        \left\langle u - p,\, \frac{1}{a_i} e_i \right\rangle_a
        \frac{1}{a_i} e_i
\end{eqnarray}
and since $T_{-q} \circ F_2 \circ T_p$ is linear and continuous,
we use (\ref{eq:2021-10-21-1}), (\ref{eq:2021-11-8}), and the
facts
${\rm GS}^2(E_1,p) = {\rm GS}^2(\hat{J},p) = {\rm GS}(\hat{J},p)$
and
$\Lambda(\hat{J}) = \Lambda_\alpha = \Lambda({\rm GS}^2(E_1,p))$
to have
\begin{eqnarray*}
\begin{split}
(T_{-q} \circ G_2 \circ T_p)(u - p)
& = \sum_{i \in \Lambda(\tilde{J})}
    \left\langle u - p,\, \frac{1}{a_i} e_i \right\rangle_a
    \beta_i \\
& = \sum_{i \in \Lambda(\hat{J})}
    \left\langle u - p,\, \frac{1}{a_i} e_i \right\rangle_a
    \beta_i \\
& = \sum_{i \in \Lambda(\hat{J})}
    \left\langle u - p,\, \frac{1}{a_i} e_i \right\rangle_a
    \frac{1}{a_i} (T_{-q} \circ F_2 \circ T_p)(e_i) \\
& = \lim_{n \to \infty}
    \sum_{i \in \Lambda_n(\hat{J})}
    \left\langle u - p,\, \frac{1}{a_i} e_i \right\rangle_a
    \frac{1}{a_i} (T_{-q} \circ F_2 \circ T_p)(e_i) \\
& = \lim_{n \to \infty}
    (T_{-q} \circ F_2 \circ T_p)\!
    \left(
    \sum_{i \in \Lambda_n(\hat{J})}
    \left\langle u - p,\, \frac{1}{a_i} e_i \right\rangle_a
    \frac{1}{a_i} e_i
    \right) \\
& = (T_{-q} \circ F_2 \circ T_p)\!
    \left(
    \sum_{i \in \Lambda(\hat{J})}
    \left\langle u - p,\, \frac{1}{a_i} e_i \right\rangle_a
    \frac{1}{a_i} e_i
    \right) \\
& = (T_{-q} \circ F_2 \circ T_p)(u - p),
\end{split}
\end{eqnarray*}
where we set
$\Lambda_n(\hat{J}) = \{ i \in \Lambda(\hat{J}) : i < n \}$
for every $n \in \mathbb{N}$.

Therefore, it follows that
$(T_{-q} \circ G_2 \circ T_p)(u - p)
 = (T_{-q} \circ F_2 \circ T_p)(u - p)$ for all
$u \in {\rm GS}^2(E_1,p)$, \textit{i.e.}, $G_2(u) = F_2(u)$
for all $u \in {\rm GS}^2(E_1,p)$.
In other words, $G_2$ is an extension of $F_2$.
\hfill$\Box$


\section{Cylinders}

For each positive integer $n$, let $\mathcal{B}_n$ be the set
of all basic cylinders $J = \prod\limits_{i=1}^\infty J_i$
defined by Definition \ref{def:basic} for some disjoint finite
subsets $\Lambda_1$, $\Lambda_2$, $\Lambda_3$, $\Lambda_4$ of
$\{ 1, 2, \ldots, n \}$ and $0 < p_{1i} < p_{2i} < 1$ for
$i \in \Lambda_1 \cup \Lambda_2 \cup \Lambda_3$ and
$0 \leq p_{1i} \leq 1$ for $i \in \Lambda_4$.
This definition of basic cylinders is a slight modification of
Definition \ref{def:basic}, but the two definitions are
essentially the same.
We note that
$\Lambda_1 \cup \Lambda_2 \cup \Lambda_3 \cup \Lambda_4
 \subset \{ 1, 2, \ldots, n \}$ and at most $n$ edges of each
basic cylinder of $\mathcal{B}_n$ have a Euclidean length of
less than $1$, and all remaining edges have a Euclidean length
of $1$.

For $i \in \{ 1, 2, \ldots, n \}$, each $J_i$ is a closed
subinterval of $[0,1]$ with respect to the relative topology
for $[0,1]$. 
Thus, the infinite dimensional interval
\begin{eqnarray*}
\prod_{1 \leq j < i} [0,1] \times J_i \times
\prod_{j>i}\; [0,1]
\end{eqnarray*}
is a closed subset of the Hilbert cube $I^\omega$ for each
$i \in \{ 1, 2, \ldots, n \}$.
In particular, when $i = 1$, we read the last expression as
\begin{eqnarray*}
\prod_{1 \leq j < 1} [0,1] \times J_1 \times
\prod_{j>1}\; [0,1]
= J_1 \times \prod_{j>1}\; [0,1].
\end{eqnarray*}

Since every basic cylinder $J \in \mathcal{B}_n$ is expressed
as
\begin{eqnarray*}
J = \prod\limits_{i=1}^\infty J_i
  = \bigcap_{i=1}^n
    \bigg( \prod_{1 \leq j < i} [0,1] \times J_i \times
           \prod_{j>i}\; [0,1]
    \bigg),
\end{eqnarray*}
we see that each $J \in \mathcal{B}_n$ is a closed subset of
$I^\omega$ as the intersection of closed sets.
Since $I^\omega$ is a closed subset of $M_a$ by Remark
\ref{rem:2.1} $(iv)$, we use \cite[Theorem 80.4]{kasriel} to
conclude that $J$ is a closed subset of $M_a$.

Now, let us define
\begin{eqnarray*}
\mathcal{B} = \{ \emptyset \} \cup
              \bigcup_{n=1}^\infty \mathcal{B}_n
~~~\mbox{and}~~~
\mathcal{B}_\delta
= \{ J \in \mathcal{B} \,:\, d_a(J) < \delta \}
\end{eqnarray*}
for every $\delta > 0$, where $d_a(J)$ is the \emph{diameter}
of $J$ defined as $\sup \{ d_a(x,y) : x, y \in J \}$.
We note that every basic cylinder in $\mathcal{B}$ is a closed
subset of $M_a$.

We denote by $\mathcal{C}$ the set of every subset $K$ of
$M_a$, for which there exist a basic cylinder
$J \in \mathcal{B}$ and a surjective $d_a$-isometry
$f : J \to K$, and we define
$\mathcal{C}_\delta
 = \{ K \in \mathcal{C} \,:\, d_a(K) < \delta \}$ for any
$\delta > 0$.
We notice that $\mathcal{B}_\delta \subset \mathcal{C}_\delta$
for every $\delta > 0$.
We note that the family $\mathcal{B}$ includes not only
non-degenerate basic cylinders but also degenerate ones.

Assume that a basic cylinder $J$ and a cylinder $K$ are given
such that $J$ is $d_a$-isometric to $K$ through a surjective
$d_a$-isometry $f : J \to K$.
Since $J$ is a compact subset of $M_a$ as a closed subset of a
compact set $I^\omega$, $K$ is also a compact subset of $M_a$
as the continuous image of a compact set $J$
(see \cite[Theorem 91.7]{kasriel}).
Moreover, $K$ is a closed subset of $M_a$ as a compact subset
of the Hausdorff space $M_a$
(see \cite[Theorem 91.2]{kasriel}).
Hence, we come to an important consequence.


\begin{rem} \label{rem:3.1}
Every cylinder $K \in \mathcal{C}$ is closed in $M_a$.
\end{rem}

Let $J$ be a basic cylinder that is $d_a$-isometric to a
cylinder $K$ via a surjective $d_a$-isometry $f : J \to K$.
Assume that $p = (p_1, p_2, \ldots)$ is the lower left corner
of $J$, $q$ is an element of $K$ with $q = f(p)$, ${\rm GS}(J,p)$
is the first-order generalized linear span of $J$ with respect
to $p$, and that $F : {\rm GS}(J,p) \to M_a$ is the extension
of $f$ given in Definition \ref{def:ext}.
If $x = (x_1, x_2, \ldots) \in {\rm GS}(J,p)$, then it follows
from Theorem \ref{thm:set9} that
$x-p =   \sum\limits_{j=1}^\infty (x_j - p_j) e_j
     =   \sum\limits_{j \in \Lambda(J)} (x_j - p_j) e_j
     \in M_a$.
By Theorem \ref{thm:new4}, we get
\begin{eqnarray} \label{eq:2018-9-27}
\begin{split}
(T_{-q} \circ F \circ T_p)(x - p)
& = \sum_{i \in \Lambda(J)}
    \left\langle x - p,\, \frac{1}{a_i} e_i \right\rangle_a
    \frac{1}{a_i} (T_{-q} \circ F \circ T_p)(e_i) \\
& = \sum_{i \in \Lambda(J)} a_i (x_i - p_i)
    \frac{1}{a_i} (T_{-q} \circ F \circ T_p)(e_i).
\end{split}
\end{eqnarray}

For any $y = (y_1, y_2, \ldots) \in {\rm GS}(J,p)$, it follows
from Theorem \ref{thm:set9} that
\begin{eqnarray*}
y =   p + (y - p)
  =   p + \sum\limits_{i=1}^\infty (y_i - p_i) e_i
  =   p + \sum\limits_{i \in \Lambda(J)}
      a_i (y_i - p_i) \frac{1}{a_i} e_i
  \in {\rm GS}(J,p)
\end{eqnarray*}
and it further follows from (\ref{eq:2018-9-27}) that
\begin{eqnarray*}
\begin{split}
F(y) & = q + ( T_{-q} \circ F \circ T_p )(y-p) \\
     & = q + \sum_{i \in \Lambda(J)} a_i (y_i - p_i)
         \frac{1}{a_i} ( T_{-q} \circ F \circ T_p )(e_i) \\
     & \in {\rm GS}(K,q),
\end{split}
\end{eqnarray*}
since ${\rm GS}(K,q) = F({\rm GS}(J,p))$ by Theorem
\ref{thm:2018-11-7}.
Moreover, since the sequences
$\{ \frac{1}{a_i} e_i \}_{i \in \Lambda(J)}$ and
$\{ \frac{1}{a_i}
    ( T_{-q} \circ F \circ T_p )(e_i) \}_{i \in \Lambda(J)}$
are both orthonormal, the following definition may be useful.


\begin{defn} \label{def:2.2}
Every interval in $\mathcal{B}$ will be called a
\emph{basic cylinder} and each element of\/ $\mathcal{C}$ a
\emph{cylinder}.
Assume that a basic cylinder $J$ and a cylinder $K$ are given
such that $K = f(J)$ for some surjective $d_a$-isometry
$f : J \to K$.
If $\Lambda(J) = \mathbb{N}$, then $J$ and $K$ will be called
\emph{non-degenerate}.
Otherwise, they will be called \emph{degenerate}.
\end{defn}


\begin{rem} \label{rem:4.a}
The term `degenerate' or `non-degenerate' defined in relation
to basic cylinders and cylinders is similar to, but not
identical to, the term `degenerate' or `non-degenerate' for
the general sets E defined in Definition \ref{def:Lambda}.
We notice that the terms `degenerate' and `non-degenerate' are
used for convenience only and are not exact mathematical terms.
\end{rem}


\section{Elementary volumes}

Assume that both $J_1$ and $J_2$ are two distinct basic
cylinders which are $d_a$-isometric to the same cylinder $K$
via the surjective $d_a$-isometries $f_1 : J_1 \to K$ and
$f_2 : J_2 \to K$, respectively.
Moreover, assume that $u$ is the lower left corner and $x$ is
the vertex of $J_1$ diagonally opposite to $u$, \textit{i.e.},
$x$ is the upper right corner of $J_1$.
Analogously, let $v$ be the lower left corner and $y$ the
vertex of $J_2$ diagonally opposite to $v$ and
$f_1(u) = f_2(v) =: w \in K$.
Furthermore, assume that $F_{(1)} : {\rm GS}(J_1,u) \to M_a$
and $F_{(2)} : {\rm GS}(J_2,v) \to M_a$ are $d_a$-isometries
given in Definition \ref{def:ext} and that they are extensions
of $f_1$ and $f_2$, respectively.

Then, by (\ref{eq:2018-9-27}), we have
\begin{eqnarray} \label{eq:2018-9-29}
\begin{split}
(T_{-w} \circ F_{(1)} \circ T_u)(x - u)
& = (T_{-w} \circ F_{(1)} \circ T_u)
    \Bigg( \sum_{i \in \Lambda(J_1)} (x_i - u_i) e_i \Bigg) \\
& = \sum_{i \in \Lambda(J_1)}
    (x_i - u_i) (T_{-w} \circ F_{(1)} \circ T_u)(e_i), \\
(T_{-w} \circ F_{(2)} \circ T_v)(y - v)
& = (T_{-w} \circ F_{(2)} \circ T_v)
    \Bigg( \sum_{i \in \Lambda(J_2)} (y_i - v_i) e_i \Bigg) \\
& = \sum_{i \in \Lambda(J_2)}
    (y_i - v_i) (T_{-w} \circ F_{(2)} \circ T_v)(e_i).
\end{split}
\end{eqnarray}
Further, the right hand side of the first equality expresses
the vector $f_1(x) - w$, since
\begin{eqnarray*}
f_1(x) - w = F_{(1)}(x) - w
           = ( T_{-w} \circ F_{(1)} \circ T_u )(x - u).
\end{eqnarray*}
Similarly, the right hand side of the second equality in
(\ref{eq:2018-9-29}) expresses the vector $f_2(y) - w$.

According to (\ref{eq:2021-10-21-2}), Theorem \ref{thm:new4}
and (\ref{eq:2018-9-29}), the coordinates
$\big\langle x-u, \frac{1}{a_i} e_i \big\rangle_a$ remain
unchanged under the action of the $d_a$-isometry
$T_{-w} \circ F_{(1)} \circ T_u$.
Moreover, the points $w$ and $f_1(x)$ are the diagonally
opposite vertices of the cylinder $K$.
The same is true for $w$ and $f_2(y)$.
Thus, we can conclude that $f_1(x) - w = f_2(y) - w$.

Since
$(T_{-w} \circ F_{(1)} \circ T_u)((x_i - u_i) e_i)
 = (x_i - u_i) (T_{-w} \circ F_{(1)} \circ T_u)(e_i)$ for each
$i \in \Lambda(J_1)$, $F_{(1)}$ maps each edge of basic
cylinder $J_1$ onto the edge of $K$.
Conversely, every edge of the cylinder $K$ is an image of the
edge of $J_1$ under the $d_a$-isometry $F_{(1)}$.
The same case is also for $F_{(2)}$ and $J_2$.
Therefore, we conclude that there exists a permutation
$\sigma : \Lambda(J_1) \to \Lambda(J_2)$ that satisfies
\begin{eqnarray*}
(x_i - u_i) ( T_{-w} \circ F_{(1)} \circ T_u )(e_i)
= (y_{\sigma(i)} - v_{\sigma(i)})
  ( T_{-w} \circ F_{(2)} \circ T_v )(e_{\sigma(i)})
\end{eqnarray*}
for any $i \in \Lambda(J_1)$.

According to (\ref{eq:2018-9-28}), both
$\{ \frac{1}{a_i} ( T_{-w} \circ F_{(1)} \circ T_u ) (e_i) \}$
and
$\{ \frac{1}{a_i} ( T_{-w} \circ F_{(2)} \circ T_v ) (e_i) \}$
are orthonormal sequences.
Hence, we get
\begin{eqnarray} \label{eq:2018-10-9}
a_i | x_i - u_i |
= a_{\sigma(i)} | y_{\sigma(i)} - v_{\sigma(i)} |
\end{eqnarray}
for all $i \in \Lambda(J_1)$.
Since $J_1$ and $J_2$ are both basic cylinders, due to the
structural property of basic cylinders, we can see that there
is an $\ell_0 \in \mathbb{N}$ that satisfies
$| x_i - u_i | = | y_i - v_i | = 1$ for all $i > \ell_0$.
Thus, it follows from (\ref{eq:2018-10-9}) that there exists
an $m_0 \in \mathbb{N}$ ($m_0 \geq \ell_0$) that satisfies
$a_i = a_{\sigma(i)}$ for each $i > m_0$.

Consequently, when $J_1$ is non-degenerate, we use
(\ref{eq:2018-10-9}) to show that the basic cylinders $J_1$
and $J_2$ have the same \emph{elementary volume}:
\begin{eqnarray*}
\begin{split}
{\rm vol}(J_1)
& = \prod_{i=1}^\infty | x_i - u_i |
  = \prod_{i=1}^\infty
    \frac{a_{\sigma(i)}}{a_i}
    \big| y_{\sigma(i)} - v_{\sigma(i)} \big| \\
& = \prod_{i=1}^{m_0} \frac{a_{\sigma(i)}}{a_i}
    \big| y_{\sigma(i)} - v_{\sigma(i)} \big| \times
    \prod_{i = m_0 + 1}^\infty
    \big| y_{\sigma(i)} - v_{\sigma(i)} \big| \\
& = \prod_{i=1}^{m_0} \frac{a_{\sigma(i)}}{a_i} \times
    \prod_{i=1}^\infty
    \big| y_{\sigma(i)} - v_{\sigma(i)} \big| \\
& = \prod_{i=1}^{m_0} \frac{a_{\sigma(i)}}{a_i} \times
    \prod_{i=1}^\infty | y_i - v_i | \\
& = \prod_{i=1}^\infty | y_i - v_i | \\
& = {\rm vol}(J_2)
\end{split}
\end{eqnarray*}
Hence, it is reasonable to define the \emph{volume} of the
cylinder $K$ as the elementary volume of one of the basic
cylinders which are $d_a$-isometric to $K$, \textit{i.e.},
\begin{eqnarray*}
{\rm vol}(K) = {\rm vol}(J_1) = {\rm vol}(J_2).
\end{eqnarray*}

When $J_1$ is degenerate, we define
${\rm vol}(K) = {\rm vol}(J_1) = {\rm vol}(J_2) = 0$.



\begin{rem} \label{rem:2.2}
Let $J_1$ be a basic cylinder and $K_1$ a cylinder.
Assume that $J_1$ and $K_1$ are $d_a$-isometric to each other
via a surjective $d_a$-isometry $f : J_1 \to K_1$.
Assume that $p$ is an element of $J_1$ and $q$ is an element
of $K_1$ with $q = f(p)$.
Comparing $(\ref{eq:2021-10-21-2})$ and Theorem \ref{thm:new4}
and considering the fact that
${\rm GS}(J_1,p) = {\rm GS}^n(J_1,p)$ for any $n \in \mathbb{N}$
$($see the proof of Theorem \ref{thm:3.17}$)$, under the action
of the $d_a$-isometry
$T_{-q} \circ F \circ T_p : {\rm GS}(J_1,p) - p \to M_a$, the
following statements are true.
\begin{itemize}
\item[$(i)$]  The orthonormal sequence
              $\{ \frac{1}{a_i} e_i \}_{i \in \Lambda(J_1)}$
              is changed to
              $\{ \frac{1}{a_i}
                  ( T_{-q} \circ F \circ T_p ) (e_i)
               \}_{i \in \Lambda(J_1)}$, which is also an
              orthonormal sequence;
\item[$(ii)$] The coordinates $($or Fourier coefficients$)$
              $\big\langle u-p, \frac{1}{a_i} e_i \big\rangle_a$,
              $i \in \Lambda(J_1)$, of each element
              $u - p \in {\rm GS}(J_1,p) - p$ remain unchanged.
\end{itemize}
According to Theorem \ref{thm:new4}, the $d_a$-isometry
$T_{-q} \circ F \circ T_p$ transforms the $i$th coordinate
$\frac{1}{a_i} e_i$ into
$\frac{1}{a_i} (T_{-q} \circ F \circ T_p)(e_i)$ for every
$i \in \Lambda(J_1)$.
Moreover, by $(\ref{eq:2021-10-21-2})$ and Theorem
\ref{thm:new4}, the coordinate expression of the image of
$u-p$ under the action of $T_{-q} \circ F \circ T_p$
$($in the coordinate system
$\{ \frac{1}{a_i} (T_{-q} \circ F \circ T_p)(e_i)
 \}_{i \in \Lambda(J_1)})$
is the same as that of $u-p$ in the coordinate system
$\{ \frac{1}{a_i} e_i \}_{i \in \Lambda(J_1)}$.
Therefore, $T_{-q} \circ F \circ T_p$ preserves each $m$-face
of each basic cylinder $J$ contained in ${\rm GS}(J_1,p) - p$,
where $m \in \{ 0, 1, 2, \ldots \}$.
More precisely,
\begin{itemize}
\item[$(iii)$] $F$ maps each $m$-face of basic cylinder $J$
               contained in ${\rm GS}(J_1,p)$ onto an $m$-face
               of cylinder $F(J)$, where
               $m \in \{ 0, 1, 2, \ldots \}$.
               In particular, $F$ maps each $1$-face of $J$
               onto a $1$-face of $F(J)$.
               Consequently, the `volume' of cylinder $F(J)$
               is defined as the elementary volume of the basic
               cylinder $J$, \textit{i.e.},
               ${\rm vol}(F(J)) = {\rm vol}(J)
                = \prod\limits_{i=1}^\infty s_i$, where $s_i$
               is the Euclidean length of the $i$th edge of
               $J$.
\item[$(iv)$]  The adjacent edges of the cylinder $F(J)$ meet
               orthogonally, and the volume ${\rm vol}(F(J))$
               of $F(J)$ is the infinite product of the
               Euclidean lengths of all edges of $F(J)$.
\end{itemize}
\end{rem}

Based on Remark \ref{rem:2.2}, we can define the elementary
volume of basic cylinders and the volume of cylinders
accurately.


\begin{defn} \label{def:2018-9-23}
$(i)$ The \emph{elementary volume} of a basic cylinder $J$ is
denoted by ${\rm vol}(J)$ and defined by
\begin{eqnarray*}
{\rm vol}(J)
= \left\{\begin{array}{ll}
            \prod\limits_{i=1}^\infty s_i &
            (\mbox{for}~ \Lambda(J) = \mathbb{N}),
            \vspace*{2mm}\\
            0                             &
            (\mbox{for}~ \Lambda(J) \neq \mathbb{N}),
         \end{array}
  \right.
\end{eqnarray*}
where $s_i$ is the Euclidean length of the $i$th edge of $J$.

$(ii)$ Considering Remark \ref{rem:2.2} $(iii)$, we define the
\emph{volume} of cylinder $K$ by
\begin{eqnarray*}
{\rm vol}(K) = {\rm vol}(J)
\end{eqnarray*}
for any $K \in \mathcal{C}$ for which there exist a basic
cylinder $J \in \mathcal{B}$ and a surjective $d_a$-isometry
$f : J \to K$.
\end{defn}


\section{Construction of $d_a$-invariant measure}
\label{sec:6}

If we set ${\rm vol}(\emptyset) = 0$, then the volume
`${\rm vol}$' defined in Definition \ref{def:2018-9-23} is a
\emph{pre-measure} (see \cite[Definition 5]{6}).
According to \cite[Theorem 15]{6}, we can define an
\emph{outer measure} $\mu$ on $M_a$ by
\begin{eqnarray} \label{eq:2018-9-21-2}
\mu(E) = \lim_{\delta \to 0} \mu_\delta (E)
\end{eqnarray}
for all subsets $E$ of $M_a$, where
\begin{eqnarray*}
\mu_\delta (E) = \inf \Bigg\{ \;\sum_{i=1}^\infty {\rm vol}(C_i)
                              \,:\, C_i \in \mathcal{C}_\delta
                              ~\,\mbox{with}\,~
                              E \subset \bigcup_{i=1}^\infty C_i
                      \Bigg\}.
\end{eqnarray*}

A family $\{ C_i \}$ of sets is called a \emph{covering} of
$E$ if $E \subset \bigcup\limits_{i=1}^\infty C_i$.
If, in addition, the diameter of each $C_i$ is less than
$\delta$, then $\{ C_i \}$ is called a
\emph{$\delta$-covering} of $E$.

In this section, we assume that each of the sets $E_1$ and
$E_2$ has uncountably many elements.
We note that if $E_1$ and $E_2$ have only countably many
elements, then $\mu(E_1) = 0 = \mu(E_2)$.

One of the most important theorems in this paper is the
following theorem stating that the measure $\mu$ is
$d_a$-invariant.
However, we note that this theorem was proved for the
non-degenerate case in the paper \cite[Theorem 3.1]{jung1}.
Now we will completely prove this theorem by providing the
proof for the degenerate case also.


\begin{thm} \label{thm:5}
If $E_1$ and $E_2$ are subsets of $I^\omega$ that are
$d_a$-isometric to each other, then $\mu(E_1) = \mu(E_2)$.
\end{thm}

\noindent
\emph{Proof}.
$(a)$
We assume that $E_1$ and $E_2$ are subsets of $I^\omega$, each
of which has uncountably many elements, and they are
$d_a$-isometric to each other via the surjective $d_a$-isometry
$f : E_1 \to E_2$.
Using the definition of $F$ given in Definition \ref{def:ext}
and assuming that $p$ is an element of $E_1$ and $q$ is an
element of $E_2$ with $q = f(p)$, Theorem \ref{thm:new3} states
that $F : {\rm GS}(E_1,p) \to M_a$ is a $d_a$-isometry which
extends the surjective $d_a$-isometry $f : E_1 \to E_2$.

Let $r$ be a positive real number satisfying
$E_1 \subset B_{r}(p)$, where $B_r(p)$ denotes the open ball
defined as $B_r(p) = \{ y \in M_a \,:\, \| y - p \|_a < r \}$.
According to Theorem \ref{thm:setahob}, the function
$F_2 : {\rm GS}^2(E_1,p) \to M_a$ (defined in Definition
\ref{def:6.1}) is a $d_a$-isometry and it is an extension of
$F$ and so $F_2$ is obviously an extension of $f$.

$(b)$
We consider the case where
$\Lambda({\rm GS}^2(E_1,p)) \neq \mathbb{N}$, \textit{i.e.},
we assume that the second-order generalized linear span
${\rm GS}^2(E_1,p)$ of $E_1$ with respect to $p$ is degenerate.
In other words, according to Theorem \ref{thm:setsib},
${\rm GS}^2(E_1,p) \neq M_a$.

The translation $T_{-p} : M_a \to M_a$ is a homeomorphism and
hence, it is a closed mapping.
By Lemma \ref{lem:3.8} $(iii)$, we know that
${\rm GS}^2(E_1,p)$ is a closed proper subset of $M_a$, so is
$T_{-p}({\rm GS}^2(E_1,p))$.
That is, ${\rm GS}^2(E_1,p) - p$ is a closed (proper) subspace
of the real Hilbert space $M_a$ and hence,
${\rm GS}^2(E_1,p) - p$ is itself a real Hilbert space.

Let $\{ \alpha_i \}_{i \in \Lambda}$ be a complete orthonormal
sequence in the Hilbert space ${\rm GS}^2(E_1,p) - p$, where
$\Lambda$ is a nonempty proper subset of $\mathbb{N}$.
Moreover, we note that $M_a$ can be orthogonally decomposed
into
\begin{eqnarray*}
M_a = \big( {\rm GS}^2(E_1,p) - p \big) \oplus
      \big( {\rm GS}^2(E_1,p) - p \big)^\perp,
\end{eqnarray*}
where $({\rm GS}^2(E_1,p) - p)^\perp$ is also a real Hilbert
space as a closed subspace of the Hilbert space $M_a$.
We assume that $\beta = \{ \beta_i \}_{i \in \mathbb{N}}$ is a
complete orthonormal sequence in the Hilbert space $M_a$ such
that $\beta_i = \alpha_i$ for each $i \in \Lambda$.
We note that
\begin{eqnarray*}
\begin{split}
\Lambda & = \Lambda_\beta({\rm GS}^2(E_1,p)) \\
        & = \Big\{ i \in \mathbb{N} : \mbox{there are}~
                   z \in {\rm GS}^2(E_1,p) ~\mbox{and}~
                   \alpha \neq 0 ~\mbox{satisfying}~
                   z + \alpha \beta_i \in {\rm GS}^2(E_1,p)
            \Big\}.
\end{split}
\end{eqnarray*}

We now select a $\beta$-basic cylinder or a $\beta$-basic
interval $K$ 
that satisfies $p \in E_1 \subset K \subset I^\omega$ and
$\Lambda_\beta(K) = \Lambda_\beta({\rm GS}^2(E_1,p))$, where
\begin{eqnarray*}
\Lambda_\beta(K)
= \Big\{ i \in \mathbb{N} : \mbox{there are}~
         z \in K ~\mbox{and}~
         \alpha \neq 0 ~\mbox{satisfying}~ 
         z + \alpha \beta_i \in K
  \Big\}.
\end{eqnarray*}
(Indeed, $K$ is a cylinder. 
If we think that it is difficult to find the cylinder $K$ that
satisfies $E_1 \subset K \subset I^\omega$, consider the fact
that $E_1$ can be covered with countably many cylinders
$C_1, C_2, \ldots$, each of which is included in $I^\omega$.
We simply choose any one of them, namely $C_i$, and call it
$K$.
In this case, we replace $E_1$ with $E_1 \cap C_i$.
With these adjustments, there is no problem in proving this
theorem.)

Then, by Theorem \ref{thm:set9} with $J_\beta = K$, we get
\begin{eqnarray*}
\begin{split}
{\rm GS}(K,p) - p
& =       \Bigg\{ \sum_{i \in \Lambda_\beta(K)}
                  \gamma_i \beta_i \in M_a :
                  \gamma_i \in \mathbb{R} ~\mbox{for all}~
                  i \in \Lambda_\beta(K)
          \Bigg\} \\
& =       \Bigg\{ \sum_{i \in \Lambda_\beta({\rm GS}^2(E_1,p))}
                  \gamma_i \beta_i \in M_a :
                  \gamma_i \in \mathbb{R} ~\mbox{for all}~
                  i \in \Lambda_\beta({\rm GS}^2(E_1,p))
          \Bigg\} \\
& \subset {\rm GS}^2(E_1,p) - p,
\end{split}
\end{eqnarray*}
where the last inclusion is due to the fact that
$\{ \beta_i \}_{i \in \Lambda_\beta({\rm GS}^2(E_1,p))}$ is a
complete orthonormal sequence in the real Hilbert space
${\rm GS}^2(E_1,p) - p$.
Hence, it holds that
$K \subset {\rm GS}(K,p) \subset {\rm GS}^2(E_1,p)$.
Thus, ${\rm vol}(K) = 0$ because ${\rm GS}^2(E_1,p)$ is
degenerate.

If we divide $K$ into countably many translations of degenerate
$\beta$-basic cylinders or $\beta$-basic intervals
$\{ K_i \}_{i \in \mathbb{N}}$ whose diameters are less than
$\delta$, then we get
$\mu_\delta(E_1) \leq \sum\limits_{i=1}^\infty {\rm vol}(K_i)
 = {\rm vol}(K) = 0$ for any $\delta > 0$.
Therefore,
$\mu(E_1) = \lim\limits_{\delta \to 0} \mu_\delta(E_1) = 0$.

On the other hand, we see that
$E_2 = f(E_1) = F_2(E_1) \subset F_2(K)$ and $F_2(K)$ is a
degenerate cylinder.
Hence, we get ${\rm vol}(F_2(K)) = 0$.
And if we do what we just did before, we get $\mu(E_2) = 0$.
Therefore, we conclude that $\mu(E_1) = 0 = \mu(E_2)$.

$(c)$
Now, we consider the case where
$\Lambda({\rm GS}^2(E_1,p)) = \mathbb{N}$, or equivalently,
the case where ${\rm GS}^2(E_1,p) = M_a$
(see Theorem \ref{thm:setsib}).
Let $\delta > 0$ be given.
By the definition of $\mu_\delta$, for any $\varepsilon > 0$,
there exists a $\delta$-covering $\{ K_i \}$ of $E_1$ with
cylinders from $\mathcal{C}_\delta$ such that
\begin{eqnarray} \label{eq:20110912-1}
\sum_i {\rm vol}(K_i) \leq \mu_\delta (E_1) + \varepsilon.
\end{eqnarray}

By the definitions of $\mathcal{B}_\delta$ and
$\mathcal{C}_\delta$, there exist a basic cylinder
$J_i \in \mathcal{B}_\delta$ and a surjective $d_a$-isometry
$f_i : J_i \to K_i$ for each $i$.
Since $F_2 \circ f_i : J_i \to F_2(K_i)$ is a surjective
$d_a$-isometry, ${\rm vol}(F_2(K_i)) = {\rm vol}(J_i)$ for all
$i$ and $\{ F_2(K_i) \} = \{ (F_2 \circ f_i)(J_i) \}$ is a
$\delta$-covering of $E_2$ with cylinders from
$\mathcal{C}_\delta$.
Further, by applying Definition \ref{def:2018-9-23} to the
surjective $d_a$-isometry $f_i : J_i \to K_i$, we get
${\rm vol}(J_i) = {\rm vol}(K_i)$.
Thus, we have
\begin{eqnarray} \label{eq:2018-9-12}
{\rm vol} \big( F_2(K_i) \big) = {\rm vol}(K_i)
\end{eqnarray}
for all $i$.
Therefore, it follows from (\ref{eq:20110912-1}) and
(\ref{eq:2018-9-12}) that
\begin{eqnarray*}
\mu_\delta (E_2)
\leq \sum_i {\rm vol} \big( F_2(K_i) \big)
=    \sum_i {\rm vol}(K_i)
\leq \mu_\delta (E_1) + \varepsilon.
\end{eqnarray*}
Since we can choose a sufficiently small $\varepsilon > 0$, we
conclude that $\mu_\delta (E_2) \leq \mu_\delta (E_1)$.

Conversely, if we exchange the roles of $E_1$ and $E_2$ in the
previous part, then we get
$\mu_\delta (E_1) \leq \mu_\delta (E_2)$.
Hence, we conclude that $\mu_\delta (E_1) = \mu_\delta (E_2)$
for any $\delta > 0$.
Therefore, it follows from (\ref{eq:2018-9-21-2}) that
$\mu(E_1) = \mu(E_2)$.
\hfill$\Box$
\vspace*{5mm}

The following lemmas are the same as the lemmas
\cite[Lemma 3.2, Lemma 3.3, Lemma 3.4]{jung1}.
It is easy to prove that $\mu(I^\omega) \leq 1$.


\begin{lem} \label{lem:3.1.1}
$\mu(I^\omega) \leq 1$.
\end{lem}

\noindent
\emph{Proof}.
We apply an idea from the proof of \cite[Lemma 1]{5}.
For any $\delta > 0$, there exist positive integers $m$ and
$n$ such that $I^n = [0,1]^n$ is covered by $m^n$ isometric
$n$-cubes $C_i$ which are closed in $I^n$ with non-overlapping
interiors and $d_a(J_i) < \delta$ for all
$i \in \{ 1, 2, \ldots, m^n \}$, where $J_i$ is the cylinder
in $I^\omega$ over $C_i$, \textit{i.e.},
$J_i \in \mathcal{B}_{\delta}$ for each $i$.
Then, $\{ J_1, \ldots, J_{m^n} \}$ is a $\delta$-covering of
$I^\omega$ with non-degenerate basic cylinders from
$\mathcal{B}_\delta$ and
\begin{eqnarray*}
\mu_\delta (I^\omega) \leq \sum_{i=1}^{m^n} {\rm vol}(J_i) = 1.
\end{eqnarray*}
Hence, it follows from (\ref{eq:2018-9-21-2}) that
$\mu(I^\omega)
 = \lim\limits_{\delta \to 0} \mu_\delta(I^\omega) \leq 1$.
\hfill$\Box$


\begin{lem} \label{lem:4.3}
Given real numbers $\delta$ and $b$ with
$0 < \delta < \frac{1}{2}$ and $0 \leq b \leq 1$, define
\begin{eqnarray*}
\alpha = \left\{\begin{array}{ll}
                   0
                   & (\mbox{for}~ 0 \leq b < \delta),
                   \vspace*{1mm}\\
                   b - \delta
                   & (\mbox{for}~ \delta \leq b \leq 1)
                \end{array}
         \right.
~~~\mbox{and}~~~
\beta = \left\{\begin{array}{ll}
                  b + \delta
                  & (\mbox{for}~ 0 \leq b \leq 1 - \delta),
                  \vspace*{1mm}\\
                  1
                  & (\mbox{for}~ 1 - \delta < b \leq 1).
               \end{array}
        \right.
\end{eqnarray*}
Then $0 < \beta - \alpha \leq 2\delta$ and
$\alpha \leq b \leq \beta$.
\end{lem}

\noindent
\textit{Proof}.
$(a)$ If $0 \leq b < \delta$, then
$0 \leq b < \delta < 1 - \delta$ and we get $\alpha = 0$ and
$\beta = b + \delta < 2\delta$.
Hence, it follows that $0 < \beta - \alpha < 2\delta$ and
$\alpha \leq b < \beta$.

$(b)$ If $\delta \leq b \leq 1 - \delta$, then
$\alpha = b - \delta$ and $\beta = b + \delta$.
Thus, we have $0 < \beta - \alpha = 2\delta$ and
$\alpha < b < \beta$.

$(c)$ Finally, if $1 - \delta < b \leq 1$, then we see that
$\alpha = b - \delta$ and $\beta = 1$.
So, we have $0 < \beta - \alpha = 1 - b + \delta < 2\delta$
and $\alpha < b \leq \beta$.
\hfill$\Box$
\vspace*{5mm}

For every basic cylinder $J \in \mathcal{B}$, let $\partial J$
denote the boundary of $J$.
In the following lemma, we will prove that $\mu(\partial J) = 0$.


\begin{lem} \label{lem:6}
If $J \in \mathcal{B}$, then $\mu(\partial J) = 0$.
\end{lem}

\noindent
\emph{Proof}.
Similarly as in the proof of Lemma \ref{lem:3.1.1}, there exist
positive integers $m$ and $n$ such that $I^n = [0,1]^n$ is
covered by $m^n$ isometric $n$-cubes $C_i^n$
$(i \in \{ 1, 2, \ldots, m^n \})$ with non-overlapping
interiors, where each $C_i^n$ is a closed subset of $I^n$.

Let $\delta$ and $b$ be any real numbers with
$0 < \delta < \frac{1}{2}$ and $0 \leq b \leq 1$, $\alpha$ and
$\beta$ be defined as in Lemma \ref{lem:4.3},
$P_i = C_i^n \times [\alpha,\beta]$ be an $(n+1)$-dimensional
rectangular parallelepiped, and let $J_i$ denote the cylinder
in $I^\omega$ over $P_i$ with $d_a(J_i) < \delta$ for all
$i \in \{ 1, 2, \ldots, m^n \}$, \textit{i.e.},
$J_i \in \mathcal{B}_\delta$ for each
$i \in \{ 1, 2, \ldots, m^n \}$.

In view of Lemma \ref{lem:4.3}, $\{ J_1, J_2, \ldots, J_{m^n} \}$
is a $\delta$-covering of a hyper-plane $H$ given by
\begin{eqnarray} \label{eq:2018-5-26}
H = \big\{ (x_1, x_2, \ldots) \in I^\omega \,:\, x_{n+1} = b
    \big\}.
\end{eqnarray}
Hence, by the definition of $\mu_\delta$ and using Lemma
\ref{lem:4.3} again, it holds that
\begin{eqnarray*}
\mu_\delta(H) \leq \sum_{i=1}^{m^n} {\rm vol}(J_i)
              =    \sum_{i=1}^{m^n}
                   \frac{1}{m^n} (\beta - \alpha)
              \leq 2 \delta
\end{eqnarray*}
and further we get
\begin{eqnarray} \label{eq:20151007}
\mu(H) = \lim_{\delta \to 0} \mu_\delta(H) = 0.
\end{eqnarray}

Let $J$ be a basic cylinder in $\mathcal{B}$.
In view of Definition \ref{def:basic}, without loss of
generality, we will deal with the basic cylinder of the form
\begin{eqnarray*}
J = \big\{ (x_1, x_2, \ldots) \in I^\omega \,:\,
           p_{1i} \leq x_i \leq p_{2i} ~\,\mbox{for all}\,~
           i \in \mathbb{N}\,
    \big\}
\end{eqnarray*}
only, where there exists a positive integer $n$ such that
$0 \leq p_{1i} < p_{2i} \leq 1$ for each
$i \in \{ 1, \ldots, n \}$ and $p_{1i} = 0$, $p_{2i} = 1$ for
all $i > n$.
Then, there are at most countably many hyper-planes
$H_1, H_2, \ldots$ of the form
\begin{eqnarray*}
H_{2i-1} = \big\{ (x_1, x_2, \ldots) \in I^\omega \,:\,
                  x_i = p_{1i} \big\}
~~~\mbox{and}~~~
H_{2i} = \big\{ (x_1, x_2, \ldots) \in I^\omega \,:\,
                x_i = p_{2i} \big\}
\end{eqnarray*}
satisfying
\begin{eqnarray} \label{eq:20110917-1}
\partial J \subset \bigcup_{k=1}^{\infty} H_k.
\end{eqnarray}

Finally, it follows from (\ref{eq:20151007}) and
(\ref{eq:20110917-1}) that
\begin{eqnarray*}
\mu(\partial J)
\leq \mu\!\left( \bigcup_{k=1}^\infty H_k \right)
\leq \sum_{k=1}^\infty \mu(H_k)
= 0,
\end{eqnarray*}
which completes the proof.
\hfill$\Box$


\section{Efficient coverings}

Let $J$ be a basic cylinder given by Definition \ref{def:basic}
that is $d_a$-isometric to a cylinder $K$ via a surjective
$d_a$-isometry $f : J \to K$.
We now define
\begin{eqnarray} \label{eq:2018-6-5-3}
J^\ast = \prod_{i=1}^\infty J_i^\ast, ~~~\mbox{where}~~~
J_i^\ast = \left\{\begin{array}{ll}
                     \mbox{$[ 0, b^\ast ]$}
                     & (\mbox{for}~ i \in \Lambda_1),
                     \vspace*{2mm}\\
                     \mbox{$[ p_{1i}, p_{1i} + b^\ast ]$}
                     & (\mbox{for}~ i \in \Lambda_2 \cup \Lambda_3),
                     \vspace*{2mm}\\
                     \mbox{$\{ p_{1i} \}$}
                     & (\mbox{for}~ i \in \Lambda_4),
                     \vspace*{2mm}\\
                     \mbox{$[ 0, 1 ]$}
                     & (\mbox{otherwise})
                  \end{array}
           \right.
\end{eqnarray}
and $b^\ast$ is a sufficiently small positive real number in
comparison with each of $p_{2i}$, $p_{2i} - p_{1i}$, and
$1 - p_{1i}$ for all $i \in \Lambda_1$, $i \in \Lambda_2$, and
$i \in \Lambda_3$, respectively.
We then note that $J^\ast \subset J$.
Taking Remark \ref{rem:2.2} $(ii)$ and $(iii)$ into account,
we can cover the cylinder $K = f(J)$ with a finite number of
translations of $f(J^\ast)$ as efficiently as we wish by
choosing the $b^\ast$ sufficiently small
(see the illustration and Lemma \ref{lem:7} below).

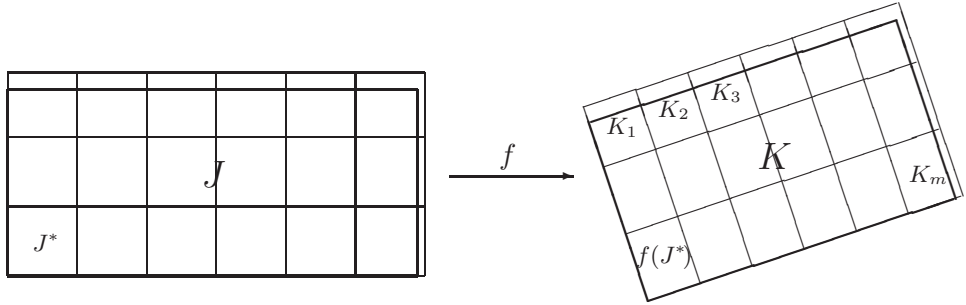
\begin{center}
\setlength{\unitlength}{0.33mm}
\begin{picture}(382,155)(0,-38)

\thicklines
\put(0,0){\line(1,0){165}}
\put(0,75){\line(1,0){165}}
\put(0,0){\line(0,1){75}}
\put(165,0){\line(0,1){75}}

\thinlines
\put(0,0){\line(1,0){168}}
\put(0,28){\line(1,0){168}}
\put(0,56){\line(1,0){168}}
\put(0,82){\line(1,0){168}}
\put(0,0){\line(0,1){82}}
\put(28,0){\line(0,1){82}}
\put(56,0){\line(0,1){82}}
\put(84,0){\line(0,1){82}}
\put(112,0){\line(0,1){82}}
\put(140,0){\line(0,1){82}}
\put(168,0){\line(0,1){82}}

\put(78,36){\Large\bf $J$}
\put(10,10){\footnotesize $J^\ast$}

\thicklines
\put(258,-10){\line(3,1){123}}
\put(235,62){\line(3,1){123}}
\put(258,-10){\line(-1,3){24}}
\put(382,31){\line(-1,3){24}}

\thinlines
\put(258,-10){\line(3,1){126}}
\put(232,68){\line(3,1){126}}
\put(258,-10){\line(-1,3){26}}
\put(385,32){\line(-1,3){26}}
\put(279,-3){\line(-1,3){26}}
\put(300,4){\line(-1,3){26}}
\put(321,11){\line(-1,3){26}}
\put(342,18){\line(-1,3){26}}
\put(363,25){\line(-1,3){26}}
\put(249,17){\line(3,1){126}}
\put(240,44){\line(3,1){126}}

\put(302,43){\Large\bf $K$}
\put(253,5){\footnotesize $f(J^\ast\!)$}
\put(241,57){\footnotesize $K_1$}
\put(262,64){\footnotesize $K_2$}
\put(283,71){\footnotesize $K_3$}
\put(363,38){\footnotesize $K_m$}

\thinlines
\put(178,40){\vector(1,0){50}}
\put(198,45){$f$}

\put(53,-33){\text{A finite number of translations of $f(J^\ast)$ cover $K = f(J)$}}

\end{picture}
\end{center}

$J^\ast$ is a basic cylinder and the restriction
$f |_{J^\ast} : J^\ast \to f(J^\ast)$ is a surjective
$d_a$-isometry.
Thus, $f(J^\ast)$ is a cylinder, \textit{i.e.},
$f(J^\ast) \in \mathcal{C}$ (see Remark \ref{rem:3.1}).
Applying this argument, Remark \ref{rem:2.2} $(ii)$, $(iii)$
and Lemma \ref{lem:6}, we obtain the following lemma that is
an improved version of \cite[Lemma 4.1]{jung1}.
This new version includes the degenerate case.


\begin{lem} \label{lem:7}
Let $\delta > 0$ and $\varepsilon > 0$ be given.
If $K$ is a cylinder from $\mathcal{C}_\delta$, then there
exist a finite number of translations $K_1, K_2, \ldots, K_m$
of some cylinder in $\mathcal{C}_\delta$
$($\/for example, $f(J^\ast)$ above and see the corresponding
illustration above\/$)$ such that
\begin{itemize}
\item[$(i)$]   $K_i \cap K_j$ $(i \neq j)$ is included in the
               union of at most countably many $(d_a$-isometric
               images of\/$)$ hyper-planes of the form
               $\big\{
               (x_1, x_2, \ldots) \in I^\omega \,:\, x_{\ell} = b
                \big\}$ for some $\ell \in \mathbb{N}$ and
               $0 \leq b \leq 1$, which have $\mu$-measure $0$;
\item[$(ii)$]  $\{ K_1, K_2, \ldots, K_m \}$ is a covering of
               $K$, \textit{i.e.},
               $K \subset \bigcup\limits_{i=1}^m K_i$;
\item[$(iii)$] $\sum\limits_{i=1}^m {\rm vol}(K_i)
                \leq {\rm vol}(K) + \varepsilon$.
\end{itemize}
\end{lem}

\noindent
\emph{Proof}.
We can choose a $J \in \mathcal{B}_\delta$ and a surjective
$d_a$-isometry $f : J \to K$, where $J$ is a basic cylinder of
the form given in Definition \ref{def:basic}.
We now define a basic cylinder $J^\ast$ by
(\ref{eq:2018-6-5-3}) such that
$J^\ast \subset J$.
Then $J$ can be covered with at most $m$ translations of the
basic cylinder $J^\ast$, where we set
\begin{eqnarray*}
m := \prod_{j \in \Lambda_1}
     \bigg( \bigg[ \frac{p_{2j}}{b^\ast} \bigg] + 1 \bigg)
     \times
     \prod_{j \in \Lambda_2}
     \bigg( \bigg[ \frac{p_{2j} - p_{1j}}{b^\ast} \bigg] + 1
     \bigg) \times
     \prod_{j \in \Lambda_3}
     \bigg( \bigg[ \frac{1 - p_{1j}}{b^\ast} \bigg] + 1 \bigg)
\end{eqnarray*}
and where $[x]$ denotes the largest integer not exceeding the
real number $x$.
This fact, together with Remark \ref{rem:2.2} $(ii)$ and
$(iii)$, implies that the cylinder $K = f(J)$ can be covered
with at most $m$ translations of the cylinder $f(J^\ast)$
which are denoted by $K_1, K_2, \ldots, K_m$
(see the previous illustration).

Moreover, in view of Remark \ref{rem:2.2} $(iii)$, we have
\begin{eqnarray*}
\begin{split}
{\rm vol}(K_i) = {\rm vol}(f(J^\ast)) = {\rm vol}(J^\ast)
= \left\{\begin{array}{ll}
            \prod\limits_{j \in \Lambda_1 \cup \Lambda_2 \cup \Lambda_3}
            b^\ast
            & (\mbox{for non-degenerate $K$}), \vspace*{2mm}\\
            0
            & (\mbox{for degenerate $K$})
         \end{array}
  \right.
\end{split}
\end{eqnarray*}
for each $i \in \{ 1, 2, \ldots, m \}$ and
\begin{eqnarray*}
\begin{split}
{\rm vol}(K)
& = {\rm vol}(J) \\
& = \left\{\begin{array}{ll}
              \prod\limits_{j \in \Lambda_1} p_{2j} \times
              \prod\limits_{j \in \Lambda_2} (p_{2j} - p_{1j})
              \times
              \prod\limits_{j \in \Lambda_3} (1 - p_{1j})
              & (\mbox{for non-degenerate $K$}), \vspace*{2mm}\\
              0
              & (\mbox{for degenerate $K$}).
           \end{array}
    \right.
\end{split}
\end{eqnarray*}

Thus, when $K$ is a non-degenerate cylinder, we have
\begin{eqnarray*}
\begin{split}
\sum_{i=1}^m {\rm vol}(K_i)
=    &\, \prod_{j \in \Lambda_1}
         \bigg( \bigg[ \frac{p_{2j}}{b^\ast} \bigg] + 1 \bigg)
         \times
         \prod_{j \in \Lambda_2}
         \bigg( \bigg[ \frac{p_{2j} - p_{1j}}{b^\ast} \bigg] + 1
         \bigg) \times \\
     &\, \times
         \prod_{j \in \Lambda_3}
         \bigg( \bigg[ \frac{1 - p_{1j}}{b^\ast} \bigg] + 1
         \bigg) \times {\rm vol}(K_1) \\
\leq &\, \prod_{j \in \Lambda_1}
         \bigg( \frac{p_{2j}}{b^\ast} + 1 \bigg) \times
         \prod_{j \in \Lambda_2}
         \bigg( \frac{p_{2j} - p_{1j}}{b^\ast} + 1 \bigg)
         \times \\
     &\, \times
         \prod_{j \in \Lambda_3}
         \bigg( \frac{1 - p_{1j}}{b^\ast} + 1 \bigg) \times
         \prod_{j \in \Lambda_1 \cup \Lambda_2 \cup \Lambda_3}
         b^\ast \\
=    &\, \prod_{j \in \Lambda_1} ( p_{2j} + b^\ast ) \times
         \prod_{j \in \Lambda_2} ( p_{2j} - p_{1j} + b^\ast )
         \times
         \prod_{j \in \Lambda_3} ( 1 - p_{1j} + b^\ast ) \\
=    &\, \prod_{j \in \Lambda_1} p_{2j} \times
         \prod_{j \in \Lambda_2} (p_{2j} - p_{1j}) \times
         \prod_{j \in \Lambda_3} (1 - p_{1j}) + O(b^\ast) \\
=    &\; {\rm vol}(K) + O(b^\ast)
\end{split}
\end{eqnarray*}
and we choose a sufficiently small $b^\ast$ such that the term
$O(b^\ast)$ becomes less than $\varepsilon$.
When $K$ is degenerate, we have ${\rm vol}(K) = 0$ and
${\rm vol}(K_i) = 0$ for any $i \in \{ 1, 2, \ldots, m \}$.
Hence, our assertion $(iii)$ holds true.
\hfill$\Box$
\vspace{5mm}

Using Lemmas \ref{lem:3.1.1} and \ref{lem:7}, we will prove
that $\mu(I^\omega) = 1$.
The following theorem is equivalent to \cite[Theorem 4.2]{jung1},
but the proof of this theorem is much more concise than that
of \cite[Theorem 4.2]{jung1}.
Hence, we will introduce the proof.


\begin{thm} \label{thm:8}
$\mu(I^\omega) = 1$.
\end{thm}

\noindent
\emph{Proof}.
Given a $\delta > 0$ and an $\varepsilon > 0$, let $\{ K_i \}$
be a $\delta$-covering of $I^\omega$ with cylinders from
$\mathcal{C}_\delta$ such that
\begin{eqnarray} \label{eq:20110919-1}
\sum_{i=1}^\infty {\rm vol}(K_i)
\leq \mu_\delta(I^\omega) + \frac{\varepsilon}{2}.
\end{eqnarray}
In view of Lemma \ref{lem:7}, for each $K_i$, there exist
translations $K_{i1}, K_{i2}, \ldots, K_{im_i}$ of some
cylinder in $\mathcal{C}_\delta$ such that
\begin{itemize}
\item[$(i)$]   $K_{ij} \cap K_{i\ell}$ $(j \neq \ell)$ is
               included in the union of at most countably many
               ($d_a$-isometric images of) hyper-planes of the
               form
               $\big\{
               (x_1, x_2, \ldots) \in I^\omega \,:\, x_{\ell} = b
                \big\}$ for some $\ell \in \mathbb{N}$ and
               $0 \leq b \leq 1$, whose $\mu$-measures are $0$;
\item[$(ii)$]  $K_i \subset \bigcup\limits_{j=1}^{m_i} K_{ij}$;
\item[$(iii)$] $\sum\limits_{j=1}^{m_i}
                {\rm vol}(K_{ij}) \leq
                {\rm vol}(K_i) + \frac{\varepsilon}{2^{i+1}}$.
\end{itemize}

We notice that each $K_{ij}$ is a cylinder from
$\mathcal{C}_\delta$.
If we replace the covering $\{ K_i \}$ with a new
$\delta$-covering $\{ K_{ij} \}$, then it follows from
$(iii)$ that
\begin{eqnarray*}
\sum_{i=1}^\infty \sum_{j=1}^{m_i} {\rm vol}(K_{ij})
\leq \sum_{i=1}^\infty {\rm vol}(K_i) + \frac{\varepsilon}{2}
\end{eqnarray*}
and it follows from this inequality and (\ref{eq:20110919-1})
that
\begin{eqnarray*}
1 = {\rm vol}(I^\omega) \leq
\sum_{i=1}^\infty \sum_{j=1}^{m_i} {\rm vol}(K_{ij})
\leq \sum_{i=1}^\infty {\rm vol}(K_i) + \frac{\varepsilon}{2}
\leq \mu_\delta(I^\omega) + \varepsilon.
\end{eqnarray*}

If we take a sufficiently small value of $\varepsilon > 0$,
then we have $\mu_\delta(I^\omega) \geq 1$ and
$\mu(I^\omega) = \lim\limits_{\delta \to 0} \mu_\delta(I^\omega)
 \geq 1$.
On the other hand, in view of Lemma \ref{lem:3.1.1}, we have
$\mu(I^\omega) \leq 1$.
Hence, we conclude that $\mu(I^\omega) = 1$.
\hfill$\Box$


\section{Ulam's conjecture on invariance of measure}

According to \cite[Theorems 16 and 19]{6}, all Borel sets in
$M_a$ are $\mu$-measurable.
Moreover, each Borel subset of $I^\omega$ is also a Borel
subset of $M_a$, \textit{i.e.}, each Borel subset of $I^\omega$
is $\mu$-measurable.

In view of Theorems \ref{thm:5} and \ref{thm:8}, the measure
$\mu$ is $d_a$-invariant with $\mu(I^\omega) = 1$.
The proof of the following lemma is the same as that of
\cite[Lemma 5.1]{jung1}. Hence, we omit the proof.


\begin{lem} \label{lem:9}
The measure $\mu$ coincides with the standard product
probability measure $\pi$ on the Borel subsets of $I^\omega$.
\end{lem}

According to Theorem \ref{thm:5}, the measure $\mu$ is
$d_a$-invariant.
Using Lemma \ref{lem:9}, we obtain our main result:


\begin{thm} \label{thm:10}
For any sequence $a = \{ a_i \}$ of positive real numbers
satisfying $(\ref{eq:20110922-1})$, the standard product
probability measure $\pi$ on $I^\omega$ is $d_a$-invariant.
More precisely, if $E_1$ and $E_2$ are Borel subsets of
$I^\omega$ that are $d_a$-isometric to each other, then
$\pi(E_1) = \pi(E_2)$.
\end{thm}
\vspace*{5mm}



\end{document}